\documentclass[review]{elsarticle}

\usepackage{lineno,hyperref}

\usepackage{array,color}
\usepackage{amsbsy}
\usepackage{amsmath}
\usepackage{amssymb}
\usepackage{amsfonts}

\usepackage{epsf}
\usepackage{epsfig}
\usepackage{graphicx}
\usepackage{fancybox,fancyhdr,times,subfigure}
\usepackage{epsfig, color, graphicx, rotate, colordvi}
\usepackage{a4wide}
\usepackage{natbib}
\usepackage{mathrsfs}
\usepackage{bm}
\usepackage{multirow}
\usepackage{pgfplots}

\newcommand{\bbB}{{\bf B}}
\newcommand{\bbZ}{{\bf Z}}

\newcommand{\bbX}{{\bf X}}
\newcommand{\bbS}{{\bf S}}

\newcommand{\bbI}{{\bf I}}

\newcommand{\bbF}{{\bf F}}

\newcommand{\bqn}{\begin{eqnarray*}}
\newcommand{\eqn}{\end{eqnarray*}}

\newcommand{\bqa}{\begin{eqnarray}}
\newcommand{\eqa}{\end{eqnarray}}

\newcommand{\al}{\alpha}

\newtheorem{thm}{Theorem}
\newtheorem{lem}[thm]{Lemma}
\newdefinition{rmk}{Remark}
\newproof{pf}{Proof}
\newproof{pot}{Proof of Theorem \ref{thm2}}

\modulolinenumbers[5]

\journal{}

\begin{document}

\begin{frontmatter}

\title{Modified Pillai's trace statistics for two high-dimensional sample covariance matrices}



\author{Qiuyan Zhang}
\ead{zhangqy919@nenu.edu.cn}

\author{Jiang Hu*}
\ead{huj156@nenu.edu.cn}

\author{Zhidong Bai}
\ead{baizd@nenu.edu.cn}

\address{Key Laboratory for Applied Statistics of the Ministry of Education,
School of Mathematics and Statistics,
Northeast Normal University, China.}

\cortext[mycorrespondingauthor]{Corresponding author.}

\begin{abstract}
The goal of this study was to test the equality of two covariance matrices by using modified Pillai's trace statistics under a high-dimensional framework, i.e., the dimension and sample sizes go to infinity proportionally. In this paper, we introduce two modified Pillai's trace statistics and obtain their asymptotic distributions under the null hypothesis. The benefits of the proposed statistics include the following: (1) the sample size can be smaller than the dimensions; (2) the limiting distributions of the proposed statistics  are universal; and  (3) we do not restrict the structure of the population covariance matrices. 
The theoretical results are established under mild and practical assumptions, and their  properties  are demonstrated numerically by simulations and a real data analysis.

\end{abstract}


\end{frontmatter}

\linenumbers

\section{Introduction}
High-dimensional data are common in modern scientific domains, such as finance and wireless communication. Hence, the testing of covariance matrices under high-dimensional settings constitutes an important issue in these areas.
The following three main tests have been investigated widely by
statisticians for one-sample tests: (i) sphericity test, (ii) identity matrix test, and (iii) diagonal matrix test.
Ledoit and Wolf \cite{Ledoit2002} investigated the properties of the sphericity and identity matrix tests when the sample size and the dimension converge to infinity proportionally, and Birke and Dette \cite{BirkeDette2005} generalized Ledoit and Wolf's \cite{Ledoit2002} conclusion to the case where the sample size and dimension are not of the same order. Srivastava \cite{Srivastava2005} proved the asymptotic null and alternative distributions of the testing statistics for normally distributed data. Furthermore, Chen et al. \cite{Chen2010} proposed a nonparametric method and reported that its data could come from any distribution with a specified data structure.
Cai and Ma \cite{CaiMa2013} developed an identity matrix test procedure based on minimax analysis and showed that the power of their test uniformly dominates the power of the corrected likelihood ratio test by Bai et al. \cite{Bai2009} over the entire asymptotic regime. Under the alternative hypothesis, Chen and Jiang \cite{ChenJiang2018} demonstrated the central limit theorem (CLT) of the likelihood ratio test (LRT) statistic. Schott \cite{Schott2006}, Fisher et al. \cite{FisherSunGallagher2010}, Srivastava et al. \cite{Srivastava2011}, Qiu and Chen \cite{QiuChen2012} and Wu and Li \cite{WuLi2015} also analyzed this issue in depth.

Moreover, testing procedures for the equivalence of high-dimensional two-sample covariance matrices are also frequently considered.
Regarding the hypothesis test problem,
\begin{equation}\label{hypo}
H_0: \Sigma_1=\Sigma_2\quad   \mathrm{v.s.} \quad  H_1:\Sigma_1\neq\Sigma_2,
\end{equation}
where $\Sigma_1$ and $\Sigma_2$ are two population covariance matrices. Shott \cite{Schott2007} proposed a statistic based on the idea of an unbiased estimation of the squared Frobenius norm of $\Sigma_1-\Sigma_2$, and showed its asymptotic distribution under the condition that the sample sizes and the dimension converge to infinity proportionally.
A similar idea was adopted by  Li and Chen \cite{Li2012}  and Gao et al. \cite{Gao2013}. 
In addition, Srivastava \cite{Srivastava2008} considered the lower bound of this Frobenius norm, 
and Zhang et al. \cite{Zhang2017} generalized Li and Chen's statistic to multiple samples. Srivastava and Yanagihara \cite{Srivastava2009} considered the distance measure $\mathbf{tr}\Sigma_1^2/(\mathbf{tr}\Sigma_1)^2-\mathbf{tr}\Sigma_2^2/(\mathbf{tr}\Sigma_2)^2$ and proposed a test based on a consistent estimation of this distance.
Moreover, Cai et al. \cite{Cai2013} developed an estimator to find the maximum difference between entries in two-sample covariance matrices. Bai et al. \cite{Bai2009}, Zhang et al. \cite{ZhangHuBai2017} and Jiang et al. \cite{Jiang2012} presented the asymptotic distribution of the correctional LRT under high-dimensional assumptions.
 Later, 
Zheng et al. \cite{ZhengBaiYao2015} extended the results of Bai et al. \cite{Bai2009} to general populations with unknown means. 

The goal of this study is to test the hypothesis \eqref{hypo}. Assume that our samples $\{z_{1}^{(l)},z_{2}^{(l)},\dots,z_{n_l}^{(l)},l=1,2\}$ are drawn independently from populations $z^{(l)}$ with mean $\mu_{l}$ and covariance matrices $\Sigma_{l}$. 
We denote $\bbS_l^{z}:=\frac{1}{n_l-1}\sum_{j=1}^{n_l}(z_j^{(l)}-\bar{z}^{(l)})
(z_j^{(l)}-\bar{z}^{(l)})'$, where $\bar{z}^{(l)}=\frac{1}{n_l}\sum_{j=1}^{n_l}z_{j}^{(l)}$.
For the test problem \eqref{hypo}, we choose Pillai's classic trace statistic {$$\mathbf{tr}\bbS_1^{z}(\bbS_1^{z}+\frac{n_2}{n_1}\bbS_2^{z})^{-1},$$} which was first proposed by Pillai \cite{Pillai1954}. For convenience, we subsequently denote $\bbB_{n}^{z}(\bbZ_{1},\bbZ_{2}):=\bbS_1^{z}(\bbS_1^{z}+\frac{n_2}{n_1}\bbS_2^{z})^{-1}$, which is called the Beta matrix and was proposed by Bai et al. \cite{Bai2015}.
From this definition, we note that to guarantee the reversibility of $\bbS_1^{z}+\frac{n_2}{n_1}\bbS_2^{z}$, $p$ must be smaller than $n_1+n_2$.
 The {asymptotic} property of Pillai's statistic has been obtained by using the moment method under the condition that the sample sizes diverge but the dimension is fixed.
Motivated by Bai et al. \cite{Bai2015}, in this paper, we modify Pillai's trace statistic by removing the one and zero eigenvalues of $\bbB_{n}^{z}$, that is,
$$\mathcal{L}=\sum_{\lambda_k^{\bbB_{n}^{z}(\bbZ_{1},\bbZ_{2})}\neq\{0,1\}}\lambda_k^{\bbB_{n}^{z}(\bbZ_{1},\bbZ_{2})},$$ {where $\lambda_k^{\bbB_{n}^{z}(\bbZ_{1},\bbZ_{2})}$ are eigenvalues of $\bbB_{n}^{z}(\bbZ_{1},\bbZ_{2})$}.

In a similar fashion, we modify another of {Pillai's} trace statistics, $$\widetilde{\mathbf{L}}=c_{n_1}\mathbf{tr}(\frac{1}{c_{n_1}}\bbB_{n}^{z}(\bbZ_{1},\bbZ_{2})-\bbI_{p})^{2}
+c_{n_2}\mathbf{tr}(\frac{1}{c_{n_2}}\bbB_{n}^{z}(\bbZ_{2},\bbZ_{1})-\bbI_{p})^{2}$$ {and transform $\widetilde{\mathbf{L}}$ to} $$\widetilde{\mathcal{L}}=\sum_{\lambda_k^{\bbB_{n}^{z}(\bbZ_{1},\bbZ_{2})}\neq\{0,1\}}\sum_{
\lambda_{k'}^{\bbB_{n}^{z}(\bbZ_{2},\bbZ_{1})}\neq\{0,1\}}[c_{n_1}(\frac{1}{c_{n_1}}
\lambda_{k}^{\bbB_{n}^{z}(\bbZ_{1},\bbZ_{2})}-1)^{2}
+c_{n_2}(\frac{1}{c_{n_2}}{\lambda_{k'}^{\bbB_{n}^{z}(\bbZ_{2},\bbZ_{1})}}-1)^{2}],$$ where
$c_{n_1}=\frac{n_1}{n_1+n_2},$ $c_{n_2}=\frac{n_2}{n_1+n_2}$,  $\bbB_{n}^{z}(\bbZ_{2},\bbZ_{1})=\bbS_2^{z}(\bbS_2^{z}+\frac{n_1}{n_2}\bbS_1^{z})^{-1}$ and $\lambda_{k'}^{\bbB_{n}^{z}(\bbZ_{2},\bbZ_{1})}$ are eigenvalues of $\bbB_{n}^{z}(\bbZ_{2},\bbZ_{1}).$  
In the next section, we will show the CLTs of  $\mathcal{L}$ and  $\widetilde{\mathcal{L}}$ under a high-dimensional setting under the null hypothesis.

The main technical tool employed in this paper is random matrix theory (RMT), which is a powerful method  when the dimension $p$ is large.
Marchenko and Pastur \cite{MarchenkoPastur1967}  determined the limiting spectral distribution of a large-dimensional sample covariance matrix.
Bai and Silverstein \cite{BaiSilverstein2004} proposed a CLT for the linear spectral statistics (LSS) of large-dimensional sample covariance matrices that highlights this issue.
Zheng \cite{Zheng2012} considers a CLT for the LSS of a large-dimensional F matrix, which is used to fulfill the two-sample test. However, the drawback of their method is that the dimension $p$ must be smaller than $\max\{n_1, n_2\}$.
Bai and Yao \cite{BaiYao2008} focused on the spiked model, which was first proposed by Johnstone \cite{Johnstone2001}, and established a limit theorem of extreme sample eigenvalues.
Similar works include Baik and Silverstein \cite{BaikSilverstein2006}, Paul \cite{Paul2007}, Bai et al. \cite{BaiJiangYaoEtAl2013} and Passemier et al. \cite{PassemierMckayChen2015}. Recently, Bai et al. \cite{Bai2015} proved the CLT for the LSS of the Beta matrix using the asymptotically normally distributed property of the sum of the martingale difference sequence and extended the dimension to a high-dimensional situation.

 One should notice that  the main technical tool used here is  Cauchy's residue theorem---the same technique utilized in Zhang et al. \cite{ZhangHuBai2017}; however, the difference is that the integrands in the current paper are linear functions, whereas the integrands for the LRT statistics proposed in Zhang et al. \cite{ZhangHuBai2017} are logarithmic functions. 
 Moreover, these linear functions can be implemented more rapidly and in a less source-consuming way than the abovementioned logarithmic functions, which have greater computational complexity. In addition,  it is clear that when $\frac{p}{n_1}$ or $\frac{p}{n_2}$ tend to $1$, the Beta matrix will have eigenvalues tend to $0$ or $1$ which causes logarithm function tend to infinity. Therefore, the variance of the statistic would tend to infinity which surely makes the test less powerful. However, due to the linearity of the integrand functions, Pillai's trace statistics can be utilized in a space where $\frac{p}{n_1}=1$ or $\frac{p}{n_2}=1$. Therefore, Bai et al. \cite{Bai2015}, Zhang et al. \cite{ZhangHuBai2017} and the current  paper can be viewed as a series of works aimed at improving the classic test statistics of two-sample covariance matrices (see (1.2) in Bai et al. \cite{Bai2015})  from a low-dimensional framework to a high-dimensional framework.
In Section 3, we compare the test statistics proposed in this paper with Li and Chen's \cite{Li2012} statistic,  Cai et al.'s \cite{Cai2013} statistic and Zhang et al.'s \cite{ZhangHuBai2017} statistics through simulations.

The remainder of this paper is organized as follows. Section 2 presents the main conclusions related to the proposed statistics. The results of the simulations, including the comparison with Li and Chen's \cite{Li2012} statistic, Cai et al.'s \cite{Cai2013} statistic and Zhang et al.'s \cite{ZhangHuBai2017} statistics, are presented in Section 3. Section 4 includes an analysis using real Standard and Poor's 500 index data. The proof is presented in the Appendix.

\section{Asymptotically normal property}
In this section, we present the main  results of this paper. 
In the sequel,  we assume the samples satisfy     
\begin{align}\label{linear model}
	z_i^{(l)}=\Sigma_{l}^{1/2} x_{i}^{(l)}+\mu_{l},
\end{align}
where $x_i^{(l)}=(x_{i1}^{(l)},\dots,x_{ip}^{(l)})'$ and $l=1, 2$. Because the trace statistic is invariant under the null hypothesis $\Sigma_1=\Sigma_2$, we have 
{\begin{equation*}
\mathbf{tr}\bbS_1^{z}(\bbS_1^{z}+\frac{n_2}{n_1}\bbS_2^{z})^{-1}
=\mathbf{tr}\bbS_1^{x}(\bbS_1^{x}+\frac{n_2}{n_1}\bbS_2^{x})^{-1},
\end{equation*}}
where $\bbS_l^{x}:=\frac{1}{n_l-1}\sum_{j=1}^{n_l}(x_j^{(l)}-\bar{x}^{(l)})
(x_j^{(l)}-\bar{x}^{(l)})'$ and $\bar{x}^{(l)}=\frac{1}{n_l}\sum_{j=1}^{n_l}x_{j}^{(l)}$, $l=1,2.$ Then, we redefine the Beta matrix as $\bbB_{n}^{x}(\bbX_{1},\bbX_{2})=\bbS_1^{x}(\bbS_1^{x}+\frac{n_2}{n_1}\bbS_2^{x})^{-1}$, $\bbB_{n}^{x}(\bbX_{2},\bbX_{1})=\bbS_2^{x}(\bbS_2^{x}+\frac{n_1}{n_2}\bbS_1^{x})^{-1}$ and transform our statistics into the following forms
$$\mathcal{L}=\sum_{\lambda_k^{\bbB_{n}^{x}(\bbX_{1},\bbX_{2})}\neq\{0,1\}}\lambda_k^{\bbB_{n}^{x}(\bbX_{1},\bbX_{2})},$$
$$\widetilde{\mathcal{L}}=\sum_{\lambda_k^{\bbB_{n}^{x}(\bbX_{1},\bbX_{2})}\neq\{0,1\}}\sum_{
\lambda_{k'}^{\bbB_{n}^{x}(\bbX_{2},\bbX_{1})}\neq\{0,1\}}
[c_{n_1}(\frac{1}{c_{n_1}}
\lambda_k^{\bbB_{n}^{x}(\bbX_{1},\bbX_{2})}-1)^{2}
+c_{n_2}(\frac{1}{c_{n_2}}(1-\lambda_{k'}^{\bbB_{n}^{x}(\bbX_{1},\bbX_{2})})-1)^{2}].$$
Under the following mild assumptions
\begin{enumerate}[(1)]
\item $\{x_{ij}^{(l)},~i=1,\dots,p,~j=1,\dots,n_l\}$ are independent and identically distributed real random variables;
\item As $\min\{p,n_1,n_2\}\rightarrow \infty$, $y_{n_1}:=\frac p {n_1}\rightarrow y_1\in(0,+\infty)$, $y_{n_2}:=\frac p {n_2}\rightarrow y_2\in(0,+\infty)$ and $\alpha_n:=n_2/n_1\rightarrow \alpha >0$;
\item As $\min\{p,n_1,n_2\}\rightarrow \infty$, $h_{n}:=\sqrt{y_{n_1}+y_{n_2}-y_{n_1}y_{n_2}}\to \sqrt{y_{1}+y_{2}-y_{1}y_{2}}>0$;
\item  {$\textbf{E}x_{ij}^{(l)}=0$, $\textbf{E}(x_{ij}^{(l)})^2=1$}, $\Delta_{1}=:\textbf{E}(x_{ij}^{(1)})^4-3<\infty$ and $\Delta_2:=\textbf{E}(x_{ij}^{(2)})^4-3<\infty$,
    \end{enumerate}
we draw the following conclusion about the modified Pillai's statistic $T_1$.

\begin{figure}[h]
	\centering
	\subfigure{
		\label{fig:subfig:a} 
		\includegraphics[width=4.8cm,height=4.8cm]{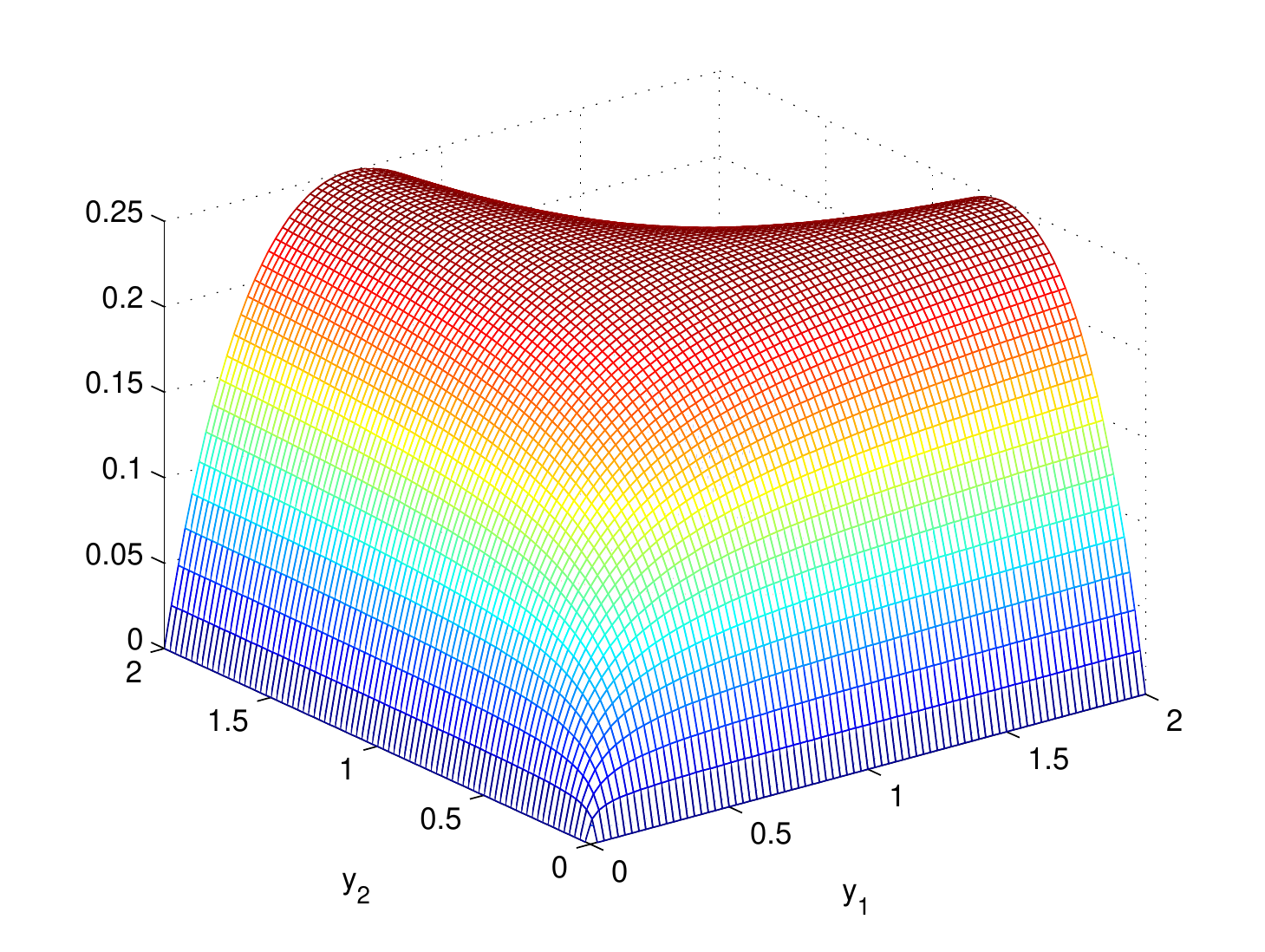}}
	\subfigure{
		\label{fig:subfig:b} 
		\includegraphics[width=4.8cm,height=4.8cm]{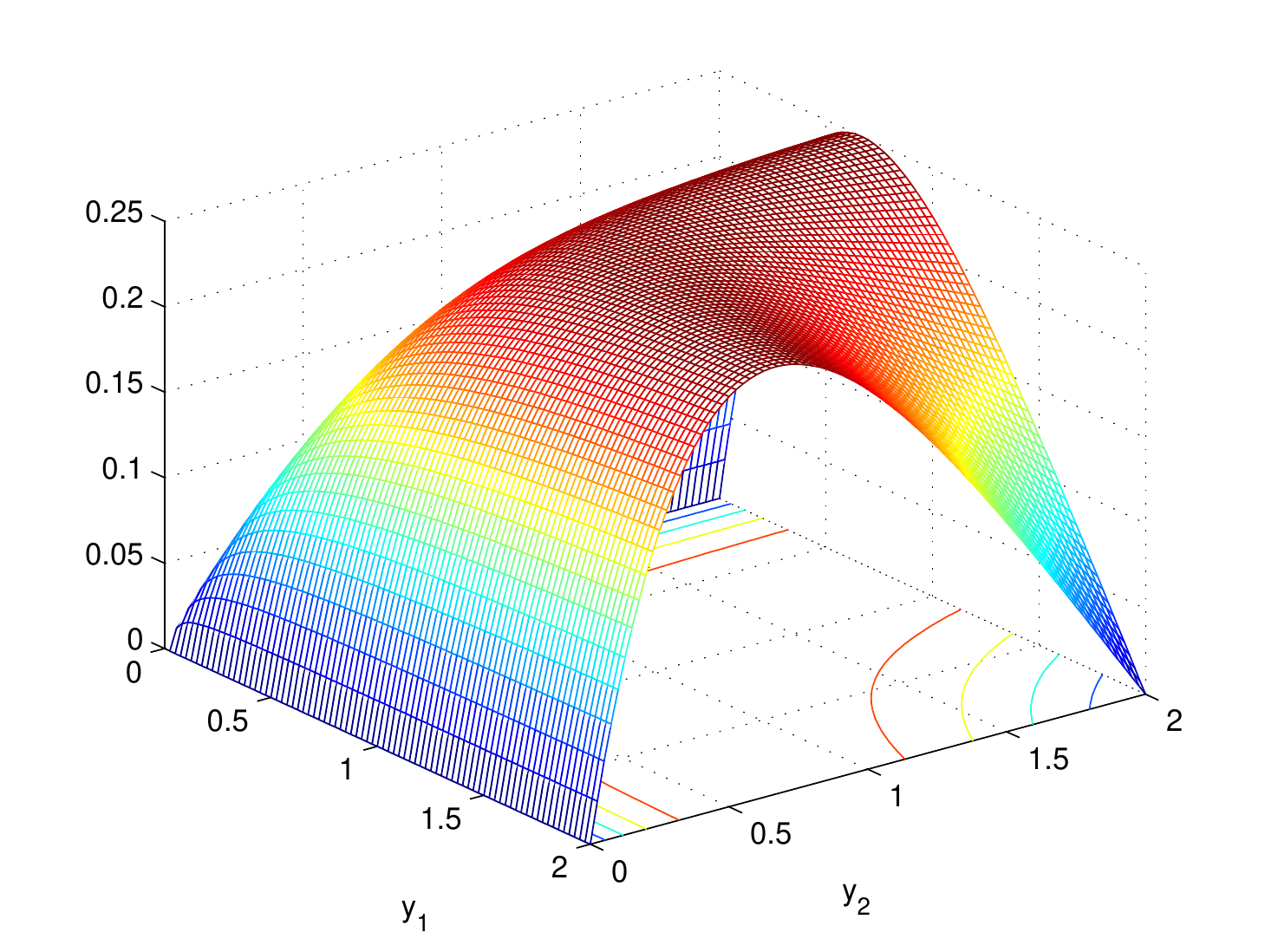}}
	\subfigure{
		\label{fig:subfig:b} 
		\includegraphics[width=4.8cm,height=4.8cm]{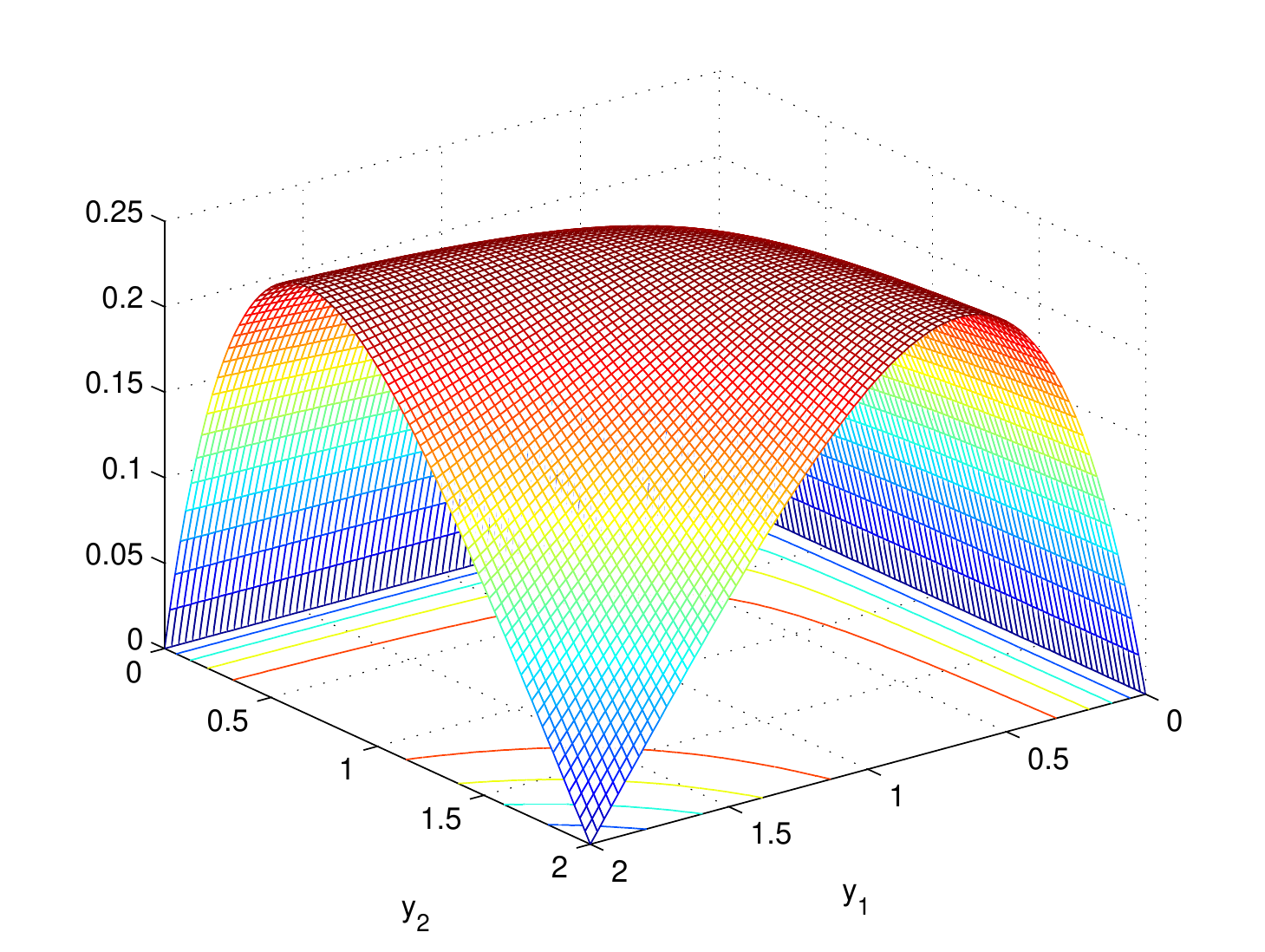}}
	\caption{Three-dimensional shaded figure of $\tilde{\mu}_{n}$ from three angles when $y_1\in(0,2)$ and $y_2\in(0,2)$. }\label{3Dfigure1}
\end{figure}

\begin{figure}[h]
	\centering
	\subfigure{
		\label{fig:subfig:a} 
		\includegraphics[width=4.8cm,height=4.8cm]{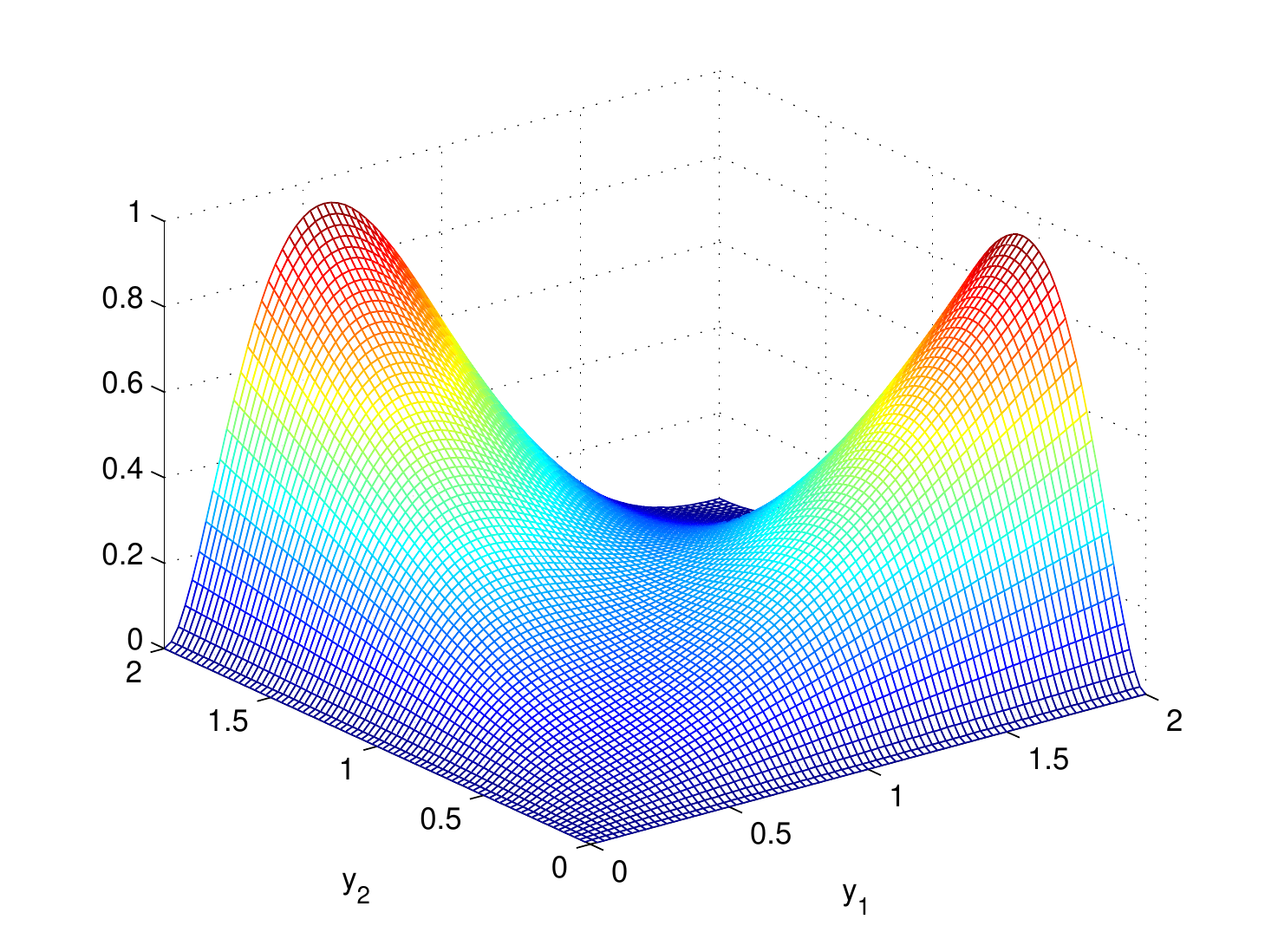}}
	\subfigure{
		\label{fig:subfig:b} 
		\includegraphics[width=4.8cm,height=4.8cm]{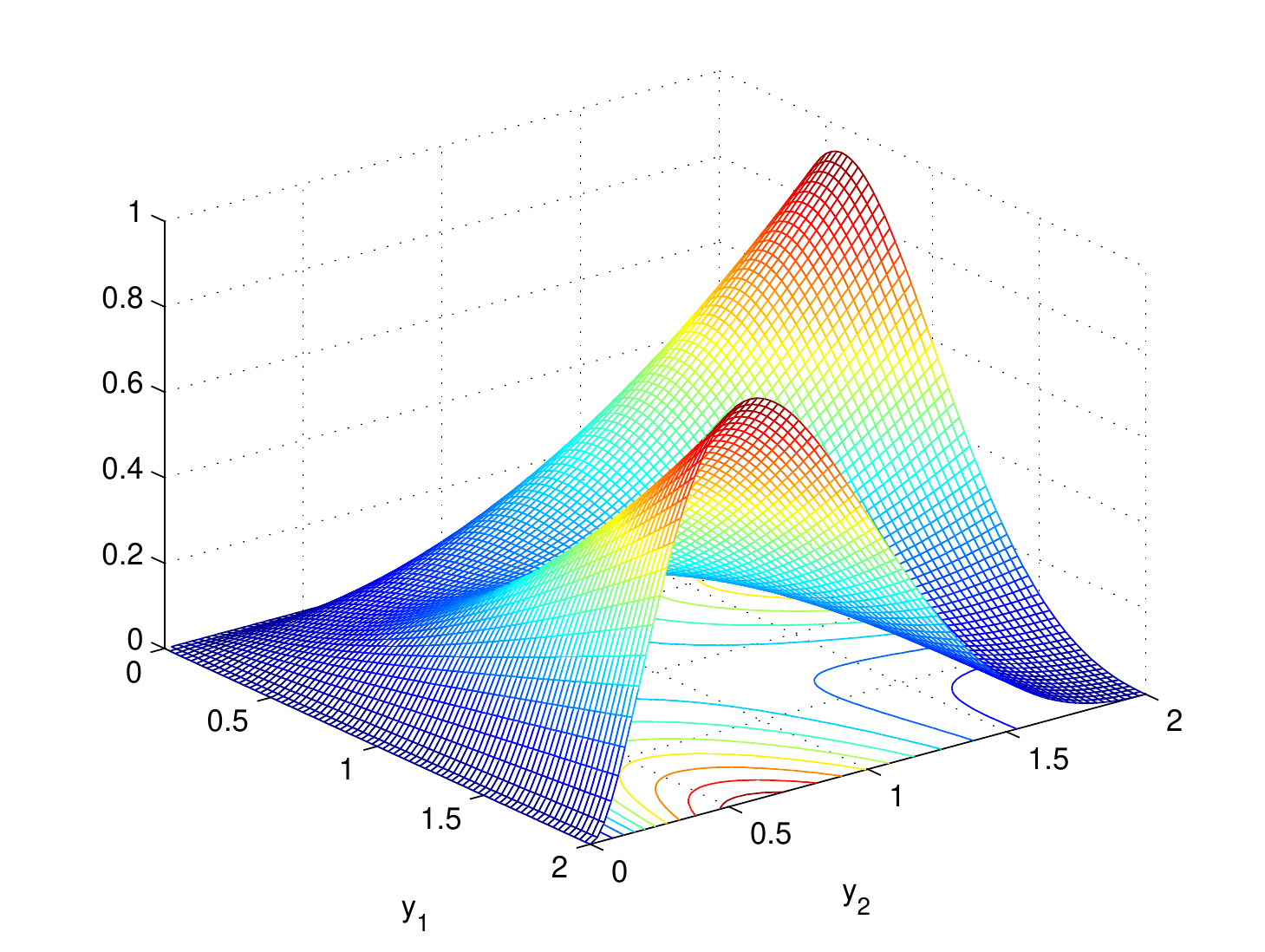}}
	\subfigure{
		\label{fig:subfig:b} 
		\includegraphics[width=4.8cm,height=4.8cm]{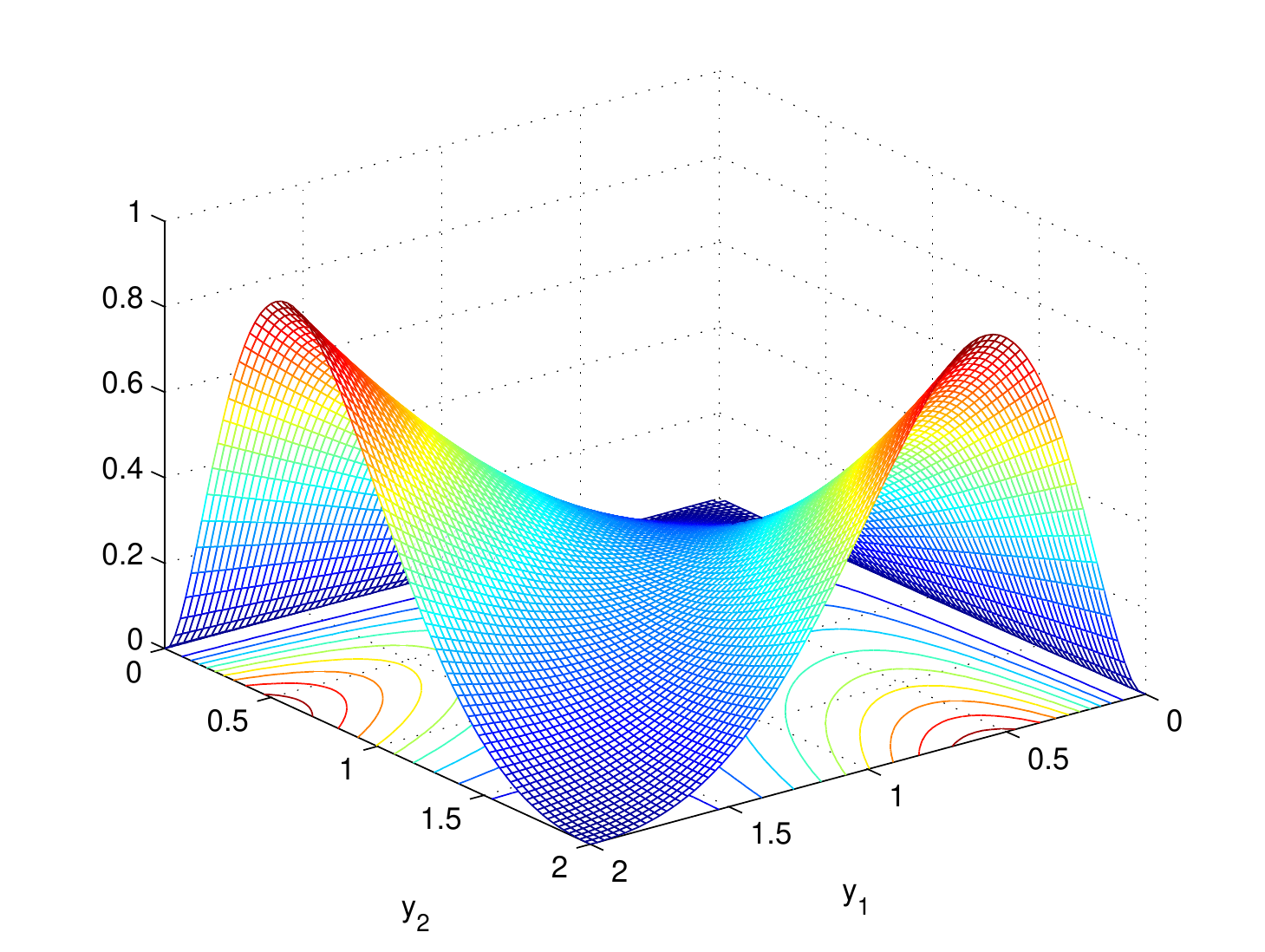}}
	\caption{Three-dimensional shaded figure of $\tilde{\nu}_n^2$ from three angles when $y_1\in(0,2)$ and $y_2\in(0,2)$. }\label{3Dfigure2}
\end{figure}

\begin{thm} \label{th1}
Under assumptions $(1)-(4)$, as  $\min\{p,n_1,n_2\}$ tends to infinity, we have
$$T_1:=\frac{\mathcal{L}-pl_{n}-\mu_{n}}{\nu_n}\stackrel{D}{\to} \mathcal{N}(0, 1),$$ where
$$
l_{n}=\frac{y_{n_2}}{y_{n_1}+y_{n_2}}-\frac{y_{n_2}-1}{y_{n_2}}{\delta_{(y_{n_2}>1)},}
~~~~~\mu_{n}=-\frac{\Delta_1y_{n_1}^2y_{n_2}^2h_n^2}{(y_{n_1}+y_{n_2})^4}-
\frac{\Delta_2y_{n_1}^2y_{n_2}^2h_n^2}{(y_{n_1}+y_{n_2})^4},$$and
$$\nu_{n}^2=\frac{2y_{n_1}^2y_{n_2}^2h_n^2}{(y_{n_1}+y_{n_2})^4}+
\frac{(y_{n_1}\Delta_1+y_{n_2}\Delta_2)y_{n_1}^2y_{n_2}^2h_n^4}{(y_{n_1}+y_{n_2})^6}.$$
Here, $\delta_{(\cdot)}$ denotes the indicator function and $\stackrel{D}{\to}
$ denotes convergence  in  distribution.
\end{thm}
%
The proof of this theorem is in the Appendix. The following theorem is based on $\widetilde{\mathcal{L}}$.
\begin{thm} \label{th2}
Under assumptions $(1)-(4)$, as $\min\{p,n_1,n_2\}\rightarrow \infty$, we have
$$T_2:=\frac{\tilde{\mathcal{L}}-p\tilde{l}_{n}-\tilde{\mu}_{n}}{\tilde{\nu}_n}\stackrel{D}{\to} \mathcal{N}(0, 1),  $$
where
\begin{gather*}
\tilde l_n =\frac{y_{n_1}y_{n_2}}{y_{n_1}+y_{n_2}}+\frac{(1-y_{n_1})y_{n_2}}{y_{n_1}^2}{\delta_{(y_{n_1}>1)}}+\frac{y_{n_1}(1-y_{n_2})}
{y_{n_2}^2}{\delta_{(y_{n_2}>1)}},
\end{gather*}
\begin{equation*}
\tilde{\mu}_{n}=\frac{y_{n_1}y_{n_2}h_n^2}{(y_{n_1}+y_{n_2})^2}+\frac{\Delta_1y_{n_1}^2y_{n_2}
h_n^2(h_n^2+2y_{n_2}(y_{n_2}-y_{n_1}))}{(y_{n_1}+y_{n_2})^4}+
\frac{\Delta_2y_{n_2}^2y_{n_1}h_n^2(h_n^2+2y_{n_1}(y_{n_1}-y_{n_2}))}{(y_{n_1}+y_{n_2})^4},
\end{equation*}
and
\begin{equation*}
\tilde{\nu}_n^2=\frac{4y_{n_1}^2y_{n_2}^2h_n^2(h_n^2+2(y_{n_1}-y_{n_2})^2)}{(y_{n_1}+y_{n_2})^4}
+\frac{(y_{n_1}\Delta_1+y_{n_2}\Delta_2)4y_{n_1}^2y_{n_2}^2h_n^4(y_{n_1}-y_{n_2})^2}{(y_{n_1}
+y_{n_2})^6}.
\end{equation*}
\end{thm}
The proof of this theorem is also in the Appendix.
\begin{rmk}
	In contrast to the LRT statistics in {Zhang et al. \cite{ZhangHuBai2017}}, the modified Pillai's trace statistics $T_1$ and $T_2$ proposed in this paper are feasible when $y_1=1$ or $y_2=1$, because all the limits, means and variances of the modified Pillai's trace statistics are continuous  in their definitional domain. We can also find from Figure 1 in Zhang et al. \cite{ZhangHuBai2017} that  when $y_1$ or $y_2$ are close to 1, the mean and variance of the LRT statistics  increase rapidly, resulting  in poor power. For  illustration purposes, three-dimensional shaded figures of $\tilde{\mu}_{n}$ and $\tilde{\nu}_n^2$ with $\Delta_1=\Delta_2=0$ are shown in Figure \ref{3Dfigure1} and Figure \ref{3Dfigure2}, respectively.
\end{rmk}
\begin{rmk}
	When $\Delta_l\neq0$,  Zhang et al. \cite{ZhangHuBai2017} provided the following estimators
 \begin{align*}
\hat\Delta_1=&(1-y)^2\frac{\sum_{j=1}^{n_1}[(z^{(1)}_j-\bar{z}^{(1)})'(c_{11}\bbS^
{z}_{1j}+c_{12}\bbS_2^{z})^{-1}(z^{(1)}_j-\bar{z}^{(1)})-\frac{p}{1-y}]^2}{pn_1}\nonumber\\
&-\frac{2}{1-y}\\
\hat\Delta_2=&(1-y)^2\frac{\sum_{j=1}^{n_2}[(z^{(2)}_j-\bar{z}^{(2)})'(c_{21}\bbS^
{z}_{1}+c_{22}\bbS_{2j}^{z})^{-1}(z^{(2)}_j-\bar{z}^{(2)})-\frac{p}{1-y}]^2}{pn_2},\nonumber\\
&-\frac{2}{1-y}
\end{align*}
where $y=\frac{p}{n_1+n_2-1}$, $c_{11}=\frac{n_1-1}{n_1+n_2-1}$, $c_{12}=\frac{n_2}{n_1+n_2-1}$,  $c_{21}=\frac{n_1}{n_1+n_2-1}$, $c_{22}=\frac{n_2-1}{n_1+n_2-1}$ and $\bbS^{z}_{lj}$ is the sample covariance matrix by removing the vector $z_{j}^{(l)}$ from the  $l$-th sample, $l=1,2$; Zhang et al. also showed that these estimators 
are weakly consistent and asymptotically unbiased under the null hypothesis. Under the alternative hypothesis, when $p<n_1$ and  $p<n_2$, the estimators are also applicable after applying a small modification. However, if  the linear model setting \eqref{linear model} does not hold or if $p\geq n_1+n_2$, then to the best of our knowledge, no consistent estimator of  $\Delta_i$ exists.
\end{rmk}
\section{Simulation}
In this section, we compare the modified Pillai's trace statistics $T_1$ and $T_2$ with four other statistics: $T_{lc}$, $T_{clx}$, $T_{zhb}^{1}$ and $T_{zhb}^{2}$ 
proposed by Li and Chen \cite{Li2012}, Cai et al. \cite{Cai2013} and Zhang et al. \cite{ZhangHuBai2017}, respectively. In the first subsection, we  compare the empirical sizes and powers of the proposed statistics $T_1$ and $T_2$ with $T_{lc}$, $T_{clx}$, $T_{zhb}^{1}$ and $T_{zhb}^{2}$ in some different settings. Because the properties of $T_1$ and $T_2$ are universal  and invariant under the null hypothesis $\Sigma_1=\Sigma_2$, and their powers depend only on the eigenvalues of $\Sigma_1\Sigma_2^{-1}$,  we organize four different targeted models and two different distributions. In the second subsection, we use the Jarque-Bera (J-B) and  Kolmogorov-Smirnov (K-S) tests to illustrate how well the proposed statistics fit their limiting distribution with a finite sample. 

\subsection{Comparisons of empirical sizes and powers}
First, we  consider a comparison with  $T_{lc}$ and $T_{clx}$. To test hypothesis \eqref{hypo}, we randomly generate $x_{j}^{(l)}$ from a standard multivariate normal distribution $N(0_p,I_p)$, and let  $z_{j}^{(l)}=\Sigma_{l}^{1/2} x_{j}^{(l)}$. Simultaneously, to realize the empirical size and power of the test, we define $$\Sigma_2=(1+\frac{\delta}{n_1})\Sigma_1.$$ When $\delta=0$, we achieve the empirical size. $\Sigma_1$ follows the following four models:
\begin{enumerate}[Model 1:]
\item $\Sigma_1=I_p$;
\item $\Sigma_1=Diag(p^2,1,...,1)$;
\item $\Sigma_1=D^{1/2}\Sigma^{\ast}D^{1/2}$, where $D=Diag(d_{1}, d_{2},...,d_{p})$, $d_{i}=$ Unif$(0.5,2.5)$  and $\Sigma^{\ast}=(\sigma_{ij}^{\ast})$, $\sigma_{ii}^{\ast}=1$, $\sigma_{ij}^{\ast}=0.5$ for $5(k-1)+1\leq i\neq j\leq5k$, $k=1,...,\lceil p/5 \rceil$; otherwise, $\sigma_{ij}^{\ast}=0$, $i,j=1,2,...,p$;
\item $\Sigma_1=(0.5 I_p+0.5\textbf{1}_p \textbf{1}_p ')$, where $\textbf{1}_p$ is an all-ones vector.
\end{enumerate}
We set the actual size to $5\%$ for 1,000 repetitive simulations.
The sample sizes $(n_1,n_2)$ increase from $(25, 35)$ to $(400, 560)$.
To fulfill all the conditions, $p$ is selected under the assumption that $p<\min\{n_1,n_2\}$, $\min\{n_1,n_2\}<p<\max\{n_1,n_2\}$ or $p>\max\{n_1,n_2\}.$
The $\Sigma_{l}$ in Model 1 comes from Li and Chen \cite{Li2012} and satisfies the corresponding assumption that $\Sigma_{l}$ has a moderate eigenvalue.
However, Model 2 fails to satisfy this demand.
Model 3 comes from Cai et al. \cite{Cai2013} and is quite sparse. Because $T_{clx}$  is established on the corresponding elements of the two covariance matrices and requires a sparsity condition, for comparison, we choose the Model 4, which is unable to satisfy the sparse condition. The simulation results for Models 1--4 are reported in Tables \ref{tab1}-\ref{tab4}, respectively.
\begin{table}[h]
\centering
\footnotesize
\begin{tabular}{ccccccccccccc}
\\ \hline
{($n_1$,$n_2$)}& Method   &\multicolumn{3}{c}{Size ($\delta=0$)}& &\multicolumn{3}{c}{Power ($\delta=5$)}& &\multicolumn{3}{c}{Power ($\delta=10$)}\\ \hline
&{p}&20&30&40& &20&30&40& &20&30&40   \\ \cline{3-5} \cline{7-9}\cline{11-13}
\multirow{3}{3em}{(25,35)}
&$T_1$    &0.051 &0.055 &0.047 & &0.948 &0.963 &0.943  & &1     &1     &1          \\
&$T_2$    &0.046 &0.052 &0.042 & &0.325 &0.436 &0.607  & &0.860 &0.952 &0.999      \\
&$T_{lc}$ &0.064 &0.070 &0.061 & &0.220 &0.205 &0.207  & &0.770 &0.818 &0.806      \\ \hline	
	
& {p}&40&60&80& &40&60&80& &40&60&80  \\ \cline{3-5} \cline{7-9}\cline{11-13}
\multirow{3}{3em}{(50,70)}

&$T_1$    &0.049 &0.051 &0.048 & &0.955 &0.964 &0.959   & &1     &1     &1          \\
&$T_2$    &0.044 &0.054 &0.042 & &0.260 &0.431 &0.602   & &0.832 &0.963 &0.994      \\
&$T_{lc}$ &0.062 &0.053 &0.073 & &0.099 &0.111 &0.112   & &0.516 &0.471 &0.501      \\ \hline	

& {p}&80&120&160& &80&120&160& &80&120&160  \\ \cline{3-5} \cline{7-9}\cline{11-13}
\multirow{3}{3em}{(100,140)}

&$T_1$    &0.047 &0.046 &0.063 & &0.973 &0.982 &0.971   & &1     &1     &1       \\
&$T_2$    &0.045 &0.055 &0.054 & &0.256 &0.434 &0.605   & &0.768 &0.965 &0.992   \\
&$T_{lc}$ &0.059 &0.062 &0.051 & &0.069 &0.071 &0.059   & &0.216 &0.208 &0.215   \\ \hline	

& {p}&160&240&320& &160&240&320& &160&240&320  \\ \cline{3-5} \cline{7-9}\cline{11-13}
\multirow{3}{3em}{(200,280)}

&$T_1$    &0.043 &0.055 &0.049 & &0.973 &0.983 &0.978  & &1     &1     &1         \\
&$T_2$    &0.052 &0.044 &0.049 & &0.240 &0.417 &0.600  & &0.756 &0.950 &0.991     \\
&$T_{lc}$ &0.055 &0.048 &0.064 & &0.054 &0.070 &0.058  & &0.098 &0.088 &0.102     \\ \hline	

& {p}&320&480&640& &320&480&640& &{320}&480&640  \\ \cline{3-5} \cline{7-9}\cline{11-13}
\multirow{3}{3em}{(400,560)}

&$T_1$    &0.056 &0.054 &0.045 & &0.971 &0.981 &0.974 & &1     &1     &1       \\
&$T_2$    &0.046 &0.051 &0.051 & &0.253 &0.434 &0.566 & &0.715 &0.959 &0.990    \\
&$T_{lc}$ &0.059 &0.053 &0.049 & &0.052 &0.047 &0.051 & &0.048 &0.061 &0.052   \\ \hline	

\end{tabular}
\caption{Empirical size and power from 1,000 repeated simulations comparing $T_1$, $T_2$ and $T_{lc}$ based on Model 1 under the normal assumption.}\label{tab1}
\end{table}

\begin{table}[h]
\centering
\footnotesize
\begin{tabular}{ccccccccccccc}
\\ \hline
{($n_1$,$n_2$)}& Method   &\multicolumn{3}{c}{Size ($\delta=0$)}& &\multicolumn{3}{c}{Power ($\delta=5$)}& &\multicolumn{3}{c}{Power ($\delta=10$)}\\ \hline
& {p}&20&30&40& &20&30&40& &20&30&40   \\ \cline{3-5} \cline{7-9}\cline{11-13}
\multirow{3}{3em}{(25,35)}

&$T_1$    &0.048 &0.038 &0.046 & &0.949 &0.966 &0.948 & &1     &1     &1         \\
&$T_2$    &0.045 &0.043 &0.036 & &0.323 &0.439 &0.614 & &0.853 &0.969 &0.992    \\
&$T_{lc}$ &0.091 &0.086 &0.093 & &0.205 &0.226 &0.216 & &0.477 &0.483 &0.514  \\ \hline	

& {p}&40&60&80& &40&60&80& &40&60&80  \\ \cline{3-5} \cline{7-9}\cline{11-13}
\multirow{3}{3em}{(50,70)}

&$T_1$    &0.040 &0.065 &0.045 & &0.956 &0.966 &0.960 & &1     &1     &1        \\
&$T_2$    &0.056 &0.059 &0.050 & &0.257 &0.433 &0.622 & &0.825 &0.956 &0.993    \\
&$T_{lc}$ &0.074 &0.097 &0.077 & &0.167 &0.166 &0.166 & &0.345 &0.385 &0.374  \\ \hline	

& {p}&80&120&160& &80&120&160& &80&120&160  \\ \cline{3-5} \cline{7-9}\cline{11-13}
\multirow{3}{3em}{(100,140)}

&$T_1$    &0.043 &0.046 &0.052 & &0.971 &0.970 &0.970 & &1     &1     &1        \\
&$T_2$    &0.047 &0.048 &0.053 & &0.276 &0.422 &0.584 & &0.809 &0.967 &0.994     \\
&$T_{lc}$ &0.092 &0.100 &0.076 & &0.111 &0.114 &0.119 & &0.240 &0.253 &0.263  \\ \hline

& {p}&160&240&320& &160&240&320& &160&240&320  \\ \cline{3-5} \cline{7-9}\cline{11-13}
\multirow{3}{3em}{(200,280)}

&$T_1$    &0.054 &0.051 &0.055 & &0.971 &0.981 &0.975 & &1     &1     &1        \\
&$T_2$    &0.053 &0.048 &0.054 & &0.270 &0.403 &0.602 & &0.759 &0.962 &0.998      \\
&$T_{lc}$ &0.093 &0.102 &0.074 & &0.098 &0.111 &0.131 & &0.175 &0.181 &0.171   \\ \hline	

& {p}&320&480&640& &320&480&640& &320&480&640  \\ \cline{3-5} \cline{7-9}\cline{11-13}
\multirow{3}{3em}{(400,560)}

&$T_1$    &0.051 &0.037 &0.038 & &0.979 &0.989 &0.981 & &1     &1     &1        \\
&$T_2$    &0.052 &0.046 &0.044 & &0.250 &0.417 &0.599 & &0.719 &0.961 &0.997    \\
&$T_{lc}$ &0.085 &0.081 &0.090 & &0.090 &0.082 &0.086 & &0.118 &0.126 &0.129   \\ \hline	

\end{tabular}
\caption{Empirical size and power from 1,000 repeated simulations comparing $T_1$, $T_2$ and $T_{lc}$ based on Model 2 under the normal assumption.}\label{tab2}
\end{table}

\begin{table}[h]
\centering
\footnotesize
\begin{tabular}{ccccccccccccc}
\\ \hline
{($n_1$,$n_2$)}& Method   &\multicolumn{3}{c}{Size ($\delta=0$)}& &\multicolumn{3}{c}{Power ($\delta=5$)}& &\multicolumn{3}{c}{Power ($\delta=10$)}\\ \hline
& {p}&20&30&40& &20&30&40& &20&30&40   \\ \cline{3-5} \cline{7-9}\cline{11-13}
\multirow{3}{3em}{(25,35)}

&$T_1$    &0.051 &0.054 &0.055 & &0.951 &0.964 &0.940 & &1     &1     &1         \\
&$T_2$    &0.046 &0.043 &0.048 & &0.337 &0.423 &0.623 & &0.871 &0.964 &0.995     \\
&$T_{clx}$  &0.070 &0.091 &0.075 & &0.084 &0.090 &0.087 & &0.164 &0.173 &0.158  \\ \hline	

& {p}&40&60&80& &40&60&80& &40&60&80  \\ \cline{3-5} \cline{7-9}\cline{11-13}
\multirow{3}{3em}{(50,70)}

&$T_1$    &0.051 &0.045 &0.056 & &0.964 &0.971 &0.967 & &1     &1     &1        \\
&$T_2$    &0.054 &0.042 &0.043 & &0.296 &0.431 &0.615 & &0.832 &0.975 &0.995     \\
&$T_{clx}$  &0.045 &0.054 &0.055 & &0.065 &0.049 &0.071 & &0.078 &0.081 &0.073  \\ \hline	

& {p}&80&120&160& &80&120&160& &80&120&160  \\ \cline{3-5} \cline{7-9}\cline{11-13}
\multirow{3}{3em}{(100,140)}

&$T_1$    &0.052 &0.050 &0.049 & &0.973 &0.978 &0.970 & &1     &1     &1         \\
&$T_2$    &0.053 &0.042 &0.064 & &0.281 &0.422 &0.599 & &0.805 &0.954 &0.994     \\
&$T_{clx}$  &0.042 &0.044 &0.050 & &0.052 &0.045 &0.031 & &0.060 &0.049 &0.055  \\ \hline	

& {p}&160&240&320& &160&240&320& &160&240&320  \\ \cline{3-5} \cline{7-9}\cline{11-13}
\multirow{3}{3em}{(200,280)}

&$T_1$    &0.054 &0.039 &0.048 & &0.971 &0.977 &0.976 & &1     &1     &1         \\
&$T_2$    &0.060 &0.036 &0.048 & &0.258 &0.457 &0.625 & &0.762 &0.957 &0.993    \\
&{$T_{clx}$} &0.050 &0.032 &0.043 & &0.045 &0.046 &0.043 & &0.050 &0.052 &0.044  \\ \hline	

& {p}&320&480&640& &320&480&640& &{320}&480&640  \\ \cline{3-5} \cline{7-9}\cline{11-13}
\multirow{3}{3em}{(400,560)}

&$T_1$    &0.046 &0.047 &0.040 & &0.966 &0.987 &0.971 & &1     &1     &1          \\
&$T_2$    &0.054 &0.046 &0.038 & &0.254 &0.423 &0.621 & &0.743 &0.950 &0.994    \\
&{$T_{clx}$} &0.044 &0.046 &0.044 & &0.038 &0.044 &0.036 & &0.045 &0.054 &0.050  \\ \hline

\end{tabular}
\caption{Empirical size and power from 1,000 repeated simulations comparing $T_1$, $T_2$ and $T_{clx}$ based on Model 3 under the normal assumption.}\label{tab3}
\end{table}

\begin{table}[h]
\centering
\footnotesize
\begin{tabular}{ccccccccccccc}
\\ \hline
{($n_1$,$n_2$)}& Method   &\multicolumn{3}{c}{Size ($\delta=0$)}& &\multicolumn{3}{c}{Power ($\delta=5$)}& &\multicolumn{3}{c}{Power ($\delta=10$)}\\ \hline
& {p}&20&30&40& &20&30&40& &20&30&40   \\ \cline{3-5} \cline{7-9}\cline{11-13}
\multirow{3}{3em}{(25,35)}

&$T_1$    &0.054 &0.052 &0.043 & &0.948 &0.962 &0.940 & &1     &1     &1          \\
&$T_2$    &0.045 &0.040 &0.045 & &0.325 &0.442 &0.620 & &0.869 &0.960 &0.997    \\
&$T_{clx}$  &0.045 &0.039 &0.056 & &0.065 &0.061 &0.065 & &0.148 &0.154 &0.144  \\ \hline	

& {p}&40&60&80& &40&60&80& &40&60&80  \\ \cline{3-5} \cline{7-9}\cline{11-13}
\multirow{3}{3em}{(50,70)}

&$T_1$    &0.048 &0.045 &0.053 & &0.954 &0.977 &0.948 & &1     &1     &1         \\
&$T_2$    &0.055 &0.048 &0.053 & &0.248 &0.455 &0.582 & &0.828 &0.962 &0.992    \\
&$T_{clx}$  &0.039 &0.028 &0.036 & &0.034 &0.031 &0.028 & &0.088 &0.056 &0.045  \\ \hline	

& {p}&80&120&160& &80&120&160& &80&120&160  \\ \cline{3-5} \cline{7-9}\cline{11-13}
\multirow{3}{3em}{(100,140)}

&$T_1$    &0.067 &0.042 &0.042 & &0.972 &0.980 &0.970 & &1     &1     &1          \\
&$T_2$    &0.051 &0.039 &0.052 & &0.252 &0.472 &0.606 & &0.780 &0.959 &0.996     \\
&$T_{clx}$  &0.035 &0.022 &0.021 & &0.029 &0.023 &0.010 & &0.041 &0.027 &0.032  \\ \hline	

& {p}&160&240&320& &160&240&320& &160&240&320  \\ \cline{3-5} \cline{7-9}\cline{11-13}
\multirow{3}{3em}{(200,280)}

&$T_1$    &0.047 &0.048 &0.045 & &0.974 &0.984 &0.976 & &1     &1     &1         \\
&$T_2$    &0.054 &0.050 &0.048 & &0.270 &0.441 &0.587 & &0.736 &0.948 &0.986   \\
&{$T_{clx}$} &0.015 &0.023 &0.015 & &0.012 &0.008 &0.014 & &0.033 &0.013 &0.014  \\ \hline	

& {p}&320&480&640& &320&480&640& &{320}&480&640  \\ \cline{3-5} \cline{7-9}\cline{11-13}
\multirow{3}{3em}{(400,560)}

&$T_1$    &0.053 &0.052 &0.042 & &0.971 &0.983 &0.973 & &1     &1     &1          \\
&$T_2$    &0.051 &0.053 &0.040 & &0.231 &0.439 &0.559 & &0.740 &0.952 &0.993      \\
&{$T_{clx}$} &0.021 &0.017 &0.012 & &0.017 &0.010 &0.013 & &0.014 &0.011 &0.005  \\ \hline

\end{tabular}
\caption{Empirical size and power from 1,000 repeated simulations comparing $T_1$, $T_2$ and $T_{clx}$ based on Model 4 under the normal assumption.}\label{tab4}
\end{table}
%
%

\begin{figure}[h]
	\begin{center}
        \includegraphics[width=7.4cm]{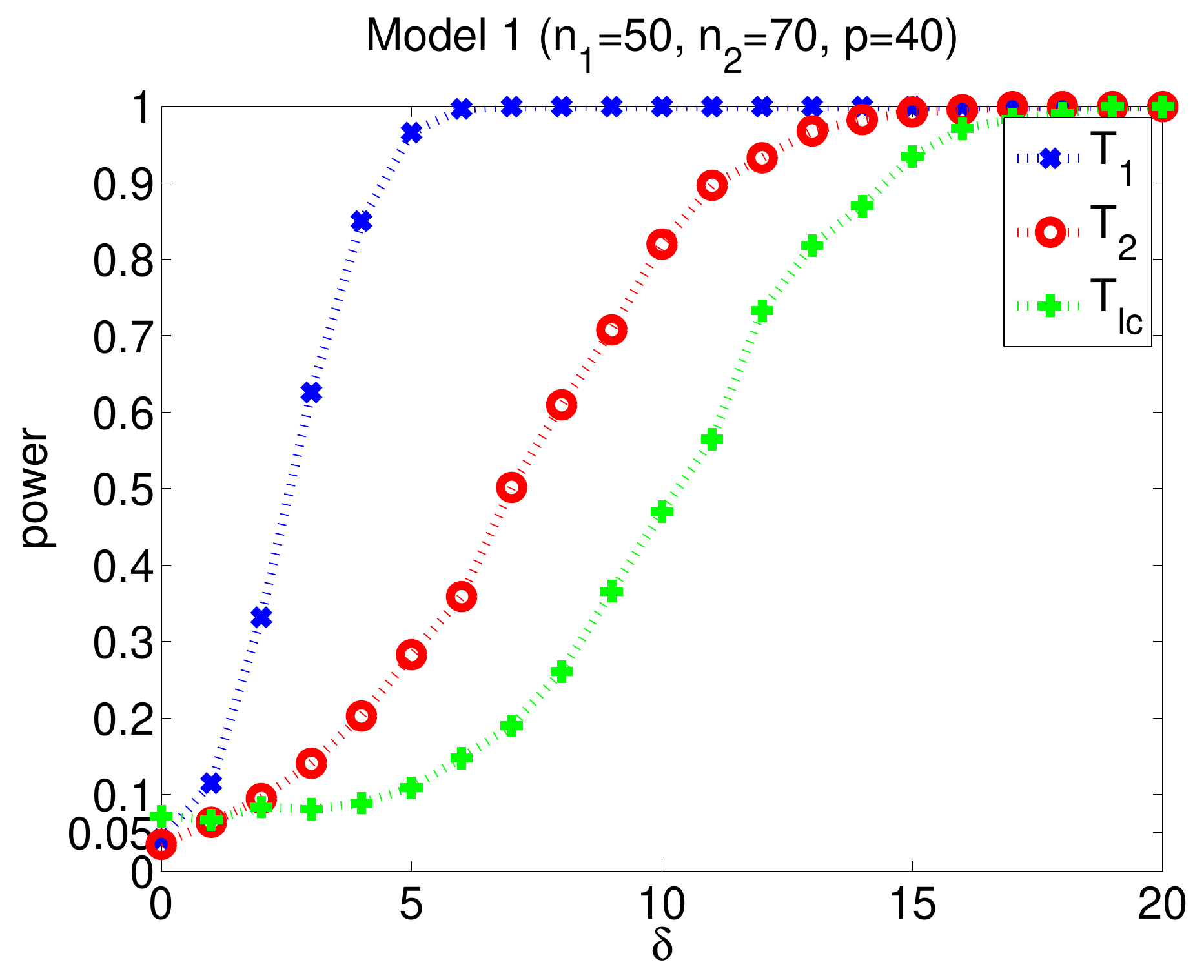}
        \includegraphics[width=7.4cm]{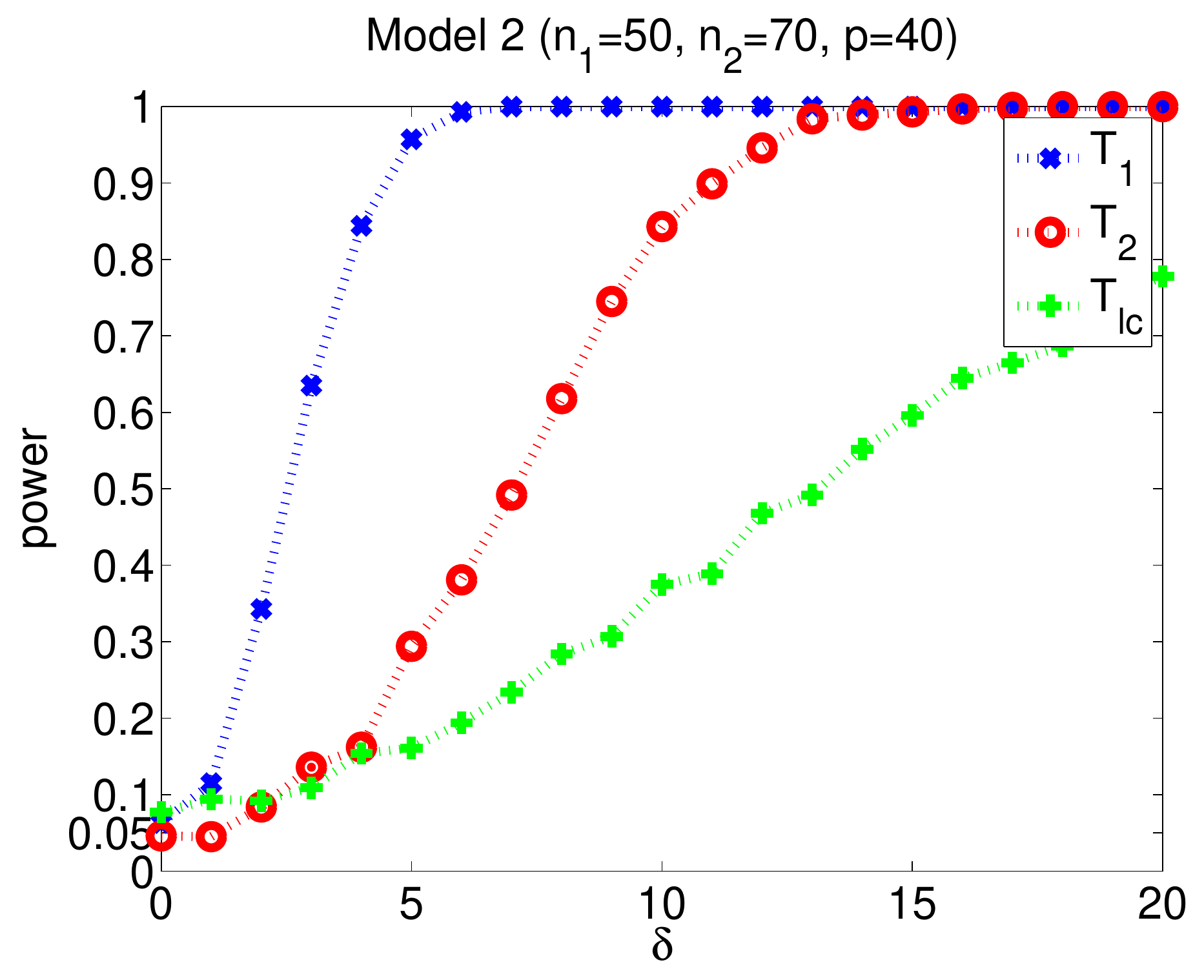}
 \caption{Scatter diagram of the empirical power for $T_1$, $T_2$ and $T_{lc}$ based on Model 1 and Model 2 under the normal assumption. }
 \label{2Dscatter1}
	\end{center}
\end{figure}

\begin{figure}[h]
	\begin{center}
		\includegraphics[width=7.4cm]{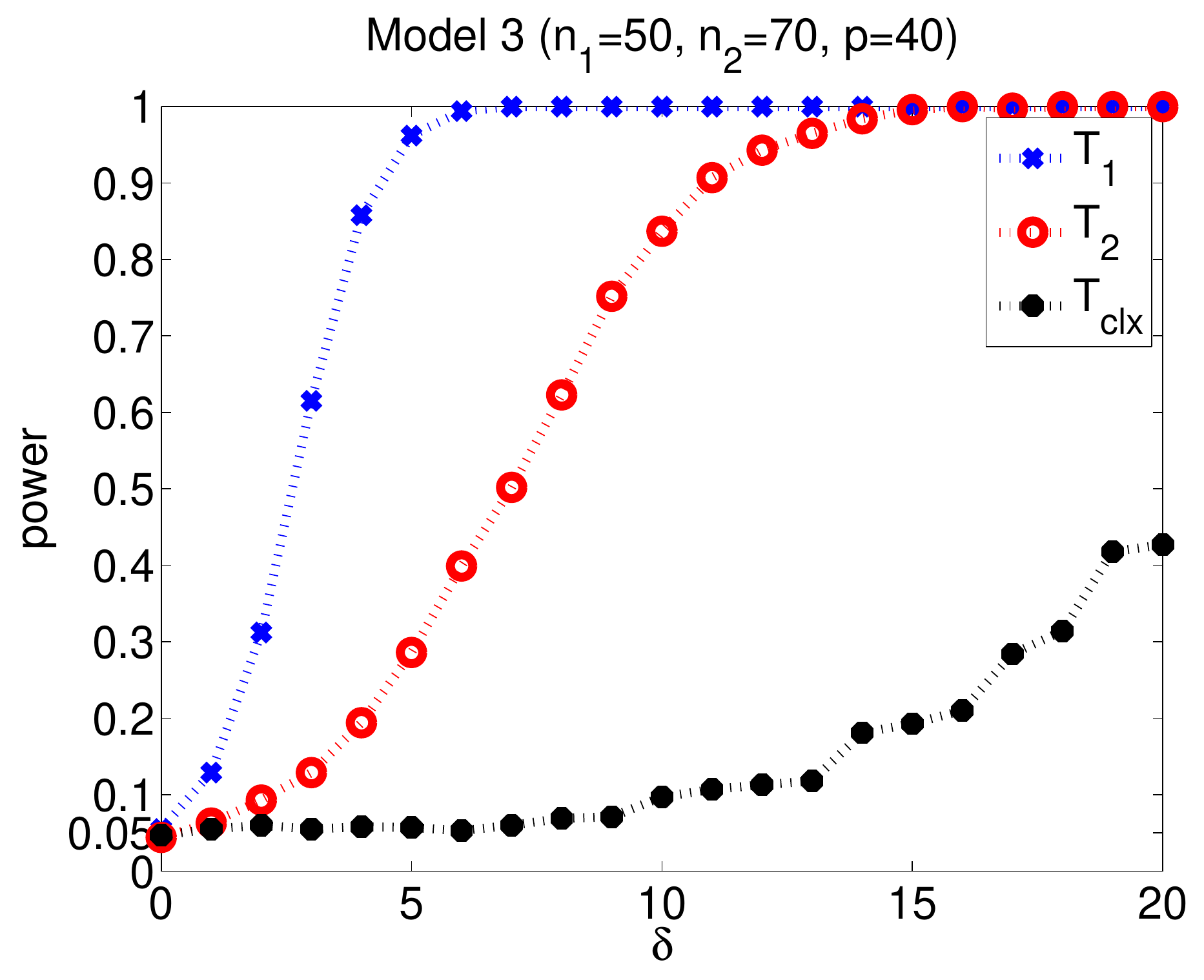}
        \includegraphics[width=7.4cm]{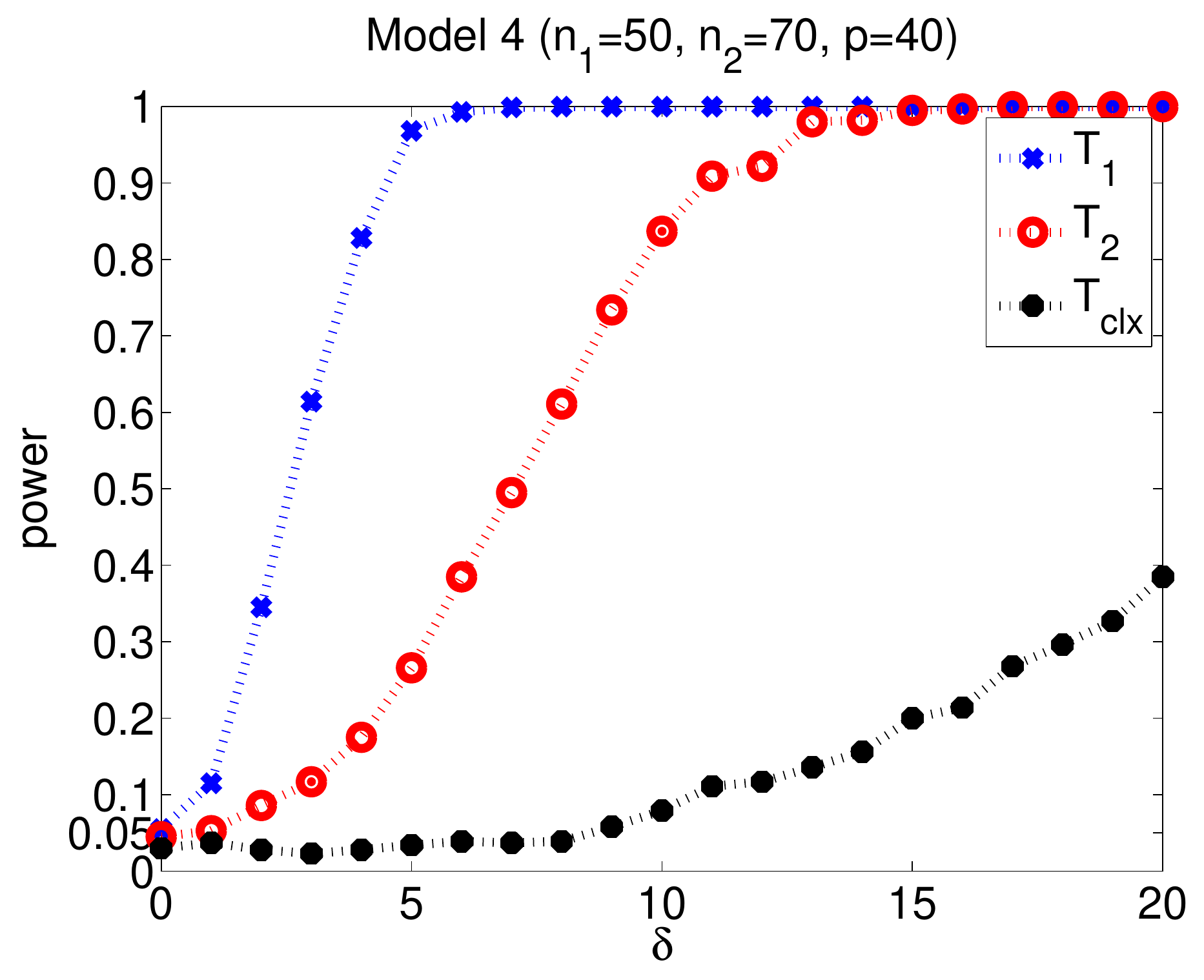}
 \caption{Scatter diagram of the empirical power for $T_1$, $T_2$ and $T_{clx}$ based on Model 3 and Model 4 under the normal assumption. }
 \label{2Dscatter2}
	\end{center}
\end{figure}

Tables \ref{tab1}-\ref{tab4} show that when $\min\{n_1,n_2,p\}$ is large, the modified Pillai's trace statistics have relatively good sizes, and in terms of  $\delta$, the modified Pillai's trace statistics are more efficient than are $T_{lc}$ and $T_{clx}$. That is, $T_1$ and $T_2$ can be utilized to distinguish the two different covariance matrices even when $\delta$ is  small.
To illustrate the outstanding efficiency of the modified Pillai's trace statistics, we show scatter plots in Figures \ref{2Dscatter1}-\ref{2Dscatter2} under $(n_1,n_2,p)=(50,70,40)$ as $\delta$ increases from $0$ to $20.$ The modified Pillai's trace statistics tend toward $1$ more quickly than do either $T_{lc}$ or $T_{clx}$.

\begin{figure}[htbp]
	\begin{center}
		\includegraphics[height=4.7cm, width=7.4cm]{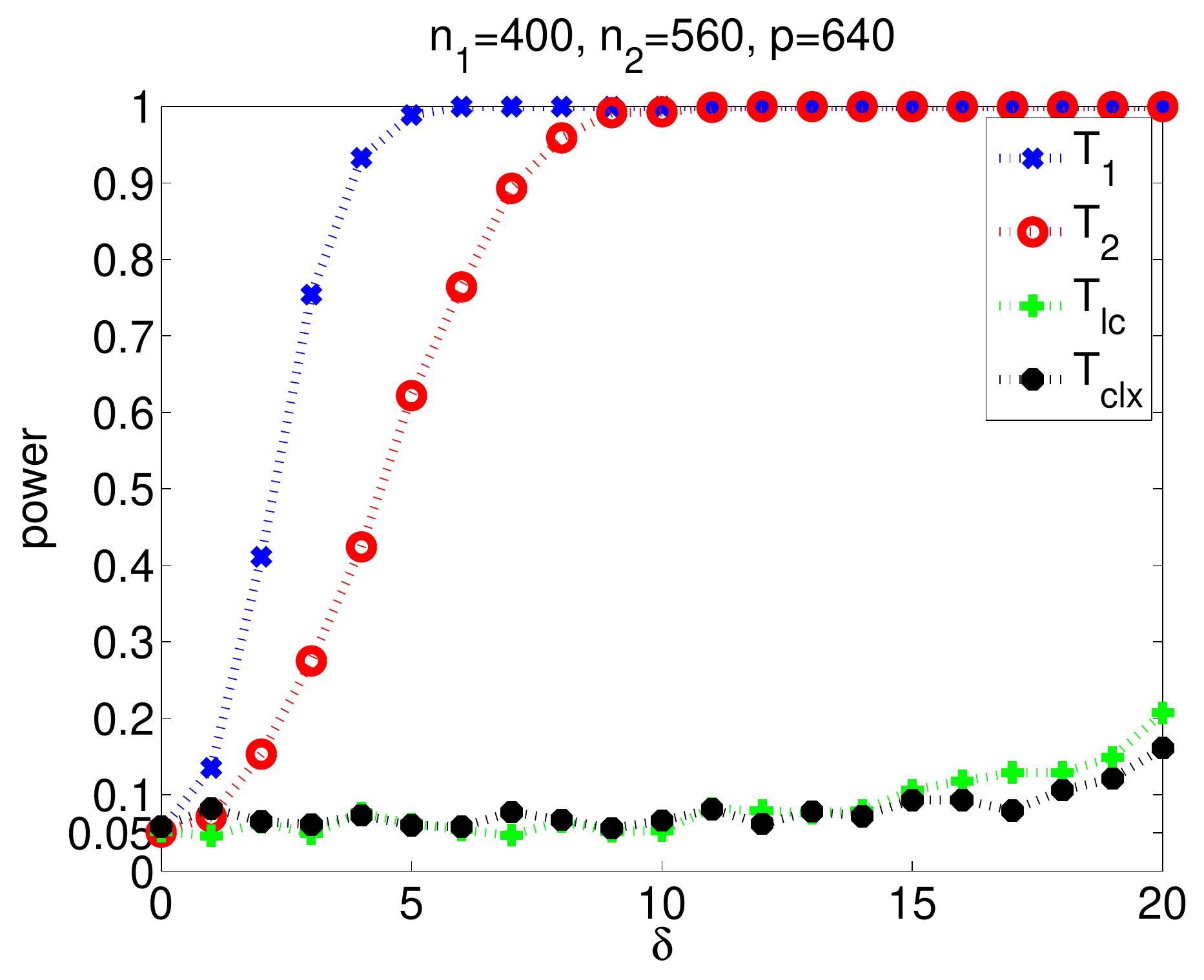}\includegraphics[height=4.7cm, width=7.4cm]{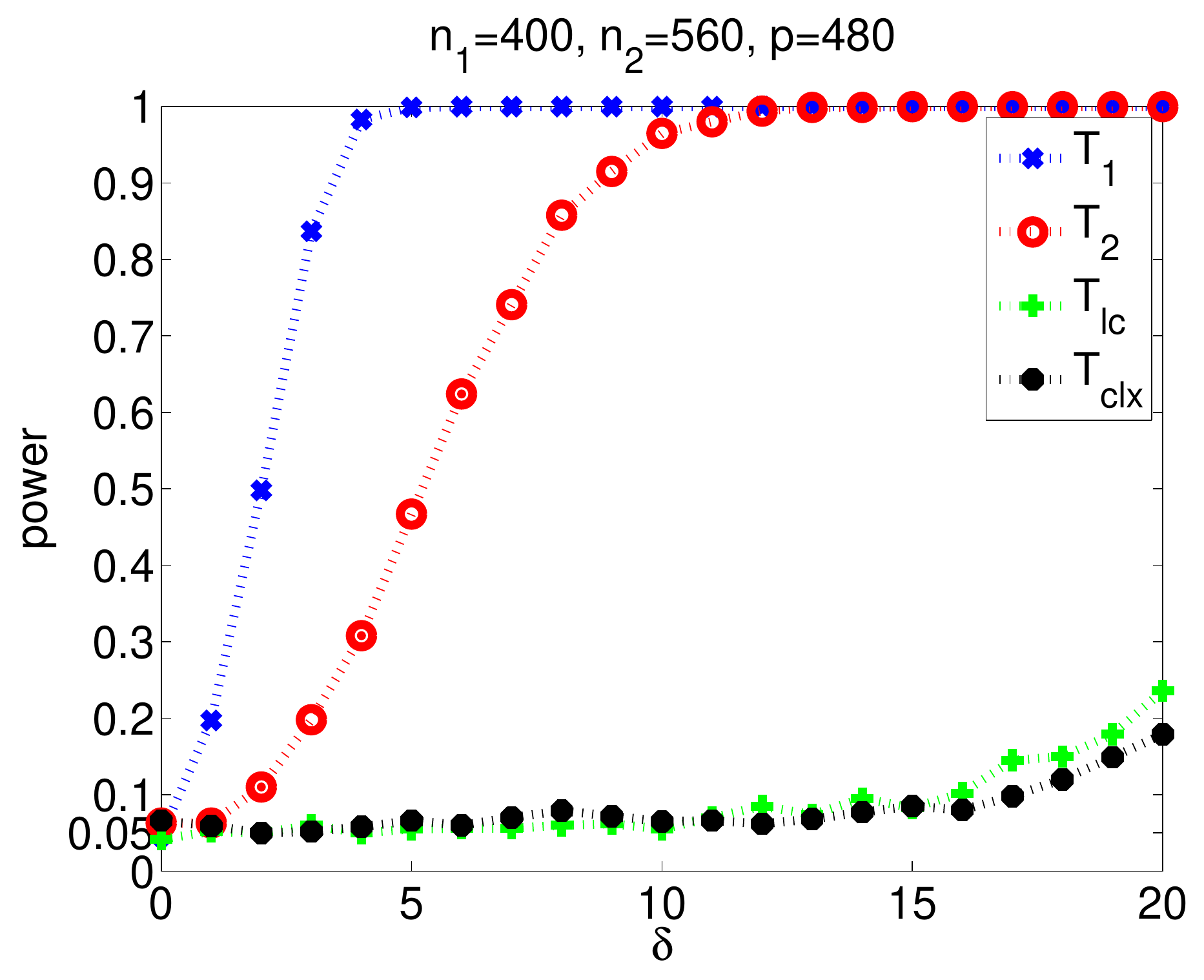}
        \includegraphics[height=4.7cm, width=7.4cm]{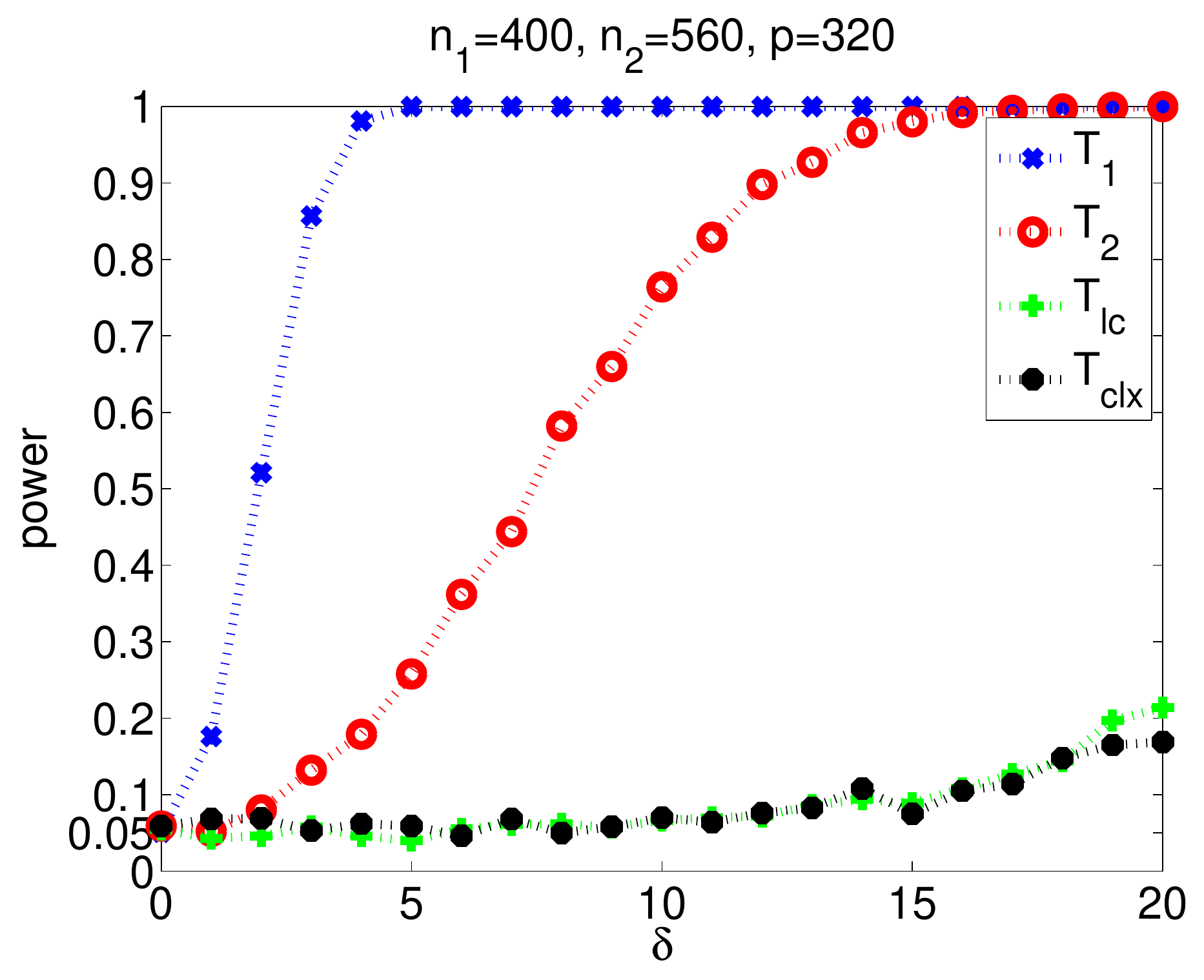}\includegraphics[height=4.7cm, width=7.4cm]{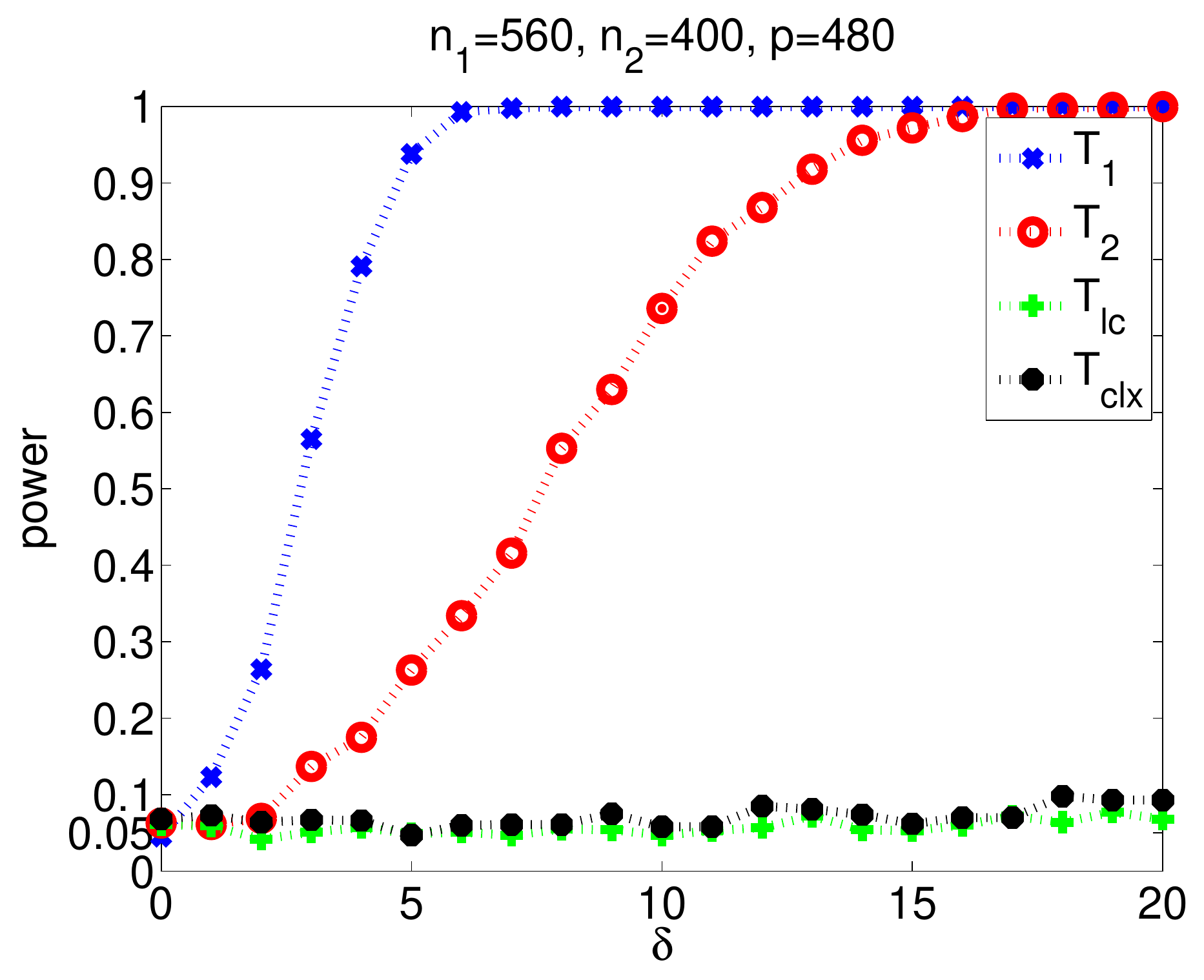}
    \caption{Scatter diagram of the empirical power for $T_1$, $T_2$, $T_{lc}$ and $T_{clx}$ based on Model 1 under the uniform assumption. }
    \label{unifMo1}
    \end{center}
\end{figure}

\begin{figure}[htbp]
	\begin{center}
		\includegraphics[height=4.7cm, width=7.4cm]{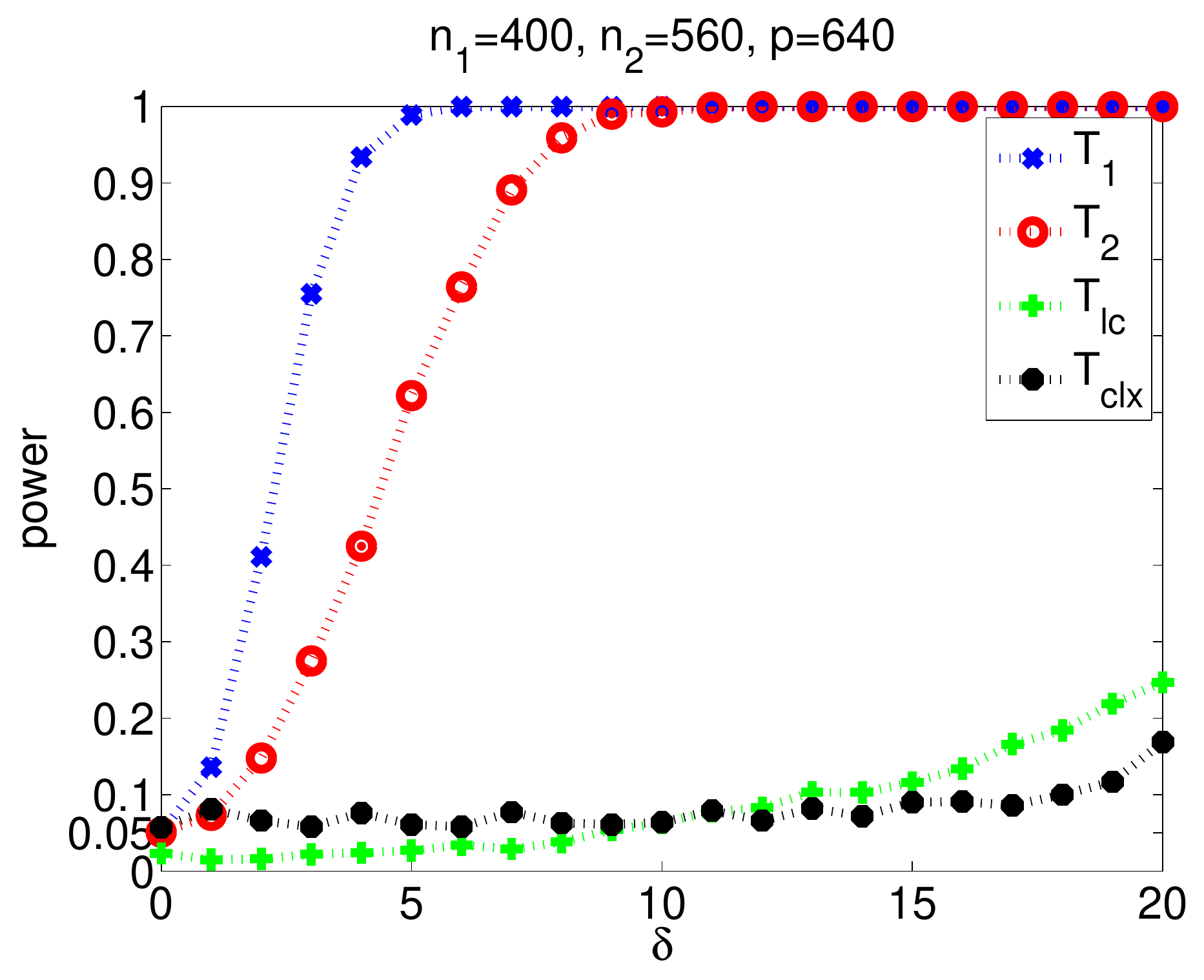} \includegraphics[height=4.7cm, width=7.4cm]{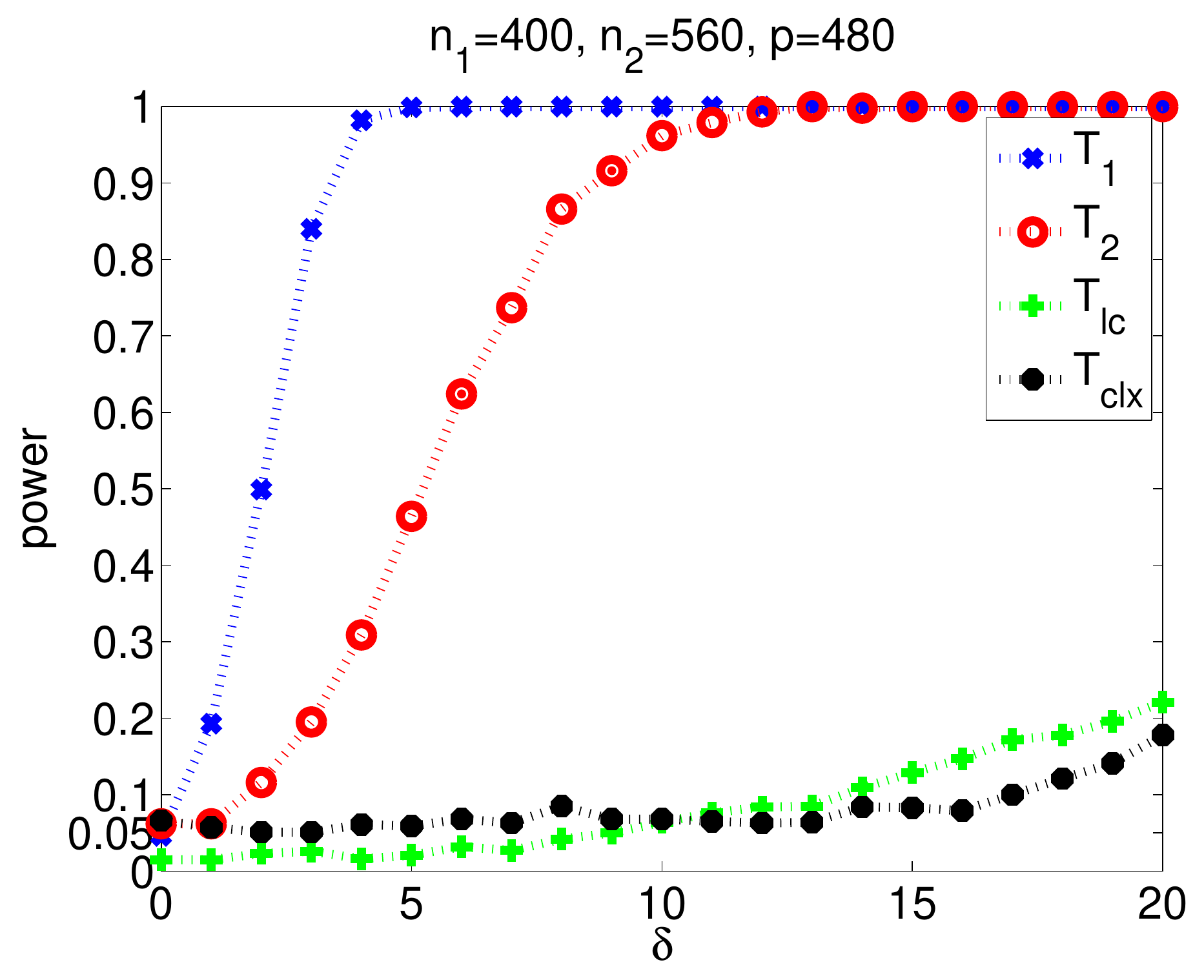}
        \includegraphics[height=4.7cm, width=7.4cm]{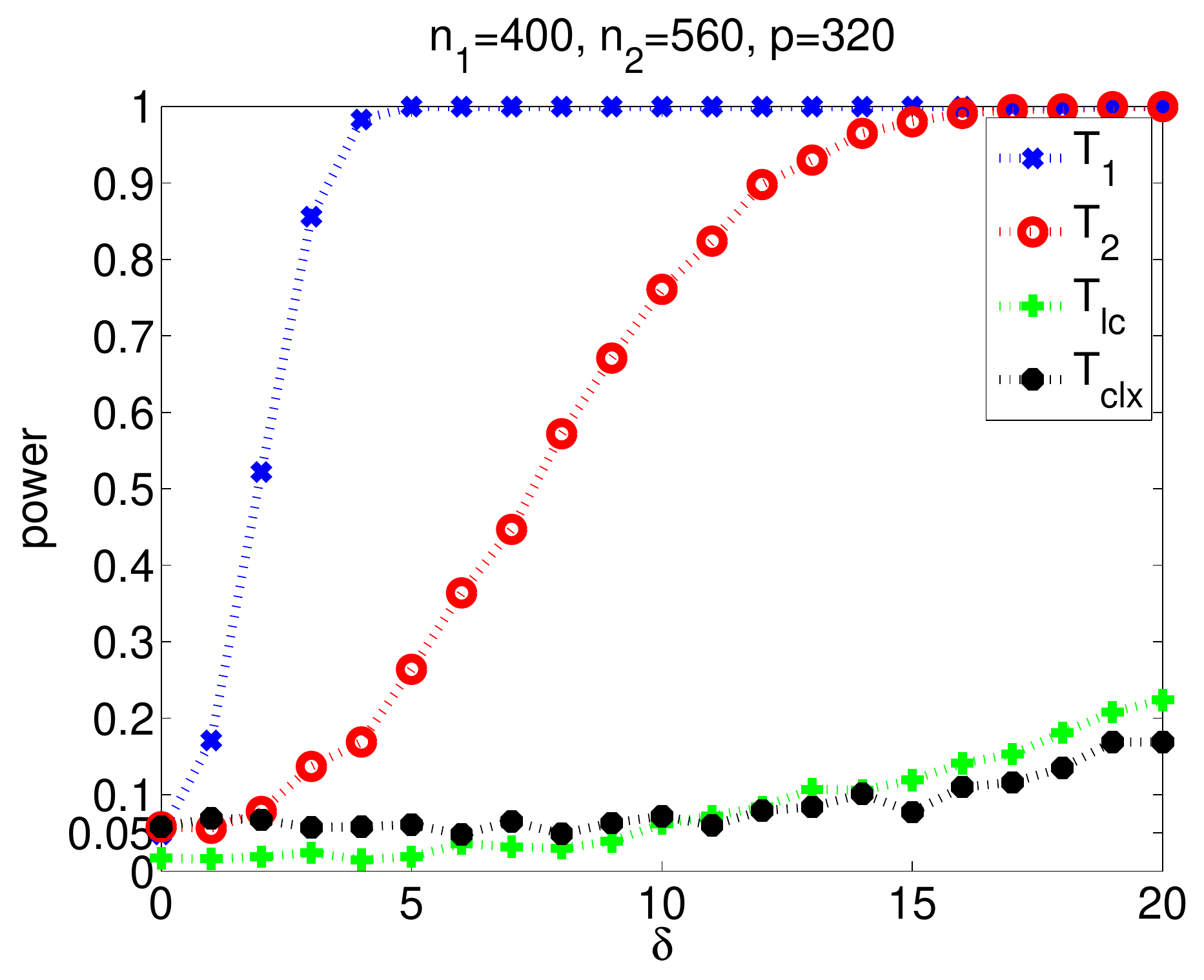} \includegraphics[height=4.7cm, width=7.4cm]{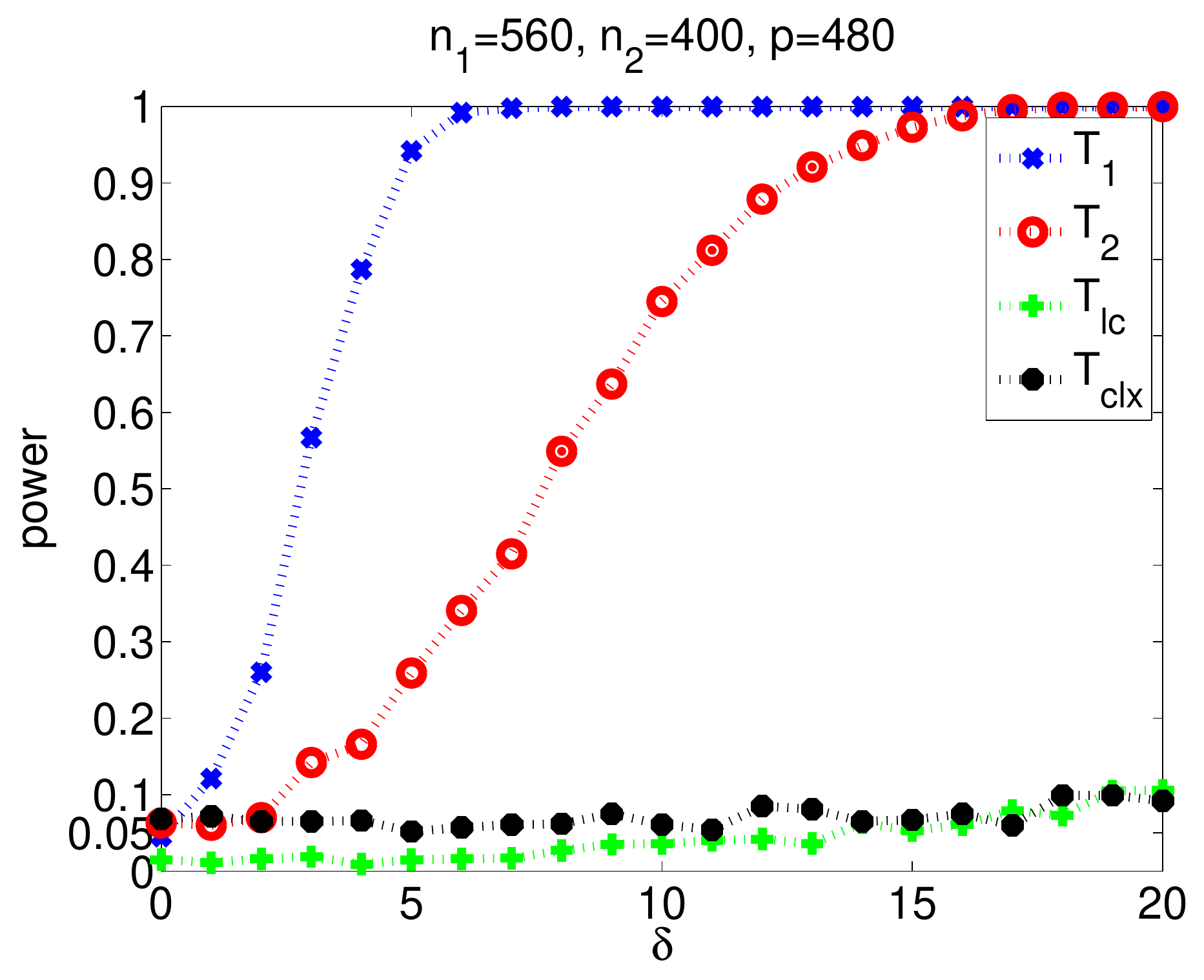}
    \caption{Scatter diagram of the empirical power for $T_1$, $T_2$, $T_{lc}$ and $T_{clx}$ based on Model 2 under the uniform assumption. }
    \label{unifMo2}
    \end{center}
\end{figure}

\begin{figure}[htbp]
	\begin{center}
		\includegraphics[height=4.7cm, width=7.4cm]{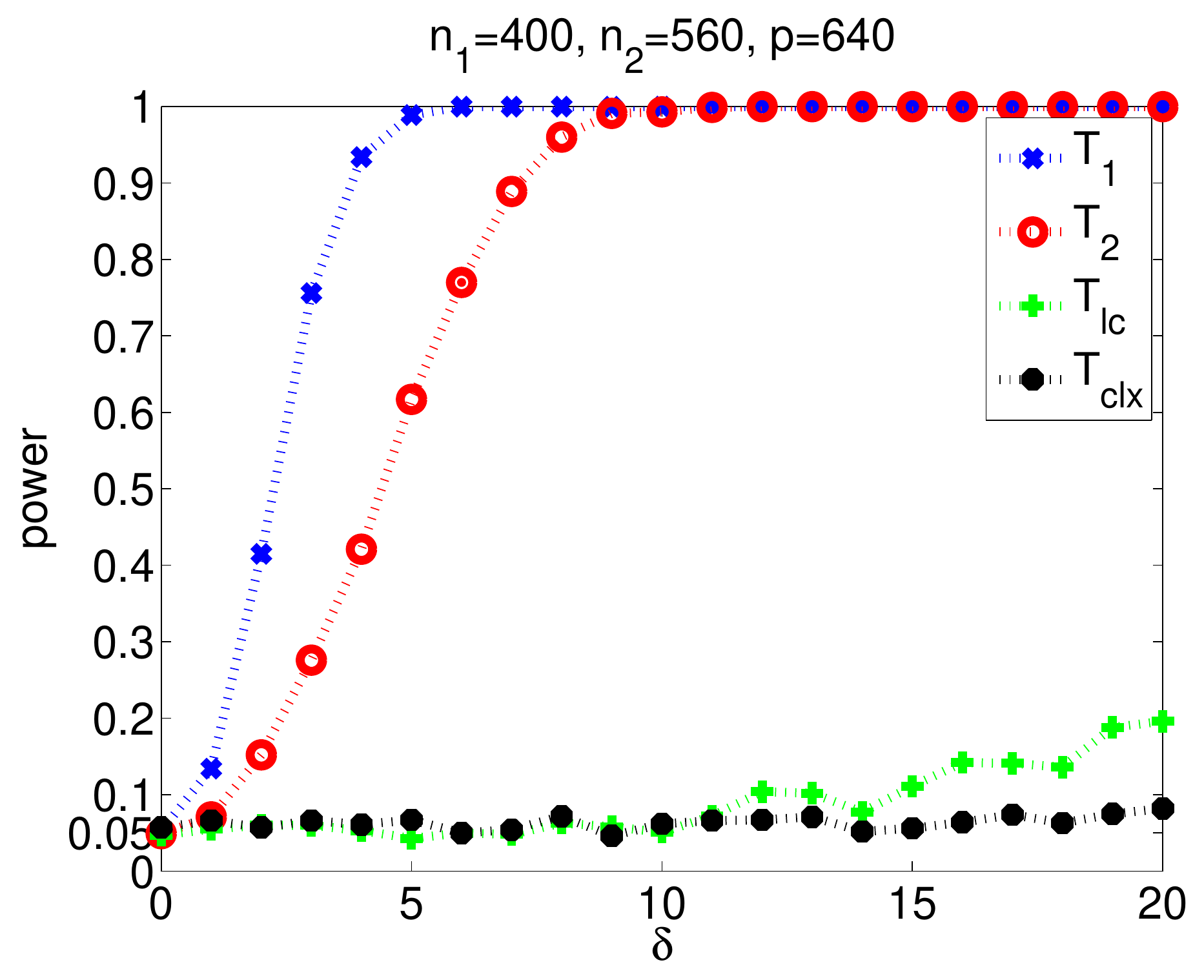}\includegraphics[height=4.7cm, width=7.4cm]{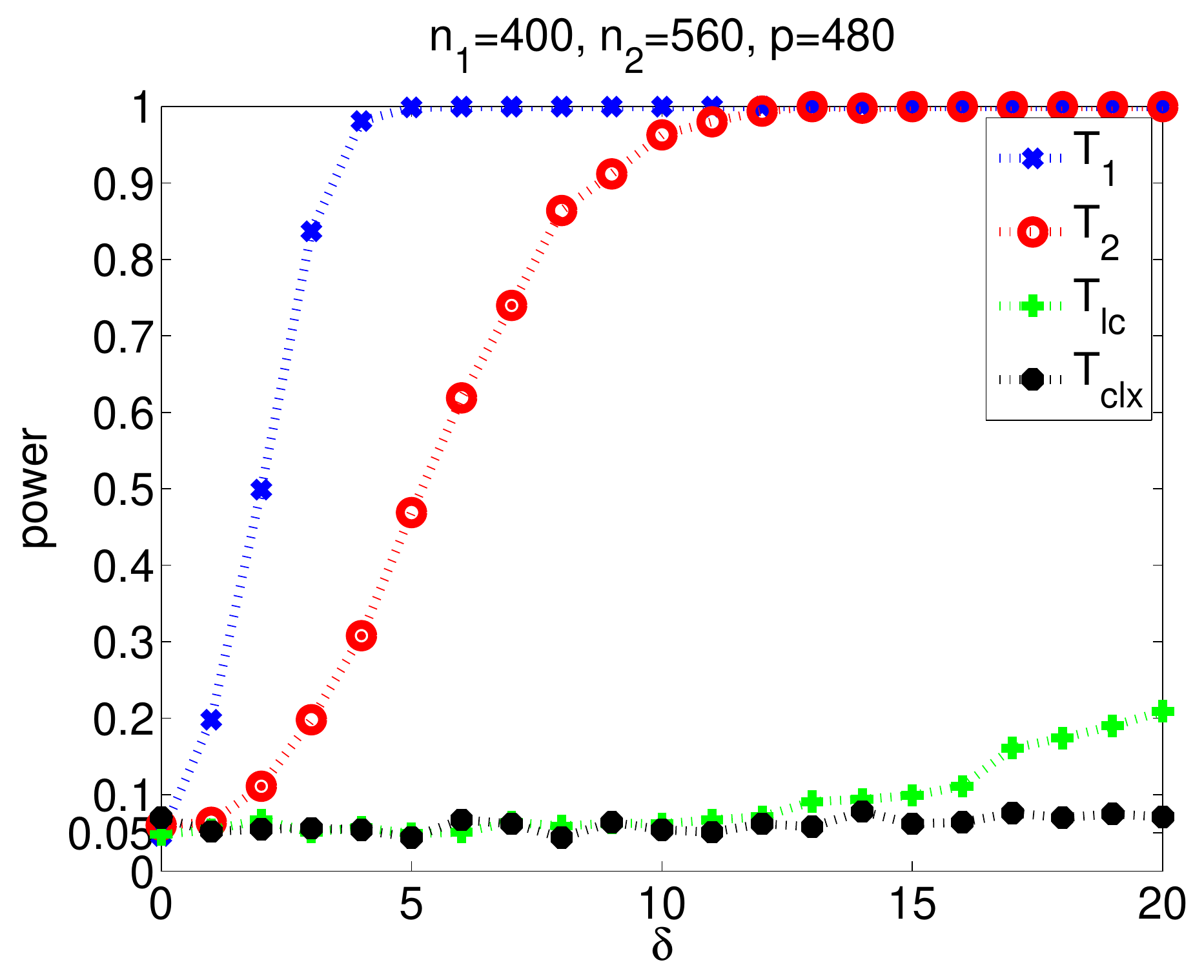}
        \includegraphics[height=4.7cm, width=7.4cm]{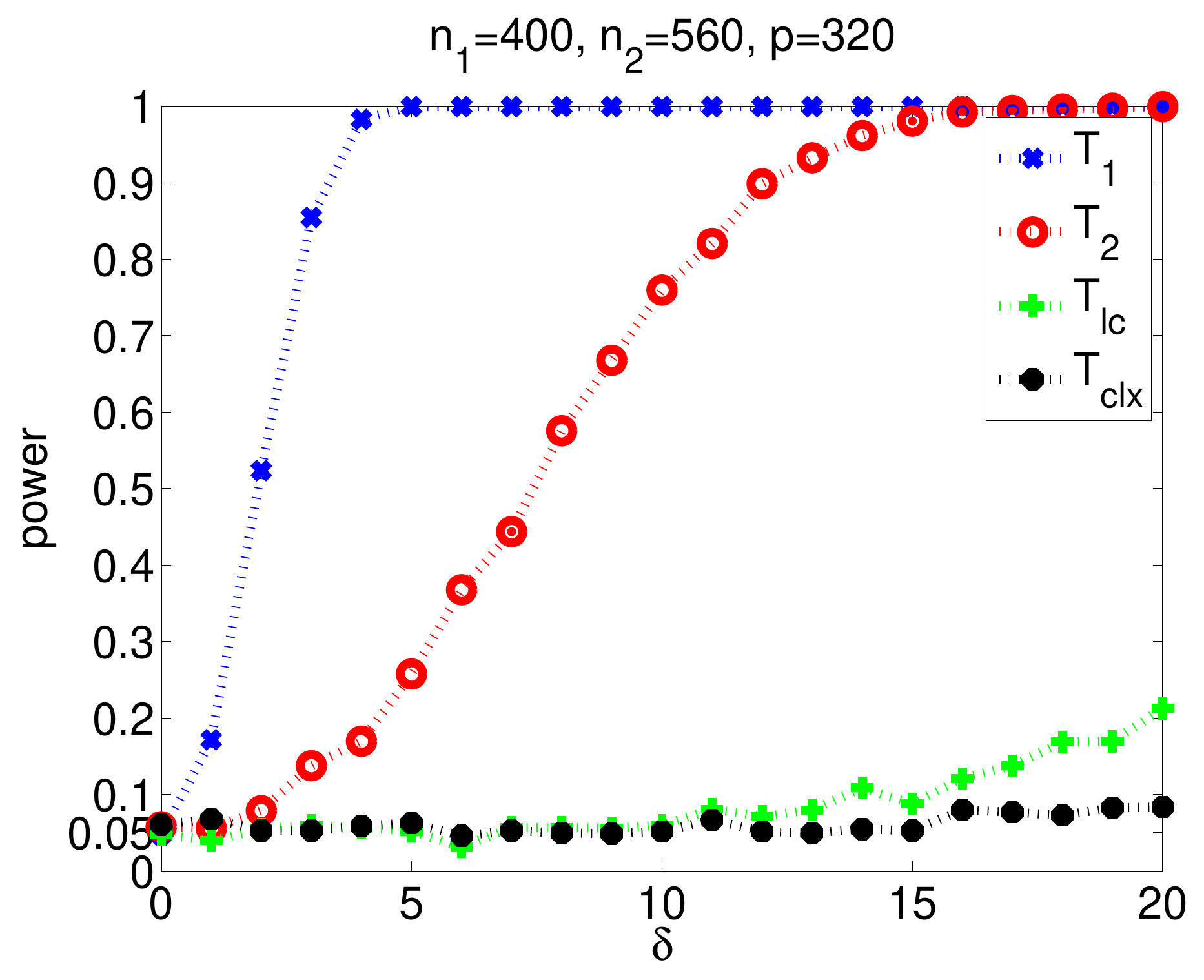}\includegraphics[height=4.7cm, width=7.4cm]{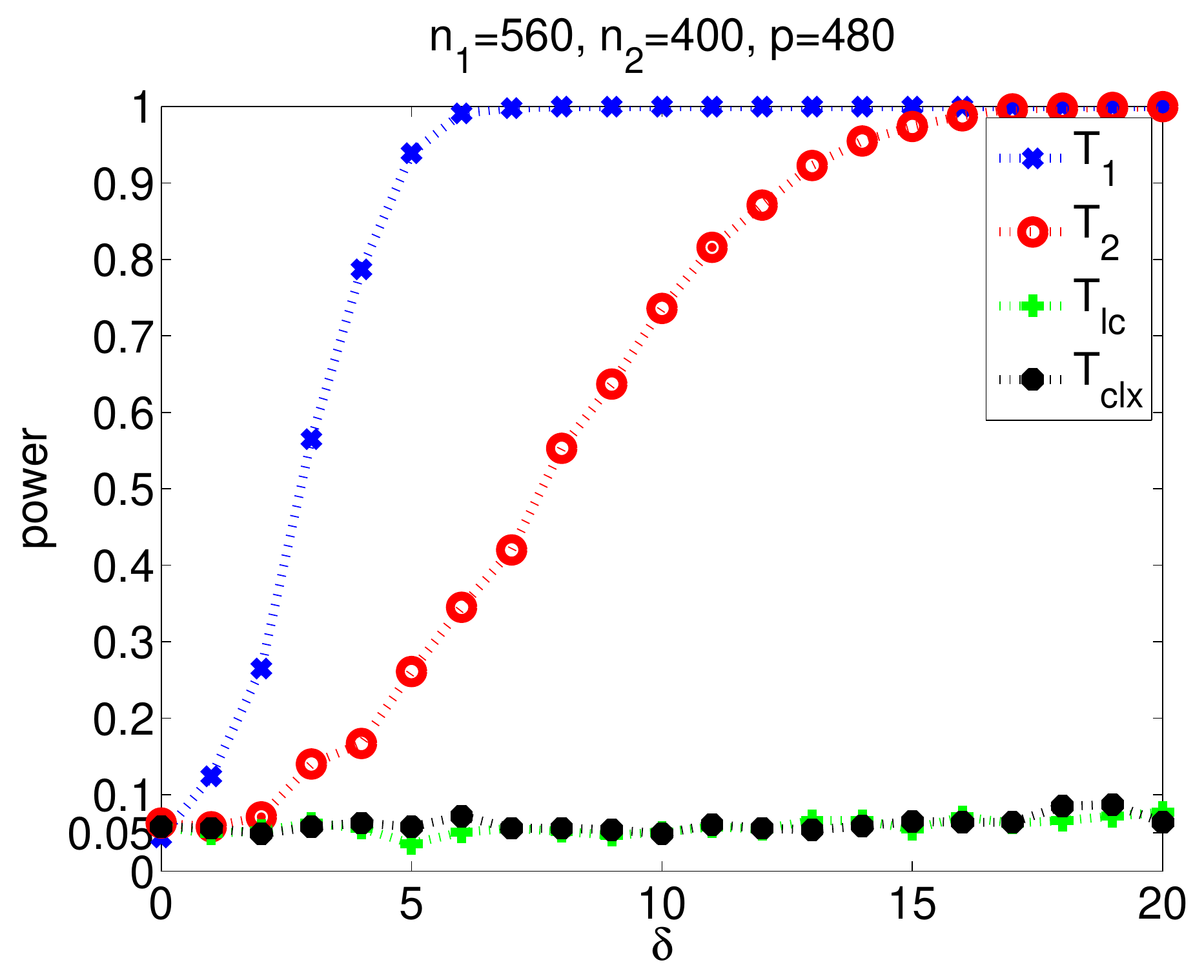}
    \caption{Scatter diagram of the empirical power for $T_1$, $T_2$, $T_{lc}$ and $T_{clx}$ based on Model 3 under the uniform assumption. }
    \label{unifMo3}
    \end{center}
\end{figure}

\begin{figure}[htbp]
	\begin{center}
		\includegraphics[height=4.7cm, width=7.4cm]{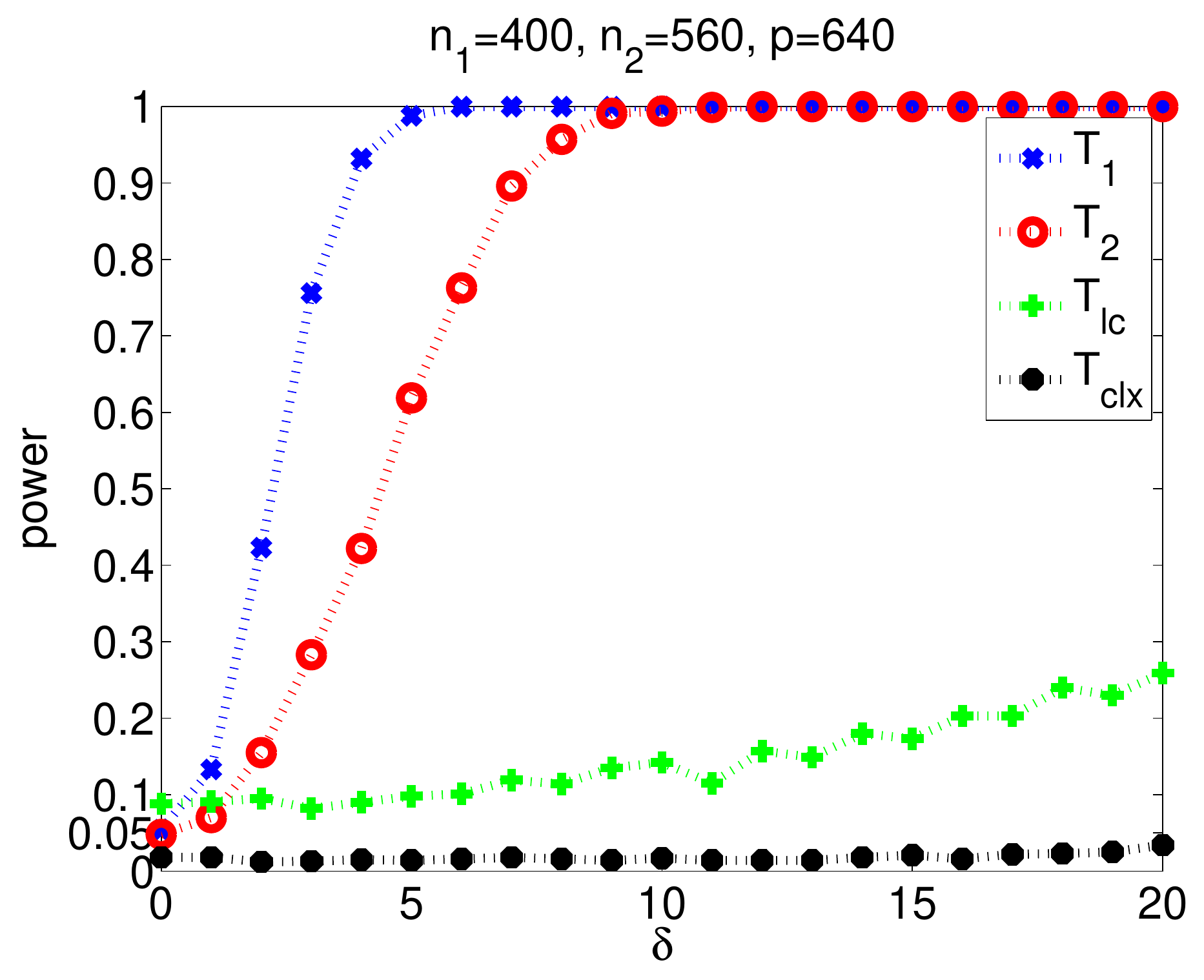}\includegraphics[height=4.7cm, width=7.4cm]{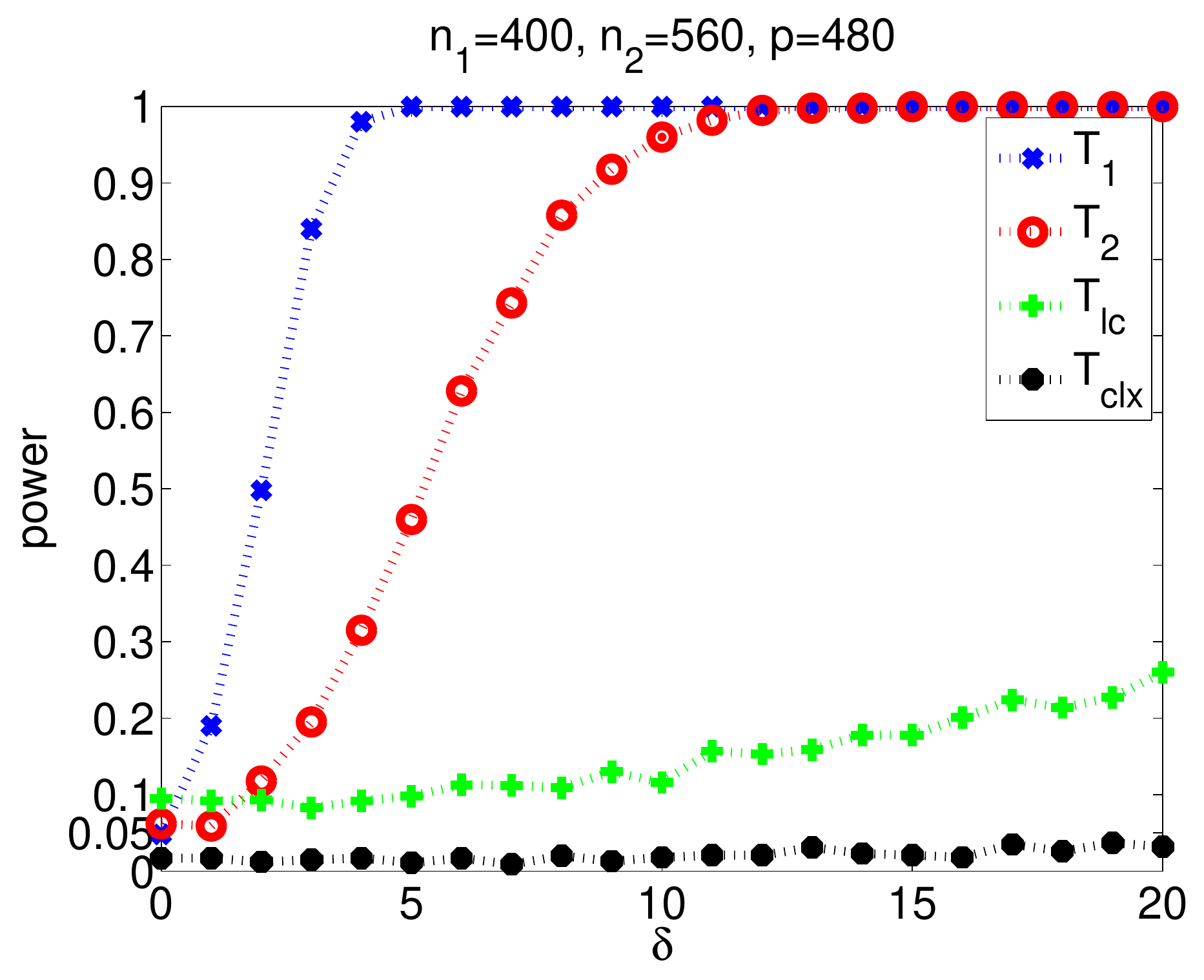}
        \includegraphics[height=4.7cm, width=7.4cm]{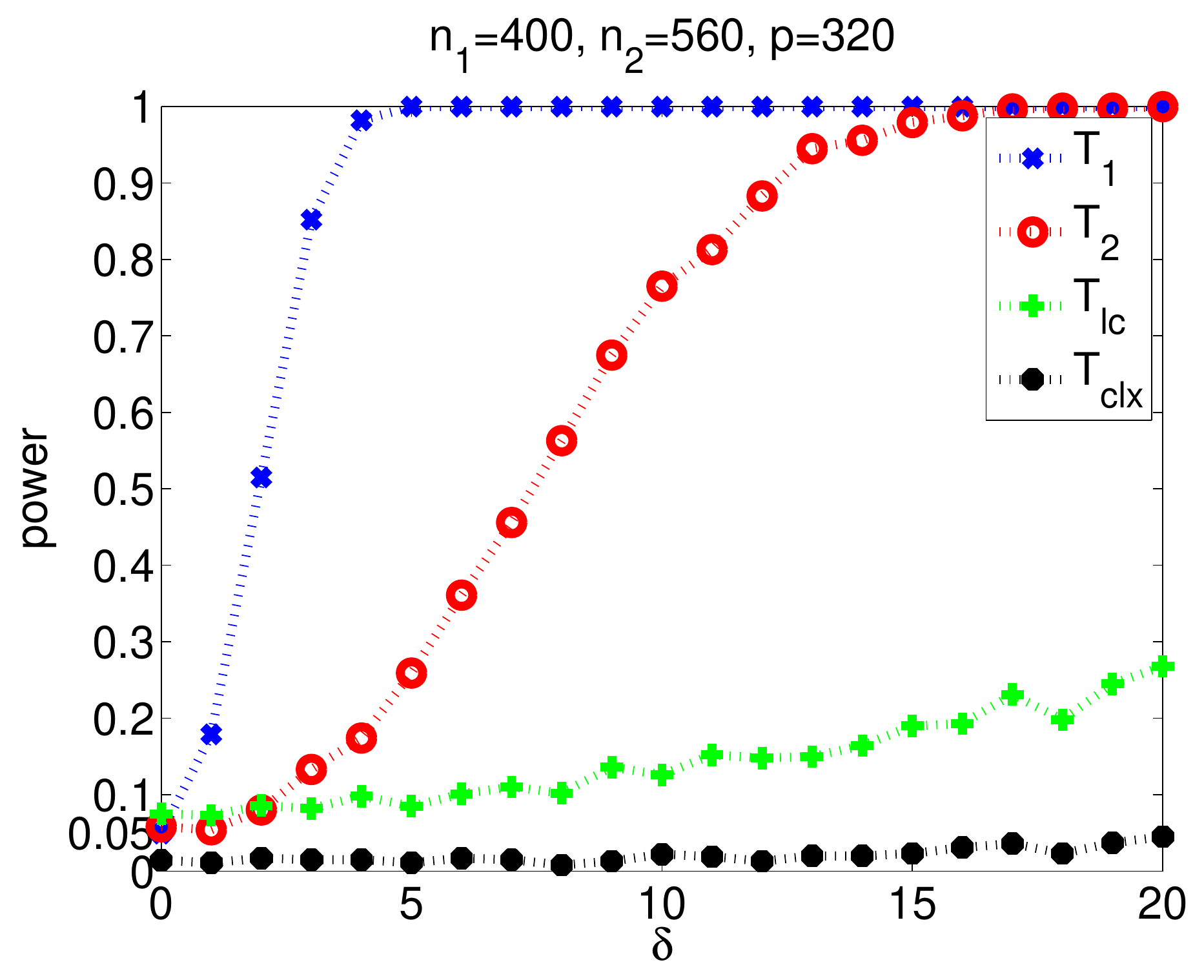}\includegraphics[height=4.7cm, width=7.4cm]{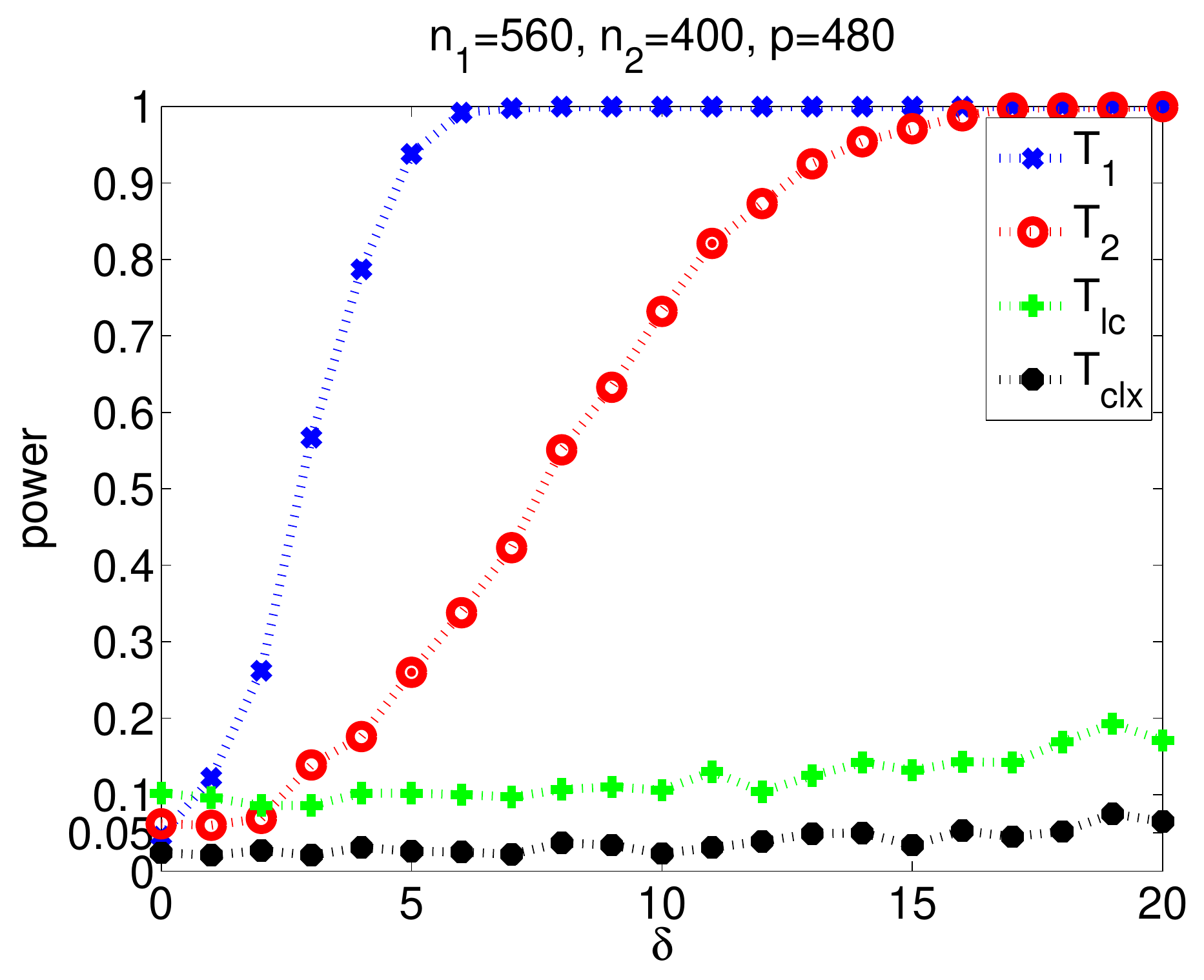}
    \caption{Scatter diagram of the empirical power for $T_1$, $T_2$, $T_{lc}$ and $T_{clx}$ based on Model 4 under the uniform assumption. }
    \label{unifMo4}
\end{center}
\end{figure}

Because the proposed statistics are also applicable under nonnormal conditions, we also simulate $x_{j}^{(l)}$ with the $p-$dimensional uniform distribution Unif$_{p}(-\sqrt{3},\sqrt{3})$; that is, all elements of $x_{j}^{(l)}$ are independently generated from the uniform distribution Unif$(-\sqrt{3},\sqrt{3})$.
The results of the four models for $\Sigma_{l}$ are presented in Tables \ref{uniM1}-\ref{uniM4} and in Figures \ref{unifMo1}-\ref{unifMo4}, respectively.
From these results, we find that the performances of $T_1$  under uniform distribution condition are very similar to its performances under the normal distribution condition, and so is to the performances of $T_2$. Therefore, that coincides with our claim that the proposed statistics are universal.

\begin{table}[h]
\scriptsize
\begin{tabular}{ccccccccccccc}
\\ \hline
 \multirow{3}{6em}{($n_1$,$n_2$,$p$)\\$y_1>1$,$y_2>1$} &\multicolumn{3}{c}{(50,70,80)} &\multicolumn{3}{c}{(100,140,160)} &\multicolumn{3}{c}{(200,280,320)}&\multicolumn{3}{c}{(400,560,640)}\\ \cline{2-13}
                     &size&\multicolumn{2}{c}{power}
                     &size&\multicolumn{2}{c}{power}&size&\multicolumn{2}{c}{power}&size&\multicolumn{2}{c}
                     {power}\\\cline{2-13}
                        &$\delta$=0&$\delta$=5&$\delta$=10&$\delta$=0 &$\delta$=5&$\delta$=10&$\delta$=0&$\delta$=5&$\delta$=10&$\delta$=0&$\delta$=5&$\delta$=10 \\ \hline
		$T_1$               &0.043 & 0.991 & 1 & 0.054 &0.987 & 1 & 0.053 & 0.993& 1&0.052&0.990 &1      \\ \hline
		$T_2$               &0.055  &0.622 & 1 & 0.049& 0.623 & 0.997 & 0.058 &0.620& 0.995  &0.058&0.643 &0.996     \\ \hline
$T_{lc}$                    &0.058 & 0.101 & 0.490  & 0.058 & 0.058 & 0.227  & 0.032 & 0.049 & 0.094 &0.049&0.070 &0.064  \\ \hline
$T_{clx}$                   &0.122 & 0.180 & 0.567 & 0.091 & 0.108& 0.214 & 0.058 & 0.069 & 0.081 &0.067&0.065 &0.077 \\ \hline
		
		\multirow{3}{6.3em}{{\centering ($n_1,n_2,p$)\\$y_1>1$,$y_2<1$}} &\multicolumn{3}{c}{(50,70,60)} &\multicolumn{3}{c}{(100,140,120)} &\multicolumn{3}{c}{(200,280,240)}&\multicolumn{3}{c}{(400,560,480)}\\ \cline{2-13}
		&size&\multicolumn{2}{c}{power}
		&size&\multicolumn{2}{c}{power}&size&\multicolumn{2}{c}{power}&size&\multicolumn{2}{c}
		{power}\\\cline{2-13}
		&$\delta$=0&$\delta$=5&$\delta$=10&$\delta$=0 &$\delta$=5&$\delta$=10&$\delta$=0&$\delta$=5&$\delta$=10&$\delta$=0&$\delta$=5&$\delta$=10 \\ \hline
		$T_1$                  &0.059 & 0.996 & 1 & 0.047 & 0.998 & 1 & 0.047 & 0.998 & 1 &0.057&0.998    &1     \\ \hline
		$T_2$                  &0.041 & 0.454 & 0.975 & 0.048 & 0.462 & 0.971 & 0.042 & 0.460 & 0.973 &0.045&0.489&0.961    \\ \hline
$T_{lc}$  &	 0.055 & 0.099 & 0.518 & 0.048 & 0.062 & 0.200 & 0.054 & 0.072 & 0.113 &0.049&0.055&0.062 \\ \hline
$T_{clx}$ &	 0.120 & 0.161 & 0.587 & 0.076 & 0.121 & 0.230 & 0.078 & 0.083 & 0.122 &0.060&0.048&0.072 \\ \hline
		
		\multirow{3}{6.3em}{{\centering ($n_1,n_2,p$)\\$y_1<1$,$y_2>1$}} &\multicolumn{3}{c}{(70,50,60)} &\multicolumn{3}{c}{(140,100,120)} &\multicolumn{3}{c}{(280,200,240)}&\multicolumn{3}{c}{(560,400,480)}\\ \cline{2-13}
		&size&\multicolumn{2}{c}{power}
		&size&\multicolumn{2}{c}{power}&size&\multicolumn{2}{c}{power}&size&\multicolumn{2}{c}
		{power}\\\cline{2-13}
		&$\delta$=0&$\delta$=5&$\delta$=10&$\delta$=0 &$\delta$=5&$\delta$=10&$\delta$=0&$\delta$=5&$\delta$=10&$\delta$=0&$\delta$=5&$\delta$=10 \\ \hline
		$T_1$                   & 0.052 & 0.928 & 1 & 0.057 & 0.945 & 1 & 0.048 & 0.946 & 1 &0.053&0.931   &1     \\ \hline
		$T_2$                   & 0.048 & 0.224 & 0.572 & 0.048 & 0.251& 0.684  & 0.032 & 0.241 & 0.715 &0.060&0.241 &0.748 \\ \hline
$T_{lc}$ &	 0.054 & 0.050 &0.149  & 0.053 & 0.052 & 0.059 & 0.045 & 0.048& 0.064  &0.046	&0.032 &0.046		\\ \hline
$T_{clx}$& 0.120 & 0.167 & 0.419 & 0.105 & 0.110& 0.158  & 0.052 & 0.064 & 0.077 &0.075&0.045	&0.055 \\ \hline
		\multirow{3}{6.3em}{{\centering ($n_1,n_2,p$)\\$y_1<1$,$y_2<1$}} &\multicolumn{3}{c}{(50,70,40)} &\multicolumn{3}{c}{(100,140,80)} &\multicolumn{3}{c}{(200,280,160)}&\multicolumn{3}{c}{(400,560,320)} \\ \cline{2-13}
		&size&\multicolumn{2}{c}{power}
		&size&\multicolumn{2}{c}{power}&size&\multicolumn{2}{c}{power}&size&\multicolumn{2}{c}
		{power}\\\cline{2-13}
		&$\delta$=0&$\delta$=5&$\delta$=10&$\delta$=0 &$\delta$=5&$\delta$=10&$\delta$=0&$\delta$=5&$\delta$=10&$\delta$=0&$\delta$=5&$\delta$=10  \\ \hline
		$T_1$               &0.047 & 0.994 & 1 & 0.050 & 0.999 & 1 & 0.047 & 1 & 1 &0.045	&0.999  	&  1	      \\ \hline
		$T_2$              & 0.055 & 0.279 & 0.860 & 0.044 & 0.269  & 0.810 & 0.052 & 0.281& 0.779  &0.046	&0.243  &0.755	  \\ \hline
$T_{lc}$  &	 0.060 & 0.098 & 0.497 & 0.050 & 0.077 & 0.209  & 0.040 & 0.066 & 0.094 &0.062&0.067&0.072 \\ \hline
$T_{clx}$  &0.106 & 0.169 & 0.563 & 0.072 & 0.106 & 0.237 & 0.056 & 0.069 & 0.100 &0.049&0.070&0.065 \\ \hline
\end{tabular}
\caption{Empirical size and power for 1000 repeated simulations to compare $T_1$, $T_2$, $T_{lc}$ and $T_{clx}$ based on Model 1 under the uniform assumption.}
\label{uniM1}
\end{table}

\begin{table}[h]
\scriptsize
\begin{tabular}{ccccccccccccc}
\\ \hline
 \multirow{3}{6em}{($n_1$,$n_2$,$p$)\\$y_1>1$,$y_2>1$} &\multicolumn{3}{c}{(50,70,80)} &\multicolumn{3}{c}{(100,140,160)} &\multicolumn{3}{c}{(200,280,320)}&\multicolumn{3}{c}{{(400,560,640)}}\\ \cline{2-13}
                     &size&\multicolumn{2}{c}{power}
                     &size&\multicolumn{2}{c}{power}&size&\multicolumn{2}{c}{power}&size&\multicolumn{2}{c}
                     {power}\\\cline{2-13}
                        &$\delta$=0&$\delta$=5&$\delta$=10&$\delta$=0 &$\delta$=5&$\delta$=10&$\delta$=0&$\delta$=5&$\delta$=10&$\delta$=0&$\delta$=5&$\delta$=10 \\ \hline
		$T_1$                &0.048 & 0.992 & 1 & 0.061 & 0.989 & 1 & 0.060 & 0.996 & 1 &0.046&0.991    &1     \\ \hline
		$T_2$                &0.048 & 0.623 & 0.999 & 0.053 & 0.645 & 0.997 & 0.044 & 0.591 & 1 &0.052&0.604&0.992     \\ \hline
$T_{lc}$ &0.013 & 0.112 & 0.411 & 0.013 & 0.064 & 0.211 & 0.016 & 0.033 & 0.096 &0.011&0.017&0.062 \\ \hline
$T_{clx}$&0.128 & 0.182 & 0.562 & 0.090 & 0.100 & 0.226 & 0.072 & 0.064 & 0.089  &0.076&0.069&0.068 \\ \hline
		
		\multirow{3}{6.3em}{{\centering ($n_1,n_2,p$)\\$y_1>1$,$y_2<1$}} &\multicolumn{3}{c}{(50,70,60)} &\multicolumn{3}{c}{(100,140,120)} &\multicolumn{3}{c}{(200,280,240)}&\multicolumn{3}{c}{(400,560,480)} \\ \cline{2-13}
		&size&\multicolumn{2}{c}{power}
		&size&\multicolumn{2}{c}{power}&size&\multicolumn{2}{c}{power}&size&\multicolumn{2}{c}
		{power}\\\cline{2-13}
		&$\delta$=0&$\delta$=5&$\delta$=10&$\delta$=0 &$\delta$=5&$\delta$=10&$\delta$=0&$\delta$=5&$\delta$=10&$\delta$=0&$\delta$=5&$\delta$=10 \\ \hline
		$T_1$                 &0.041 & 0.996 & 1 & 0.049 & 0.999 & 1 & 0.049 & 0.998 & 1 &0.045&0.998    &1     \\ \hline
		$T_2$                 &0.045 & 0.459 & 0.977 & 0.045 & 0.487 & 0.969 & 0.050 & 0.463 & 0.963 &0.033&0.481&0.968     \\ \hline
$T_{lc}$  &0.014 & 0.115 & 0.409 & 0.006 & 0.060 & 0.237 & 0.015 & 0.024 & 0.140 &0.018&0.031&0.062  \\ \hline
$T_{clx}$ &0.108 & 0.152 & 0.578 & 0.077 & 0.102 & 0.228 & 0.077 & 0.066 & 0.107 &0.056&0.059&0.060     \\ \hline
		\multirow{3}{6.3em}{{\centering ($n_1,n_2,p$)\\$y_1<1$,$y_2>1$}} &\multicolumn{3}{c}{(70,50,60)} &\multicolumn{3}{c}{(140,100,120)} &\multicolumn{3}{c}{(280,200,240)}&\multicolumn{3}{c}{(560,400,480)} \\ \cline{2-13}
		&size&\multicolumn{2}{c}{power}
		&size&\multicolumn{2}{c}{power}&size&\multicolumn{2}{c}{power}&size&\multicolumn{2}{c}
		{power}\\\cline{2-13}
		&$\delta$=0&$\delta$=5&$\delta$=10&$\delta$=0 &$\delta$=5&$\delta$=10&$\delta$=0&$\delta$=5&$\delta$=10&$\delta$=0&$\delta$=5&$\delta$=10  \\ \hline
		$T_1$            & 0.047 & 0.941 & 1 & 0.042 & 0.941 & 1 & 0.051 & 0.936 & 1 &0.055	&0.931  	&1  	 \\ \hline
		$T_2$            & 0.048 & 0.220 & 0.560 & 0.046 & 0.241 & 0.666 & 0.043 & 0.269 & 0.711 &0.057	&0.258	&0.756	   \\ \hline
$T_{lc}$ & 0.015 & 0.055 & 0.186 & 0.013 & 0.028 & 0.106 & 0.019 & 0.021 & 0.066 &0.015	&0.022	&0.041	 \\ \hline
$T_{clx}$& 0.122 & 0.168 & 0.392 & 0.103 & 0.099 & 0.158 & 0.073 & 0.059 & 0.076 &0.066	&0.060	&0.053	     \\ \hline
		\multirow{3}{6.3em}{{\centering ($n_1,n_2,p$)\\$y_1<1$,$y_2<1$}} &\multicolumn{3}{c}{(50,70,40)} &\multicolumn{3}{c}{(100,140,80)} &\multicolumn{3}{c}{(200,280,160)}&\multicolumn{3}{c}{(400,560,320)} \\ \cline{2-13}
		&size&\multicolumn{2}{c}{power}
		&size&\multicolumn{2}{c}{power}&size&\multicolumn{2}{c}{power}&size&\multicolumn{2}{c}
		{power}\\\cline{2-13}
		&$\delta$=0&$\delta$=5&$\delta$=10&$\delta$=0 &$\delta$=5&$\delta$=10&$\delta$=0&$\delta$=5&$\delta$=10&$\delta$=0&$\delta$=5&$\delta$=10 \\ \hline
		$T_1$                & 0.048 & 0.996 & 1 & 0.058 & 0.998 & 1 & 0.039 & 0.999 & 1 &0.042&1    &   1   \\ \hline
		$T_2$                & 0.045 & 0.283 & 0.864 & 0.050 & 0.264 & 0.812 & 0.046 & 0.282 & 0.767 &0.042&0.252&0.788    \\ \hline
$T_{lc}$ & 0.012 & 0.115 & 0.388 & 0.018 & 0.067 & 0.220 & 0.008 & 0.037 & 0.112 &0.008&0.019&0.066 \\ \hline
$T_{clx}$& 0.122 & 0.171 & 0.589 & 0.077 & 0.096 & 0.255 & 0.072 & 0.055 & 0.107 &0.059&0.066&0.070     \\ \hline
		
\end{tabular}
\caption{Empirical size and power for 1000 repeated simulations to compare $T_1$, $T_2$, $T_{lc}$ and $T_{clx}$ based on Model 2 under the uniform assumption.}
\label{uniM2}
\end{table}

\begin{table}[h]
\scriptsize
\begin{tabular}{ccccccccccccc}
 \\ \hline
 \multirow{3}{6em}{($n_1$,$n_2$,$p$)\\$y_1>1$,$y_2>1$}&\multicolumn{3}{c}{(50,70,80)} &\multicolumn{3}{c}{(100,140,160)} &\multicolumn{3}{c}{(200,280,320)}&\multicolumn{3}{c}{{(400,560,640)}} \\ \cline{2-13}
                     &size&\multicolumn{2}{c}{power}
                     &size&\multicolumn{2}{c}{power}&size&\multicolumn{2}{c}{power}&size&\multicolumn{2}{c}
                     {power}\\\cline{2-13}
                        &$\delta$=0&$\delta$=5&$\delta$=10&$\delta$=0 &$\delta$=5&$\delta$=10&$\delta$=0&$\delta$=5&$\delta$=10&$\delta$=0&$\delta$=5&$\delta$=10 \\ \hline
$T_1$   & 0.040 & 0.991 & 1 & 0.053 & 0.990 & 1 & 0.046 & 0.986 & 1 &0.060&0.992   &1      \\ \hline
$T_2$   & 0.053 & 0.621 & 1 & 0.056 & 0.643 & 0.996 & 0.054 & 0.620 & 0.996 &0.055&0.602&0.996     \\ \hline
$T_{lc}$& 0.053 & 0.113 & 0.475 & 0.060 & 0.077 & 0.245 & 0.056 & 0.058 & 0.098 &0.042&0.059&0.064     \\ \hline
$T_{clx}$ & 0.096 & 0.095 & 0.222 & 0.070 & 0.088 & 0.108 & 0.063 & 0.062 & 0.071 &0.044&0.048&0.066 \\ \hline
		
		\multirow{3}{6.3em}{{\centering ($n_1,n_2,p$)\\$y_1>1$,$y_2<1$}} &\multicolumn{3}{c}{(50,70,60)} &\multicolumn{3}{c}{(100,140,120)} &\multicolumn{3}{c}{(200,280,240)}&\multicolumn{3}{c}{{(400,560,480)}} \\ \cline{2-13}
		&size&\multicolumn{2}{c}{power}
		&size&\multicolumn{2}{c}{power}&size&\multicolumn{2}{c}{power}&size&\multicolumn{2}{c}
		{power}\\\cline{2-13}
		&$\delta$=0&$\delta$=5&$\delta$=10&$\delta$=0 &$\delta$=5&$\delta$=10&$\delta$=0&$\delta$=5&$\delta$=10&$\delta$=0&$\delta$=5&$\delta$=10 \\ \hline
$T_1$   & 0.063 & 0.996 & 1 & 0.063 & 0.999 & 1 & 0.042 &0.999 & 1 &0.043&0.998    &1     \\ \hline
$T_2$   & 0.056 & 0.469 & 0.982 & 0.045 & 0.481 & 0.974 & 0.051 & 0.478 & 0.978 &0.051&0.460&0.962     \\ \hline
$T_{lc}$& 0.063 & 0.124 & 0.483 & 0.068 & 0.077 & 0.192 & 0.043 & 0.054 & 0.104 &0.058&0.054&0.063 \\ \hline
$T_{clx}$& 0.082 & 0.109 & 0.240 & 0.088 & 0.076 & 0.098 & 0.060 & 0.049 & 0.063 &0.064&0.047&0.064  \\ \hline
		\multirow{3}{6.3em}{{\centering ($n_1,n_2,p$)\\$y_1<1$,$y_2>1$}} &\multicolumn{3}{c}{(70,50,60)} &\multicolumn{3}{c}{(140,100,120)} &\multicolumn{3}{c}{(280,200,240)}&\multicolumn{3}{c}{(560,400,480)} \\ \cline{2-13}
		&size&\multicolumn{2}{c}{power}
		&size&\multicolumn{2}{c}{power}&size&\multicolumn{2}{c}{power}&size&\multicolumn{2}{c}
		{power}\\\cline{2-13}
		&$\delta$=0&$\delta$=5&$\delta$=10&$\delta$=0 &$\delta$=5&$\delta$=10&$\delta$=0&$\delta$=5&$\delta$=10&$\delta$=0&$\delta$=5&$\delta$=10  \\ \hline

$T_1$   & 0.043 & 0.937 & 1 & 0.038 & 0.930 & 1 & 0.044 &0.946 & 1 &0.052	&0.940 &1  	 \\ \hline
$T_2$   & 0.036 & 0.228 & 0.596 & 0.038 & 0.218 & 0.664 & 0.040 & 0.238 & 0.716 &0.039	&0.277	&0.720  \\ \hline
$T_{lc}$& 0.058 & 0.065 & 0.172 & 0.058 & 0.058 & 0.095 & 0.051 & 0.053 & 0.068 &0.051	&0.064	&0.051	\\ \hline
$T_{clx}$& 0.091 & 0.113 & 0.241 & 0.069 & 0.073 & 0.107 & 0.050 & 0.055 & 0.061 &0.051	&0.065	&0.051   \\ \hline
		\multirow{3}{6.3em}{{\centering ($n_1,n_2,p$)\\$y_1<1$,$y_2<1$}} &\multicolumn{3}{c}{(50,70,40)} &\multicolumn{3}{c}{(100,140,80)} &\multicolumn{3}{c}{(200,280,160)}&\multicolumn{3}{c}{(400,560,320)} \\ \cline{2-13}
		&size&\multicolumn{2}{c}{power}
		&size&\multicolumn{2}{c}{power}&size&\multicolumn{2}{c}{power}&size&\multicolumn{2}{c}
		{power}\\\cline{2-13}
		&$\delta$=0&$\delta$=5&$\delta$=10&$\delta$=0 &$\delta$=5&$\delta$=10&$\delta$=0&$\delta$=5&$\delta$=10&$\delta$=0&$\delta$=5&$\delta$=10 \\ \hline

$T_1$   & 0.065 & 0.995 & 1 & 0.045 & 0.997 & 1 & 0.049 & 1 & 1 &0.042&1    &1     \\ \hline
$T_2$   & 0.042 & 0.293 & 0.854 & 0.054 & 0.287 & 0.774 & 0.051 & 0.281 & 0.777 &0.042&0.251&0.757     \\ \hline
$T_{lc}$&0.054 & 0.120 & 0.445 & 0.059 & 0.087 & 0.193 & 0.055 & 0.049 & 0.105 &0.040&0.059&0.059 \\ \hline
$T_{clx}$ & 0.103 & 0.094 & 0.235 & 0.063 & 0.072 & 0.093 & 0.066 & 0.058 & 0.075 &0.061&0.053&0.054 \\ \hline
		
\end{tabular}
\caption{Empirical size and power for 1000 repeated simulations to compare $T_1$, $T_2$, $T_{lc}$ and $T_{clx}$ based on Model 3 under the uniform assumption.}
\label{uniM3}
\end{table}

\begin{table}[h]
\scriptsize
\begin{tabular}{ccccccccccccc}
 \\ \hline
 \multirow{3}{6em}{($n_1$,$n_2$,$p$)\\$y_1>1$,$y_2>1$}&\multicolumn{3}{c}{(50,70,80)} &\multicolumn{3}{c}{(100,140,160)} &\multicolumn{3}{c}{(200,280,320)}&\multicolumn{3}{c}{{(400,560,640)}} \\ \cline{2-13}
                     &size&\multicolumn{2}{c}{power}
                     &size&\multicolumn{2}{c}{power}&size&\multicolumn{2}{c}{power}&size&\multicolumn{2}{c}
                     {power}\\\cline{2-13}
                        &$\delta$=0&$\delta$=5&$\delta$=10&$\delta$=0 &$\delta$=5&$\delta$=10&$\delta$=0&$\delta$=5&$\delta$=10&$\delta$=0&$\delta$=5&$\delta$=10 \\ \hline
		$T_1$                     &0.060 & 0.992 & 1 & 0.055 & 0.988 & 1 & 0.066 &0.995 & 1 &0.050&0.994 &  1      \\ \hline
		$T_2$             & 0.057 & 0.616 & 0.998 & 0.049 & 0.606 & 0.996 & 0.053 & 0.652 & 0.993  &0.051&0.603&  0.997     \\ \hline
$T_{lc}$ & 0.082 & 0.156 & 0.367 & 0.096 & 0.122 & 0.251 & 0.087 & 0.122 & 0.148 &0.073&0.098&0.516 \\ \hline
$T_{clx}$& 0.041 & 0.030 & 0.065 & 0.035 & 0.023 & 0.036 & 0.022 & 0.013 & 0.021  &0.065&0.011&0.216 \\ \hline
		
		\multirow{3}{6.3em}{{\centering ($n_1,n_2,p$)\\$y_1>1$,$y_2<1$}} &\multicolumn{3}{c}{(50,70,60)} &\multicolumn{3}{c}{(100,140,120)} &\multicolumn{3}{c}{(200,280,240)}&\multicolumn{3}{c}{(400,560,480)} \\ \cline{2-13}
		&size&\multicolumn{2}{c}{power}
		&size&\multicolumn{2}{c}{power}&size&\multicolumn{2}{c}{power}&size&\multicolumn{2}{c}
		{power}\\\cline{2-13}
		&$\delta$=0&$\delta$=5&$\delta$=10&$\delta$=0 &$\delta$=5&$\delta$=10&$\delta$=0&$\delta$=5&$\delta$=10&$\delta$=0&$\delta$=5&$\delta$=10 \\ \hline
		$T_1$                     & 0.051 & 0.996 & 1 & 0.057 & 1 & 1 & 0.050 & 0.999 & 1  &0.052& 0.999  &  1    \\ \hline
		$T_2$             & 0.053 & 0.464 & 0.979 & 0.047 & 0.444 & 0.978 & 0.053 & 0.473 & 0.968 &0.053&0.438&0.974      \\ \hline
$T_{lc}$ & 0.082 & 0.175 & 0.363 & 0.076 & 0.121 & 0.259 & 0.083 & 0.104 & 0.181  &0.094&0.121&0.508 \\ \hline
$T_{clx}$& 0.037 & 0.035 & 0.083 & 0.036 & 0.024 & 0.044 & 0.027 & 0.012 & 0.026 &0.058&0.014&0.251 \\ \hline
		\multirow{3}{6.3em}{{\centering ($n_1,n_2,p$)\\$y_1<1$,$y_2>1$}} &\multicolumn{3}{c}{(70,50,60)} &\multicolumn{3}{c}{(140,100,120)} &\multicolumn{3}{c}{(280,200,240)}&\multicolumn{3}{c}{(560,400,480)} \\ \cline{2-13}
		&size&\multicolumn{2}{c}{power}
		&size&\multicolumn{2}{c}{power}&size&\multicolumn{2}{c}{power}&size&\multicolumn{2}{c}
		{power}\\\cline{2-13}
		&$\delta$=0&$\delta$=5&$\delta$=10&$\delta$=0 &$\delta$=5&$\delta$=10&$\delta$=0&$\delta$=5&$\delta$=10&$\delta$=0&$\delta$=5&$\delta$=10  \\ \hline
		$T_1 $                   &0.046 &0.935 & 1 & 0.047 &0.939 & 1 & 0.048 &0.941 & 1  &0.054	&0.946 &1  	 \\ \hline
		$T_2$            &0.052 & 0.224 & 0.605 & 0.053 & 0.248 & 0.671 & 0.044 & 0.223 & 0.726 &0.054	&0.268	&0.906	     \\ \hline
$T_{lc}$ &0.075 & 0.113 & 0.222 & 0.081 & 0.108 & 0.184 & 0.074 & 0.099 & 0.132  &0.077	&0.094	&0.310	\\ \hline
$T_{clx}$&0.029 & 0.075 & 0.138 & 0.030 & 0.050 & 0.084 & 0.012 & 0.035 & 0.055 &0.081	&0.019	&0.282	\\ \hline
		\multirow{3}{6.3em}{{\centering ($n_1,n_2,p$)\\$y_1<1$,$y_2<1$}} &\multicolumn{3}{c}{(50,70,40)} &\multicolumn{3}{c}{(100,140,80)} &\multicolumn{3}{c}{(200,280,160)}&\multicolumn{3}{c}{(400,560,320)} \\ \cline{2-13}
		&size&\multicolumn{2}{c}{power}
		&size&\multicolumn{2}{c}{power}&size&\multicolumn{2}{c}{power}&size&\multicolumn{2}{c}
		{power}\\\cline{2-13}
		&$\delta$=0&$\delta$=5&$\delta$=10&$\delta$=0 &$\delta$=5&$\delta$=10&$\delta$=0&$\delta$=5&$\delta$=10&$\delta$=0&$\delta$=5&$\delta$=10 \\ \hline
		$T_1$                     & 0.062 & 0.995 & 1 & 0.045 & 0.998 & 1 & 0.052 & 1 & 1 &0.053&0.999&0.999      \\ \hline
		$T_2$             & 0.063 & 0.287 & 0.861 & 0.066 & 0.274 & 0.802 & 0.046 & 0.258 & 0.789 &0.047&0.261&0.452 \\ \hline
$T_{lc}$ & 0.088 & 0.153 & 0.381 & 0.076 & 0.129 & 0.229 & 0.087 & 0.096 & 0.157 &0.088&0.102&0.536 \\ \hline
$T_{clx}$& 0.048 & 0.049 & 0.109 & 0.026 & 0.026 & 0.041 & 0.024 & 0.020 & 0.024 &0.063&0.015&0.273 \\ \hline
		
\end{tabular}
\caption{Empirical size and power for 1000 repeated simulations to compare $T_1$, $T_2$, $T_{lc}$ and $T_{clx}$ based on Model 4 under the uniform assumption.}
\label{uniM4}
\end{table}

\begin{figure}[htbp]
	\begin{center}
		\includegraphics[width=7.4cm]{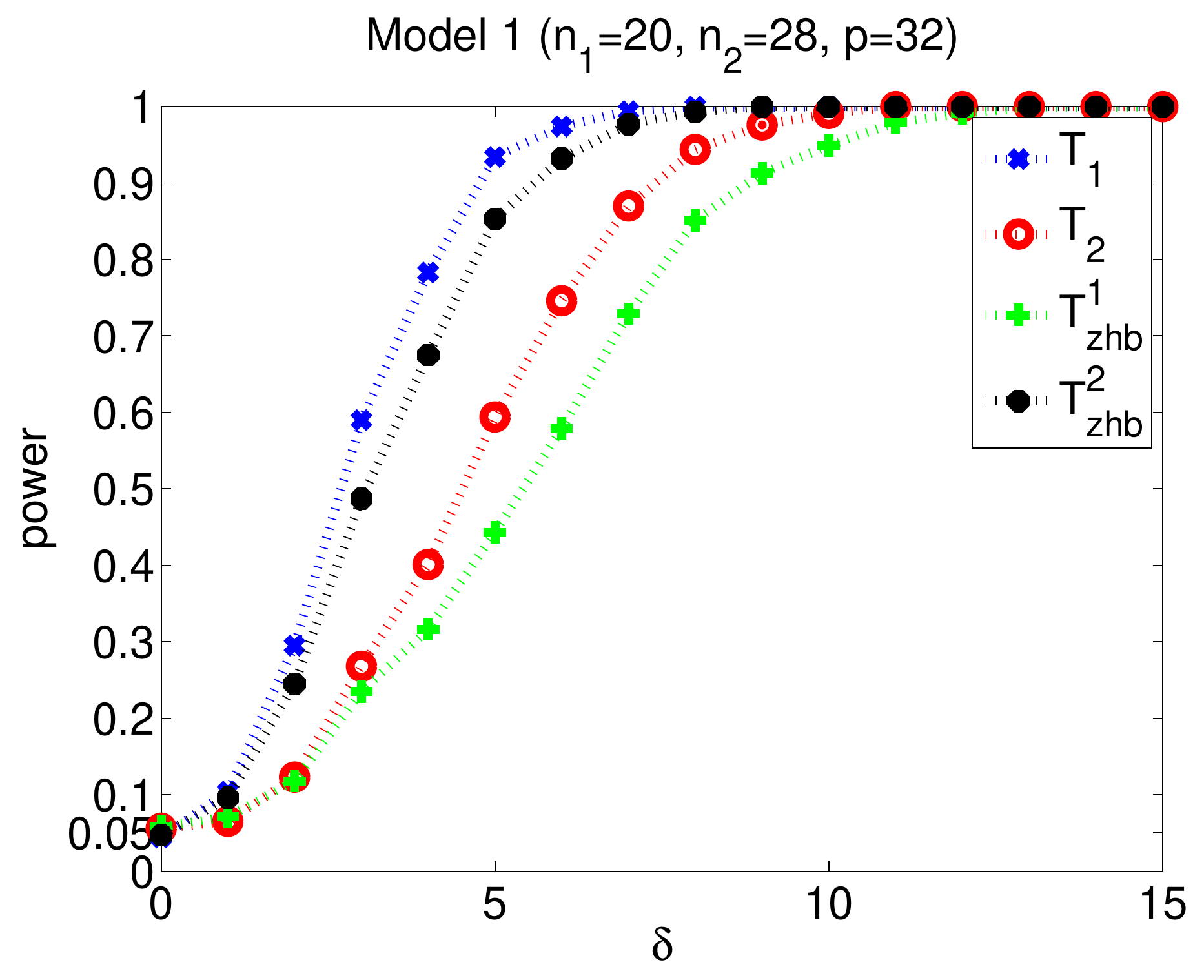}
        \includegraphics[width=7.4cm]{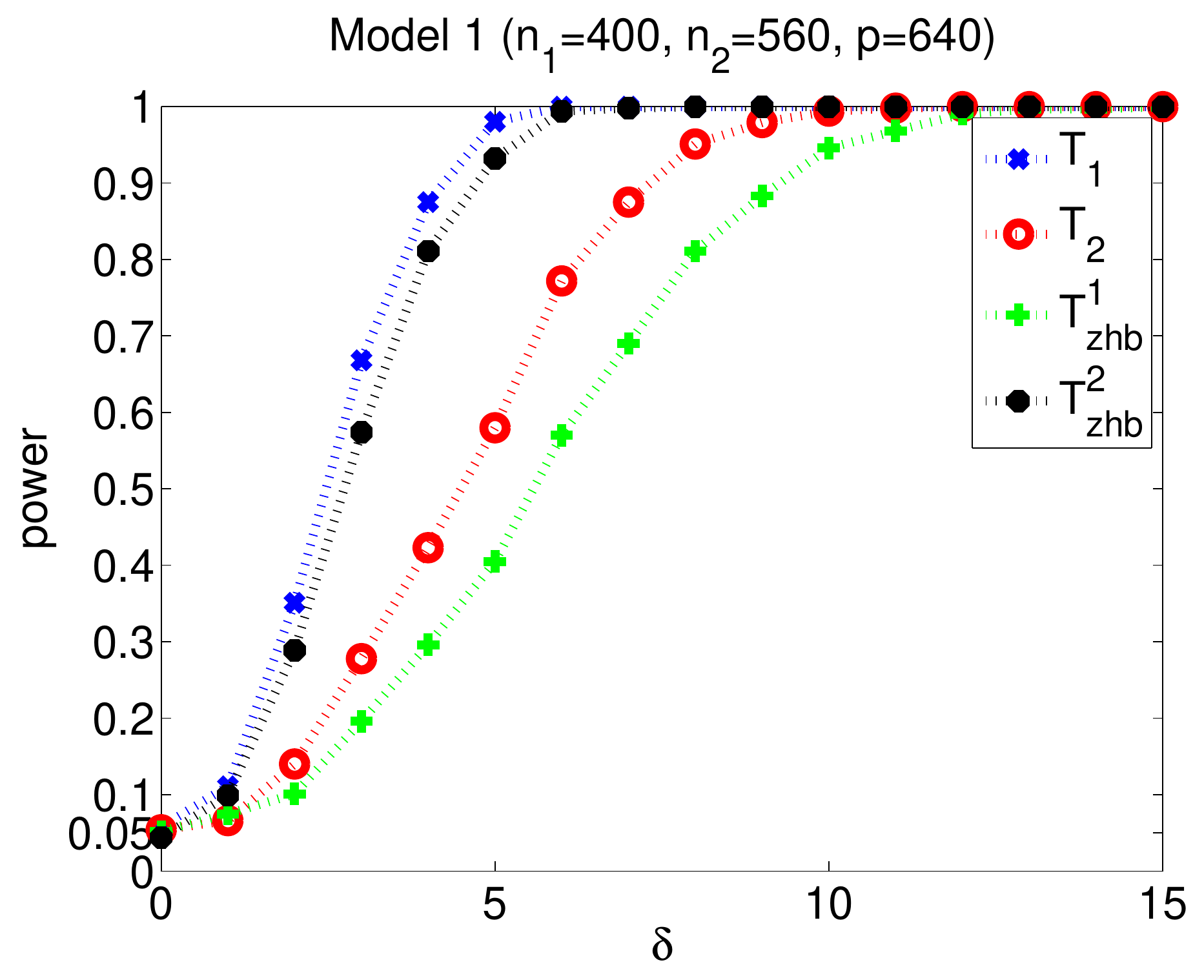}
    \caption{The empirical power for $T_1$, $T_2$, $T_{zhb}^{1}$ and $T_{zhb}^{2}$ under the normal assumption using Model 1. }
    \label{fig:sizepower}
	\end{center}
\end{figure}

\begin{figure}[htbp]
	\begin{center}
		\includegraphics[width=7.4cm]{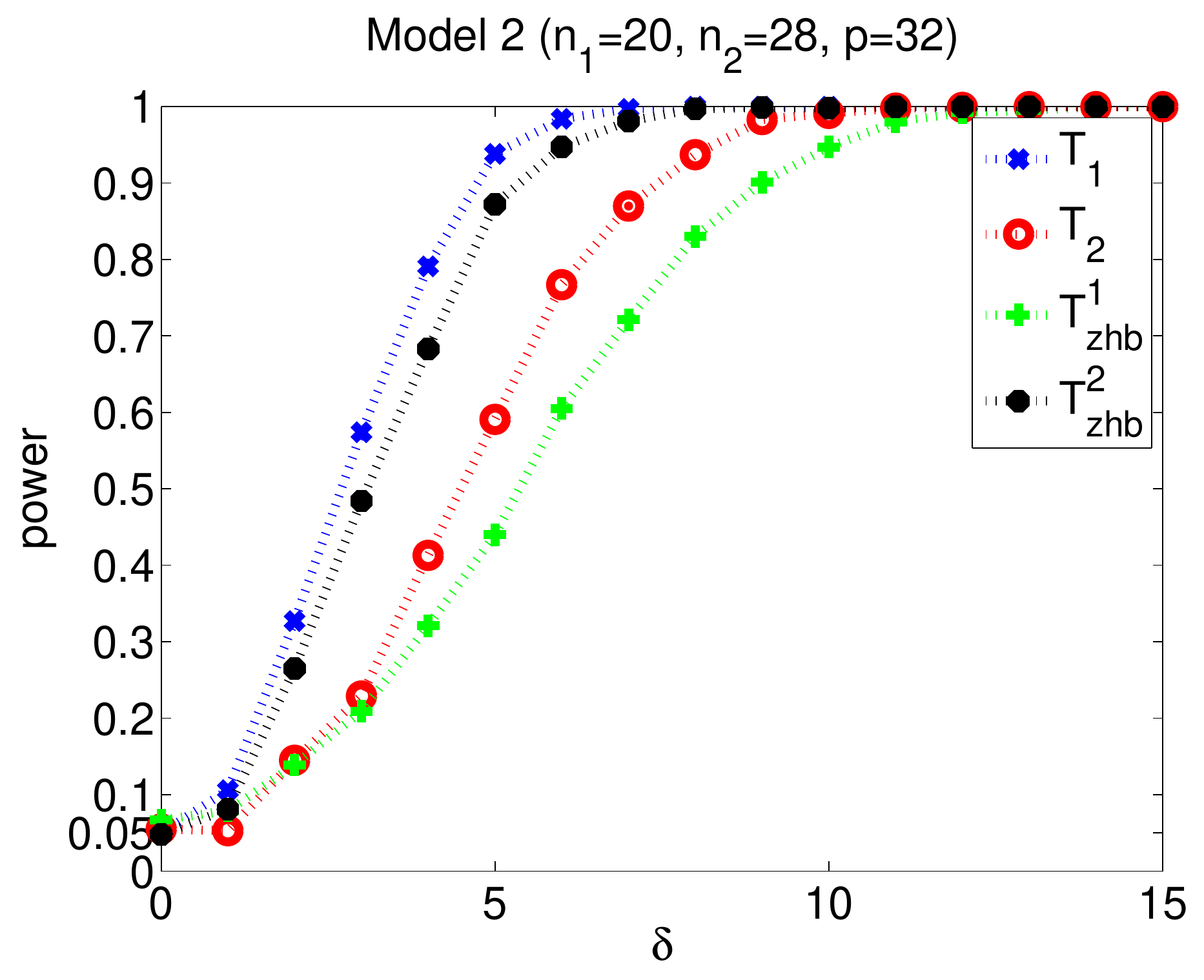}
        \includegraphics[width=7.4cm]{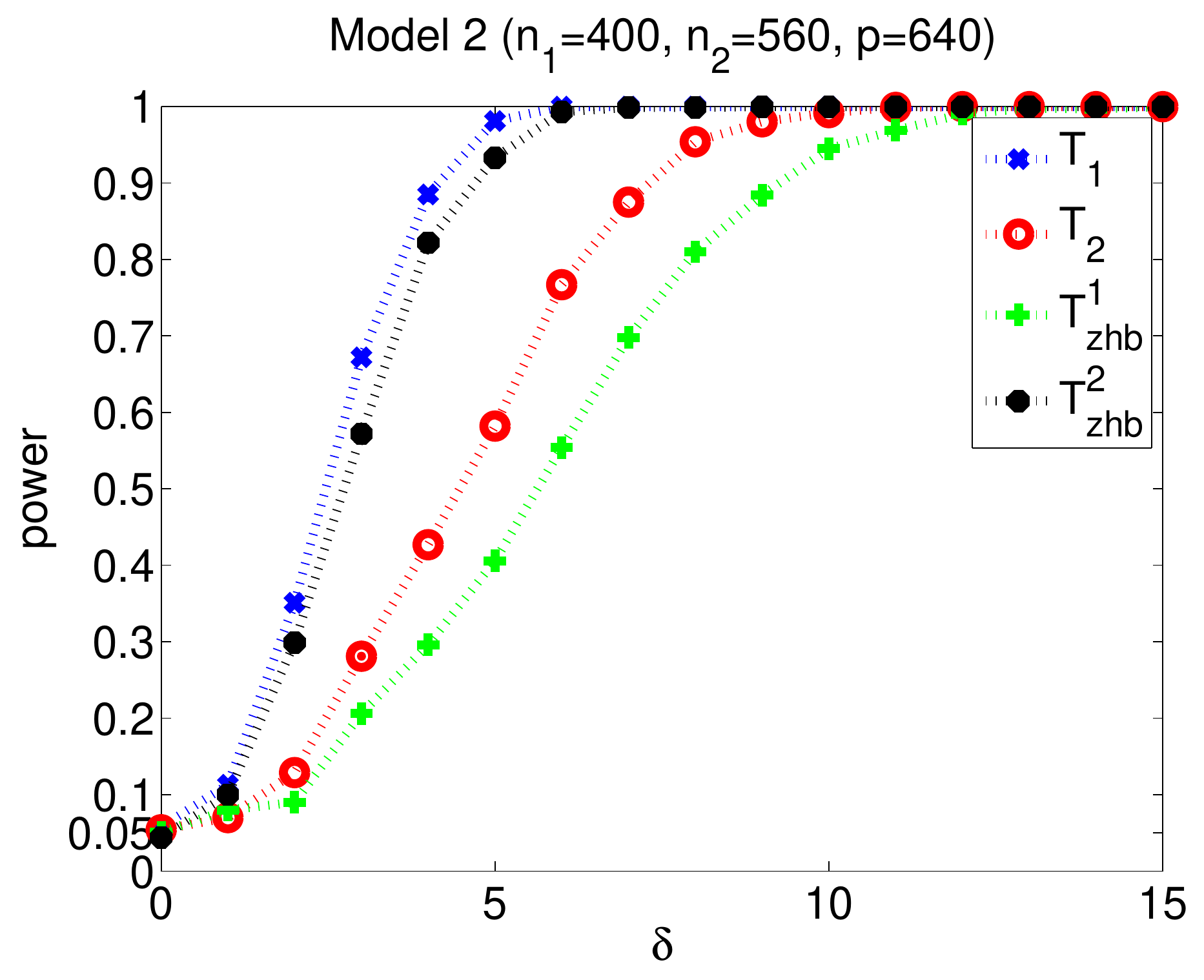}
    \caption{The empirical power for $T_1$, $T_2$, $T_{zhb}^{1}$ and $T_{zhb}^{2}$ under the normal assumption using Model 2. }
    \label{fig:sizepower2}
	\end{center}
\end{figure}

\begin{figure}[htbp]
	\begin{center}
		\includegraphics[width=7.4cm]{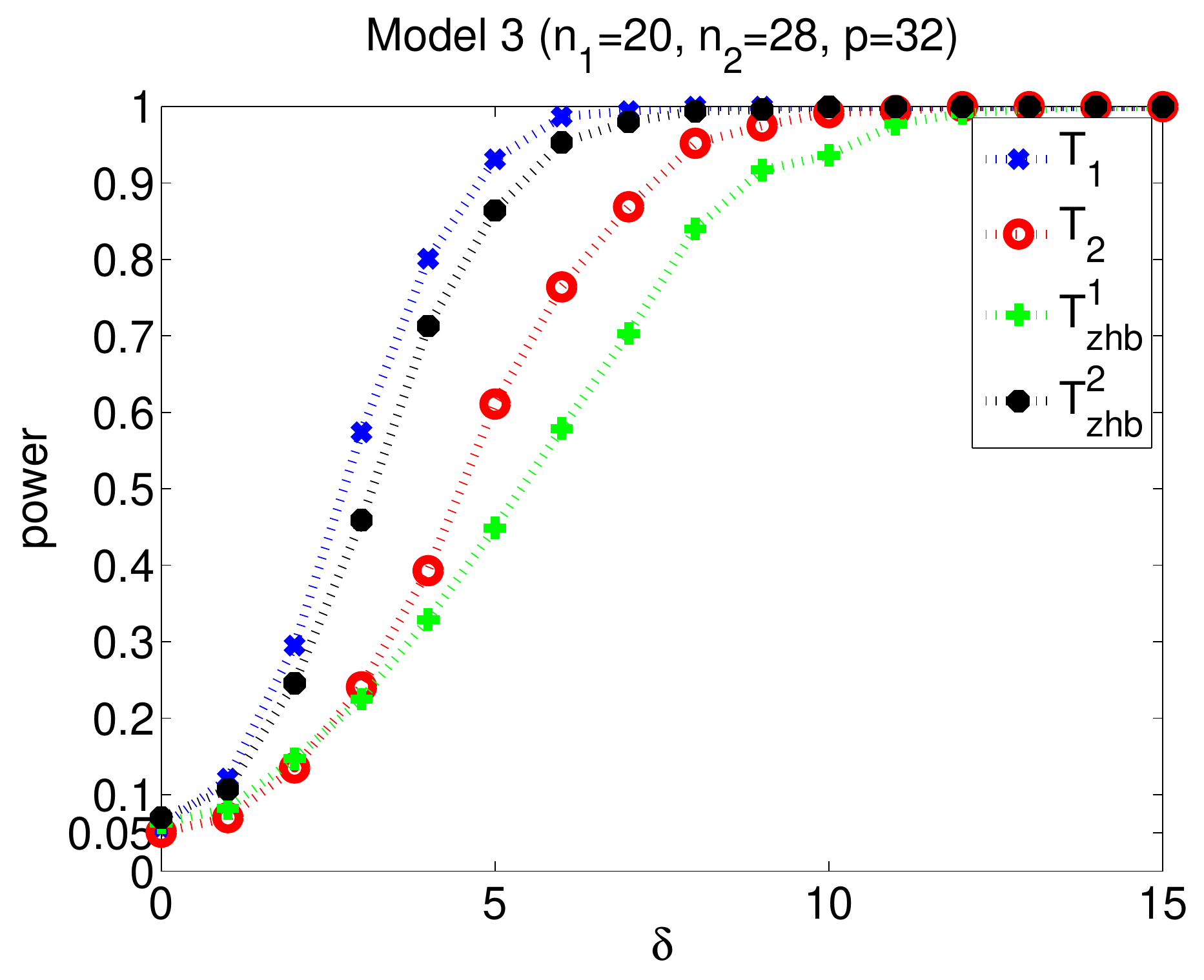}
        \includegraphics[width=7.4cm]{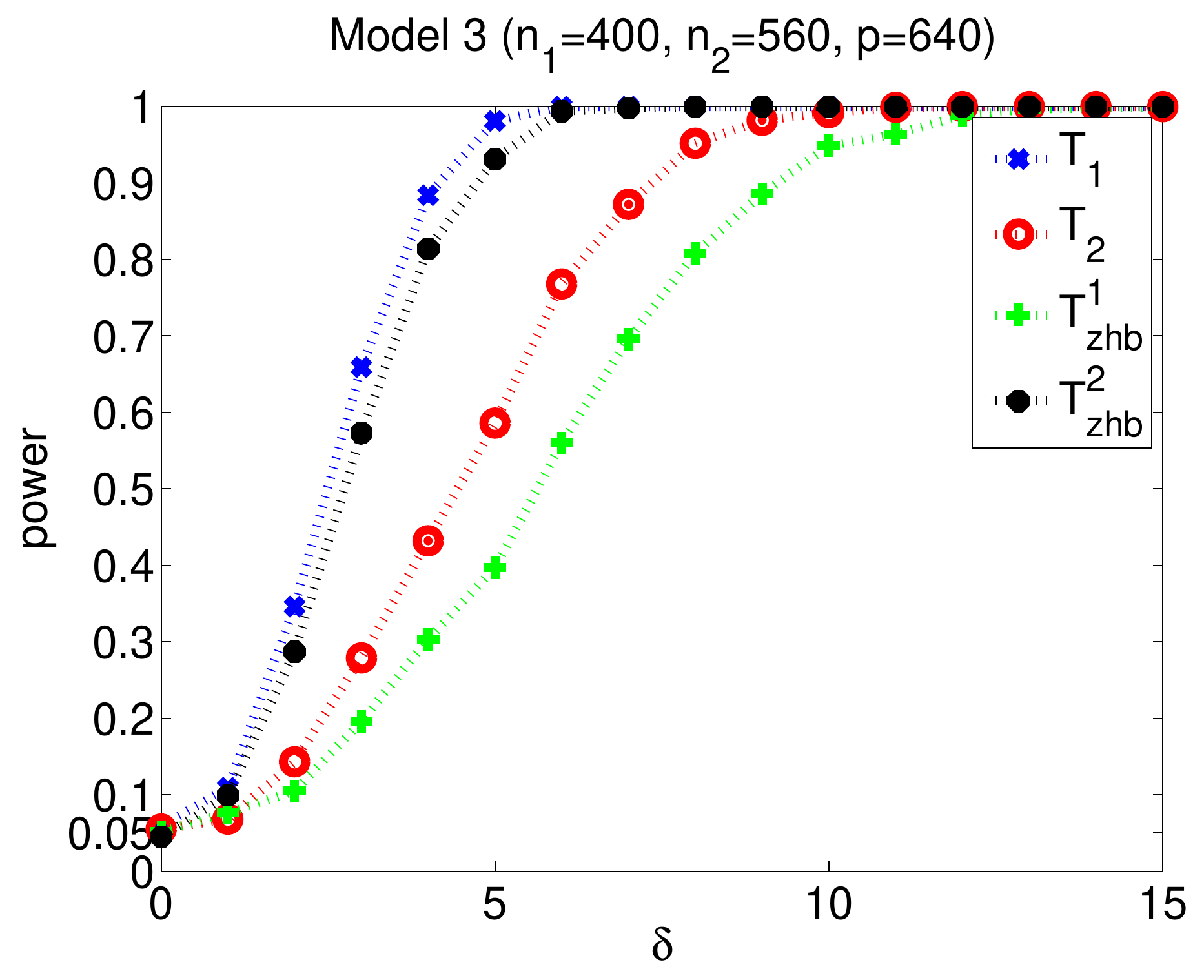}
    \caption{The empirical power for $T_1$, $T_2$, $T_{zhb}^{1}$ and $T_{zhb}^{2}$ under the normal assumption using Model 3. }
    \label{fig:sizepower3}
	\end{center}
\end{figure}

\begin{figure}[htbp]
	\begin{center}
		\includegraphics[width=7.4cm]{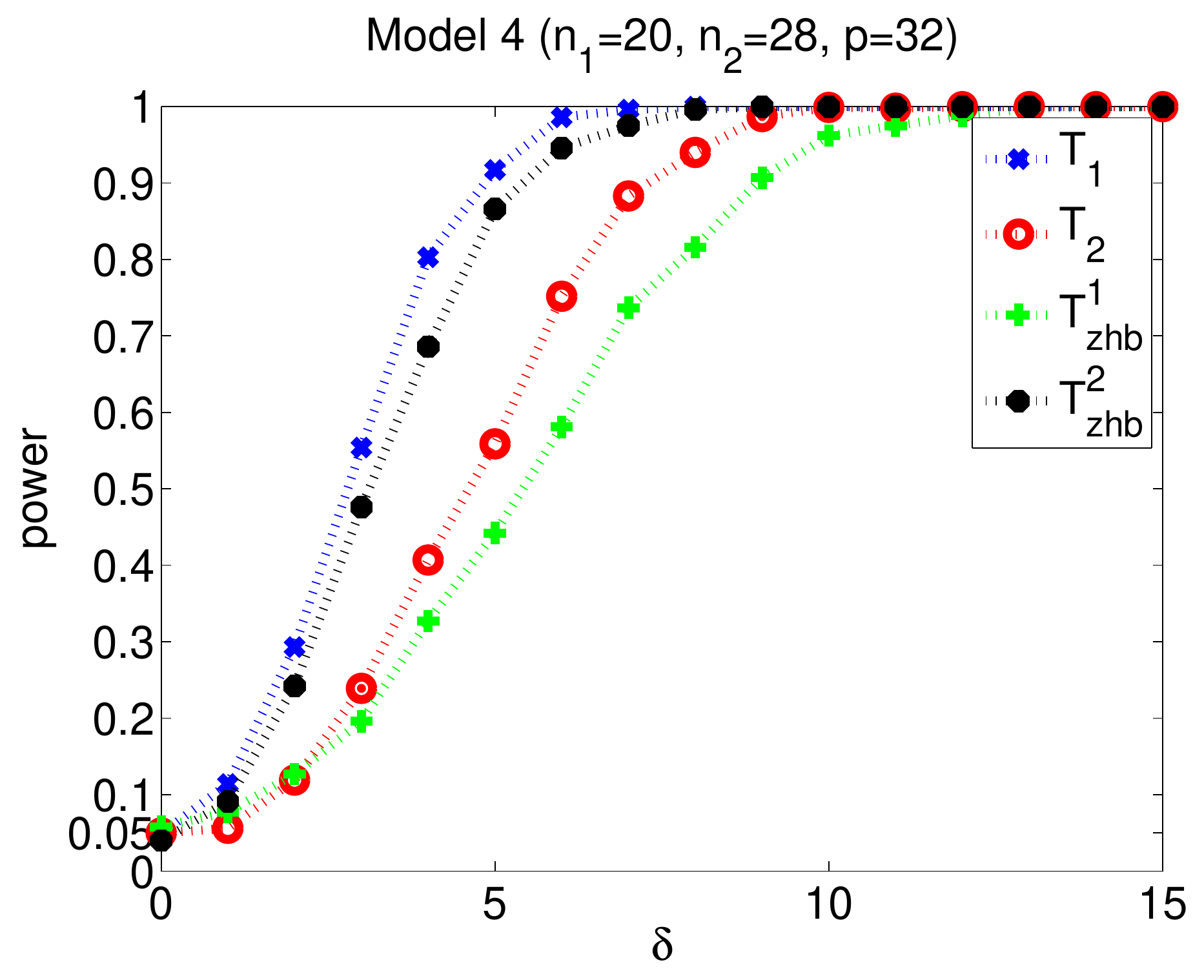}
        \includegraphics[width=7.4cm]{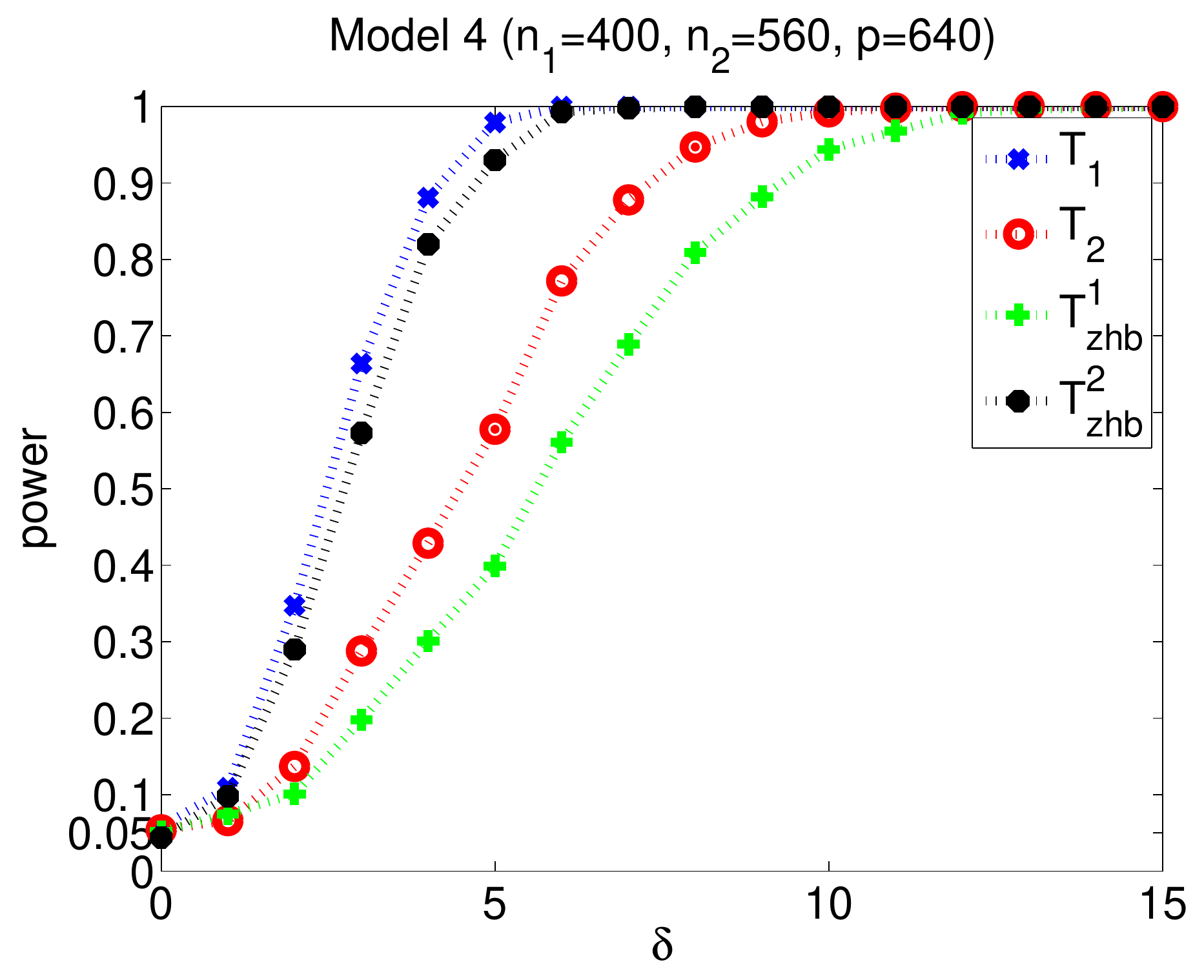}
    \caption{The empirical power for $T_1$, $T_2$, $T_{zhb}^{1}$ and $T_{zhb}^{2}$ under the normal assumption using Model 4. }
    \label{fig:sizepower4}
	\end{center}
\end{figure}

 We also compare $T_1$, $T_2$, $T_{zhb}^{1}$ and $T_{zhb}^{2}$ for four models under normal distribution condition in Figures \ref{fig:sizepower}-\ref{fig:sizepower4}. From these simulations, we find  that $T_1$ seems to be more powerful than the other three statistics under these four models. In addition, because the power of $T_{zhb}^{1}$ and $T_{zhb}^{2}$ also depend only on the eigenvalues of $\Sigma_1\Sigma_2^{-1}$, the four figures appear quite similar.

\subsection{Goodness-of-fit tests}
In this subsection, we compare the goodness-of-fit tests of the proposed statistics with those of $T_{zhb}^{1}$, $T_{zhb}^{2}$ and $T_{lc}$.
All five statistics can be used to test hypothesis \eqref{hypo}, and their asymptotic distributions are all standard normal distributions.
We generate normal $x_{j}^{(l)}$ based on Model 1 and repeat the simulation 1,000 times to obtain 1,000 random points under the null hypothesis  for all five statistics.
The J-B test statistic represents a goodness-of-fit test to determine whether the skewness and kurtosis of sample data fit a standard normal distribution, and it is used here to determine whether the 1,000 random points follow a standard normal distribution.
The J-B test statistic is expressed as follows:
$$JB=\frac{n+m-1}{6}(S^2+\frac{1}{4}(C-3)^2),$$
where $n$ is the number of observations, $S$ is the sample skewness, $C$ is the sample kurtosis, and $m$ is the number of regressors. The K-S test can also be used to compare a sample with a reference probability distribution and is defined as
$$D_n=\sup_{x}|F_n(x)-F(x)|,$$
where $F_n(x)$ is the empirical distribution function and $F(x)$ is a given cumulative distribution.
When the value of $D_n$ is small, the sample is likely to obey the given distribution $F(x)$.
Table \ref{tabJBt} provides the p-values of the J-B and K-S tests for the five compared statistics.
$(n_1,n_2,p)=(20,28,32)$, $(n_1,n_2,p)=(20,32,28)$, $(n_1,n_2,p)=(32,20,28)$, and $(n_1,n_2,p)=(32,28,20)$ represent the following cases: $y_1>1,y_2>1$; $y_1>1,y_2<1$; $y<1,y_2>1$; and $y_1<1,y_2<1$, respectively.
From Table \ref{tabJBt}, we find that the p-values exceed 0.05, in most cases, which means   we do not have sufficient evidence to reject the null hypothesis, especially when the dimension is high. Thus, we can conclude that the proposed statistics fit well for finite samples.  

\begin{table}[htbp]
\centering
\footnotesize
\begin{tabular}{cccccccccccc}
\\ \hline
 \multirow{3}{4em}{Method}&\multicolumn{8}{c}{Low-dimensional $(n_1,n_2,p)$}\\ \cline{2-12}

  &\multicolumn{2}{c}{$(20,28,32)$}&&\multicolumn{2}{c}{$(20,32,28)$}&&\multicolumn{2}{c}{$(32,20,28)$}
  &&\multicolumn{2}{c}{$(32,28,20)$}\\\cline{2-12}

  &{J-B test} &K-S test&&J-B test&K-S test&&J-B test&K-S test&&J-B test&K-S test \\ \cline{2-3}\cline{5-6}\cline{8-9}\cline{11-12}

$T_1$        &0.986&0.291&& 0.570&0.243&& 0.711&0.582&& 0.308&0.438 \\
$T_{zhb}^{1}$&0    &0.473&& 0.010&0.004&& 0.016&0.030&& 0.303&0.131 \\
$T_2$        &0.334&0.748&& 0.573&0.198&& 0.217&0.102&& 0.936&0.053 \\
$T_{zhb}^{2}$&0.059&0.002&& 0.492&0.243&& 0.115&0.021&& 0.075&0.043 \\
{$T_{lc}$}   &0.097&0.049&& 0.534&0.019&& 0.630&0.103&& 0.942&0.141 \\\hline

 \multirow{3}{4em}{Method}&\multicolumn{8}{c}{High-dimensional $(n_1,n_2,p)$}\\ \cline{2-12}

  &\multicolumn{2}{c}{$(400,560,640)$}&&\multicolumn{2}{c}{$(400,640,560)$}&&\multicolumn{2}{c}{$(640,400,560)$}
  &&\multicolumn{2}{c}{$(640,560,400)$}\\\cline{2-12}

  &J-B test&K-S test&&J-B test&K-S test&&J-B test&K-S test&&J-B test&K-S test\\ \cline{2-3}\cline{5-6}\cline{8-9}\cline{11-12}

$T_1$        &0.406&0.167&& 0.933&0.395&& 0.952&0.686&& 0.958&0.613 \\
$T_{zhb}^{1}$&0.197&0.556&& 0.628&0.505&& 0.084&0.022&& 0.468&0.951 \\
$T_2$        &0.094&0.326&& 0.366&0.534&& 0.568&0.159&& 0.581&0.978 \\
$T_{zhb}^{2}$&0.276&0.328&& 0.996&0.647&& 0.432&0.464&& 0.616&0.432 \\
{$T_{lc}$} &0.892&0.181&& 0.369&0.252&& 0.400&0.863&& 0.192&0.772 \\\hline
		
\end{tabular}
\caption{P-values of the J-B and K-S tests for $T_1$, $T_2$, $T_{zhb}^{1}$, $T_{zhb}^{2}$ and $T_{lc}$ under the normal assumption.}\label{tabJBt}
\end{table}

\section{Real data analysis}
In this section, we apply our method to an analysis of Standard and Poor's (S$\&$P) 500 index. The S$\&$P 500 index comprises 505 stocks trading on the American stock exchanges issued by 500 companies. There are ten Global Industry Classification Standard (GICS) sectors. After removing outlier data, each GICS sector includes many companies that reported 754 daily closing prices for three years.

On the stock market, volatility is the most frequently considered index, and it is critical for risk assessment.
Volatility is the degree of variation in a trading price series over time as measured by the standard deviation of logarithmic returns.
Many types of volatility exist, including historical volatility, which is a time series consisting of past market prices, current volatility, and future volatility.
In some sense, historical volatility provides a forecast of what a stock return will tend to be over a future period of time.
We focus on the historical stock return volatility in 2 sectors: Energy and Information Technology.
In each sector, our goal is to test whether the two covariance {matrices of the stock} returns are the same over time.
The diagonal elements of the covariance matrix provide a rough representation of volatility, while the off-diagonal elements indicate the reciprocal effect between two stocks.

We consider the seasonal volatility of the daily returns of the selected stocks and choose the data from the first season of 2012 as sample 1 and the data from the second season in 2012 (S2Y12) to the last season in 2014 (S4Y14) as sample 2.
We report the p-values of six tests based on the Energy sector in Table \ref{tabEn} and the Information Technology sector in Table \ref{tabRIT}.
From these results, we find that most of the p-values are smaller than $0.05$. Thus, there is strong evidence that most of the  covariance matrices are different, however,  caution is advised when utilizing the assumption that the  returns are identically distributed.
In addition, there is no evidence to reject the null hypothesis for the first seasons of 2013 and 2014 in Table \ref{tabEn}.
Thus, we suspect the presence of a periodic property in the covariance matrices.
\begin{table}[h]
\centering
\footnotesize
\begin{tabular}{cccccccccccc}
   \\ \cline{1-12}
{Season}         &{S2Y12}&{S3Y12}&{S4Y12}&{S1Y13}&{S2Y13}&{S3Y13}&{S4Y13}&{S1Y14}&{S2Y14}&{S3Y14}&{S4Y14}   \\ \hline
$T_1$            &0     & 0     & 0.005 & 0.237 & 0.116 & 0.173 & 0.290 & 0.491 & 0.177 & 0.033 & 0.165 \\ \hline
$T_2$            &0     & 0     & 0     & 0.256 & 0.378 & 0.153 & 0     & 0.030 & 0.317 & 0.016 & 0.002 \\ \hline
$T_{zhb}^{1}$    &0     & 0     & 0.001 & 0.095 & 0.078 & 0.036 & 0     & 0.001 & 0.029 & 0     & 0     \\ \hline
$T_{zhb}^{2}$    &0     & 0.001 & 0.007 & 0.340 & 0.051 & 0.119 & 0.495 & 0.400 & 0.194 & 0.044 & 0.057 \\ \hline
$T_{lc}$         &0     & 0     & 0.162 & 0.147 & 0     & 0     & 0     & 0.007 & 0     & 0     & 0     \\ \hline
$T_{clx}$        &0.004 & 0.035 & 0.129 & 0.404 & 0.002 & 0.065 & 0.003 & 0.110 & 0.027 & 0.145 & 0.003 \\ \hline
\end{tabular}
\caption{P-values from the tests of the two covariance matrices of daily stock returns in the Energy sector.}\label{tabEn}
\end{table}

\clearpage

\begin{table}[!h]
\footnotesize
\begin{tabular}{cccccccccccc}
\\ \cline{1-12}
{Season}         &{S2Y12}&{S3Y12}&{S4Y12}&{S1Y13}&{S2Y13}&{S3Y13}&{S4Y13}&{S1Y14}&{S2Y14}&{S3Y14}&{S4Y14}   \\ \hline
$T_1$            &0.1672&0.0110&0.0024&0.0025&0     &0.0011&0.3997&0     &0     &0.0015&0.4679 \\ \hline
$T_2$            &0     &0     &0     &0     &0     &0     &0     &0     &0     &0     &0      \\ \hline
$T_{zhb}^{1}$    &0     &0     &0     &0     &0     &0     &0     &0     &0     &0     &0      \\ \hline
$T_{zhb}^{2}$    &0.0102&0.0001&0.2511&0.    &0     &0     &0.0503&0     &0     &0     &0.0789 \\ \hline
$T_{lc}$         &0     &0.2761&0.0203&0     &0.1036&0.0066&0.1944&0     &0.0006&0     &0.0359 \\ \hline
$T_{clx}$        &0.0004&0.0570&0.6058&0.8945&0.1912&0.2810&0.9483&0.5678&0.1954&0.0707&0.2231 \\ \hline
\end{tabular}
\caption{P-values from the tests of the two covariance matrices of daily stock returns in the Information Technology sector.}\label{tabRIT}
\end{table}

\section{Appendix}
In this section,  we prove Theorem \ref{th1} and Theorem \ref{th2} by using the same tools. To simplify the notation, we use $c_1$, $c_2$, $y_1$, $y_2$, $h$, and $\alpha$ instead of $c_{n_1}$, $c_{n_2}$, $y_{n_1}$, $y_{n_2}$, $h_{n}$, and $\alpha_{n}$, respectively. We calculate $l_n$, $\mu_n$ and $\nu_n$ in Theorem \ref{th1} and Theorem \ref{th2}  based on Theorem 1.1 and Theorem 1.6 from Bai et al. \cite{Bai2015}, included here for convenience.

\begin{lem}[Theorem 1.1 in Bai et al. \cite{Bai2015}]\label{lem1}
Under assumptions {(1)--(4)}, the empirical spectral distribution (ESD) of the random Beta matrix $\bbB_{n}^{x}(\bbX_{1},\bbX{2})$ with probability $1$ weakly converges to a non-random distribution $\bbF^{\bbB}(x)$ whose density function is distributed as
 \begin{displaymath}
 \left\{ \begin{array}{ll}
 \frac{(\alpha+1)\sqrt{(x_r-x)(x-x_l)}}{2\pi y_1 x(1-x)} & \textrm{when $x_l<x<x_r$}\\
0 & \textrm{otherwise},\\
\end{array} \right.
\end{displaymath}
where $x_l,\:x_r=\frac{y_2(h\mp y_1)^2}{(y_1+y_2)^2}.$
\end{lem}

\begin{lem}[Theorem 1.6 in Bai et al. \cite{Bai2015}]\label{lem2}
In addition to conditions {(1)-(4)}, we further assume that $f_1,......f_k$ are analytic functions on an open region containing the interval $[c_{le},c_{ri}]$, where $c_{le}=v^{-1}(1-\sqrt{y_1})^2$, $c_{ri}=1-\alpha v^{-1}(1-\sqrt{y_2})^2$, and $v$ is defined as $v=(1+\frac{y_1}{y_2})(1-\sqrt{\frac{y_1y_2}{y_1+y_2}})^2$.

Then, as min $(n_1,n_2,p)\rightarrow\infty$, the random vector $$(\int f_{\gamma} \mathrm{d}\mathbf{G_n}(x)),\quad {\gamma=1,......,k},$$ where $\mathbf{G_n}(x)={p(\bbF^{\bbB_{n}^{x}(\bbX_{1},\bbX_{2})}(x)}-\bbF^{\bbB}(x))$ converges weakly to a Gaussian vector $(\mathbf{G}_{f_1},...,\mathbf{G}_{f_k})$ with the mean function
\begin{equation}\label{mean1}
\textbf{E}\mathbf{G}_{f_{\gamma}}=
\frac {1}{4\pi i}\oint f_{\gamma}(\frac{z}{\al+z})\mathrm{d}\log(\frac{(1-y_2)m_3^2(z)+2m_3(z)+1-y_1}{(1+m_3(z))^2})
\end{equation}
\begin{equation}\label{mean2}
+\frac {\Delta_1}{2\pi i}\oint y_1f_{\gamma}(\frac {z}{\al+z})(1+m_3)^{-3}\mathrm{d} m_3(z)
\end{equation}
\begin{equation}\label{mean3}
+\frac {\Delta_2}{4\pi i}\oint f_{\gamma}(\frac {z}{\al+z})(1-y_2m_3^2(z)(1+m_3(z))^{-2})\mathrm{d}\log(1-y_2m_3^2(z)(1+m_3(z))^{-2}),
\end{equation}
and the covariance function
\begin{equation}\label{covariance1}
\textbf{Cov}(\mathbf{G}_{f_{\gamma}},\mathbf{G}_{f_{\gamma'})=-\frac{1}{2\pi^2}}\oint
\oint\frac{f_{\gamma}(\frac{z_1}{\al+z_1})f_{\gamma'}(\frac{z_2}{\al+z_2})
\mathrm{d}m_3(z_1)\mathrm{d}m_3(z_2)}{(m_3(z_1)-m_3(z_2))^2}
\end{equation}
\begin{equation}\label{covariance2}
-\frac{y_1\Delta_1+y_2\Delta_2}{4\pi^2}\oint\oint\frac{f_{\gamma}(\frac{z_1}{\al+z_1})f_{\gamma'}
(\frac{z_2}{\al+z_2})\mathrm{d}m_3(z_1)\mathrm{d}m_3(z_2)}{(m_3
(z_1)+1)^2(m_3(z_2)+1)^2},
\end{equation}
where
\bqn
m_0(z)&=&\frac{(1+y_1)(1-z)-\alpha z(1-y_2)+\sqrt{((1-y_1)(1-z)+\alpha z(1-y_2))^2-4\alpha z(1-z)}}{2z(1-z)(y_1(1-z)+\alpha zy_2)}-\frac 1{z},\\
m_1(z)&=&\frac{\alpha}{(\alpha+z)^2}m_0(\frac{z}{\alpha+z})-\frac1{\alpha+z} , \quad m_2(z)=-z^{-1}(1-y_1)+y_1m_1(z),\\
m_{mp}^{y_2}(z)&=&\frac{1-y_2-z+\sqrt{(z-1-y_2)^2-4y_2}}{2y_2z} ,\quad  m_3(z)=y_2m_{mp}^{y_2}(-m_2(z))+(m_2(z))^{-1}(1-y_2).
\eqn
The above contour integrals can be evaluated on any contour enclosing the interval $[\frac{\alpha c_{le}}{1-c_{le}},\frac{\alpha c_{ri}}{1-c_{ri}}]$: here, $i$ represents an imaginary unit.
\end{lem}

Now, we are in position to prove Theorem \ref{th1}.

\begin{proof} {of the limit part $l_n$ in Theorem \ref{th1}.} To calculate the limit part
\begin{align}\label{limitL1}
p\int x\frac{(1+\alpha)\sqrt{(x_r-x)(x-x_l)}}{2\pi y_1x(1-x)}\mathrm{d}x,
\end{align}
we first perform the transformations $x=\frac{y_2|y_1+h\xi|^2}{(y_1+y_2)^2}$ and  $1-x=\frac{y_1|y_2-h\xi|^2}{(y_1+y_2)^2}$. Because $x$ appears in the molecular orbital
of the integrand function, no residue is related to $y_1$. We assume that $y_1>1$. 
Clearly, as $x$ moves from $\frac{y_2(h- y_1)^2}{(y_1+y_2)^2}$ to $\frac{y_2(h+y_1)^2}{(y_1+y_2)^2}$ two times, $\xi$ runs along the unit circle in the positive direction.
Thus, integral \eqref{limitL1} is equivalent to
\begin{equation*}
p\frac{h^2i}{4\pi(y_1+y_2)}\oint\frac{(\xi^2-1)^2}{\xi^2(y_2-h\xi)(\xi-\frac h{y_2})}\mathrm{d}\xi.
\end{equation*} 
According to the residue theorem, we obtain two poles $\{0, \frac{h}{y_2}\} $ in the unit disc when $y_2>1$, and the residues are $$-\frac{y_2^2+h^2}{y_2h^2},\quad \frac{y_2^2-h^2}{y_2h^2}.$$
Then, under the assumption $y_2>1$, \eqref{limitL1} yields 
\bqn
p\frac{h^2i}{4\pi(y_1+y_2)}\cdot2\pi i\cdot(-\frac{y_2^2+h^2}{y_2h^2}+\frac{y_2^2-h^2}{y_2h^2})
=p\frac{h^2}{(y_1+y_2)y_2}.
\eqn
In the same way, we obtain two poles
$\{0, \frac{y_2}{h}\}$  and two  residues
$$-\frac{y_2^2+h^2}{y_2^2h},\quad -\frac{y_2^2-h^2}{y_2^2h}$$
under the assumption $y_2<1.$
Then, we can calculate that \eqref{limitL1} is
\bqn
-p\frac{hy_2i}{4\pi(y_1+y_2)}\cdot2\pi i\cdot(\frac{y_2^2+h^2}{y_2^2h}+\frac{y_2^2-h^2}{y_2^2h})
=p\frac{y_2}{y_1+y_2},
\eqn
which completes the proof.
\end{proof}

\begin{proof} {of the mean part $\mu_n$ in Theorem \ref{th1}.}
Because $m_3$ satisfies the equation
$$z=-\frac{m_3(z)(m_3(z)+1-y_1)}{(1-y_2)(m_3(z)+\frac{1}{1-y_2})},$$ we perform an integral conversion $z=(1+hr\xi)(1+\frac{h}{r\xi})/(1-y_2)^2$, where $r$ is a number greater than but close to $1$.
For the same reason, we assume $y_1>1$ without loss of generality.
The pole related to $y_2$ of the integrand is $\frac{h}{y_2}$ when $y_2>1$. The integral value is not changed by the transformation $\xi=\frac{1}{\xi}$; however, the residue point in the unit disc becomes $\frac{y_2}{h}$, which is the residue point under the assumption $y_2<1$.
Therefore, we can assume that $y_2>1$. By solving the equation $$\frac{(1+hr\xi)(1+\frac{h}{r\xi})}{(1-y_2)^2}=-\frac{m_3(m_3+1-y_1)}{(1-y_2)m_3+1}$$ we obtain $m_3=-(1+hr\xi)/(1-y_2)$ or $m_3=-(1+\frac{h}{r\xi})/(1-y_2)$. When $z$ runs in the positive direction along the unit circle around the support of $\bbF^{\bbB}(x)$, $m_3$ runs in the opposite direction. Therefore, when $y_2>1$, we choose the outcome $m_3=-(1+\frac{h}{r\xi})/(1-y_2)$. Based on the above discussion, we have $$\frac{z}{\alpha+z}=\frac{y_2|1+hr\xi|^2}{|y_2+hr\xi|^2}.$$
Therefore, we obtain the mean part
\bqn
\eqref{mean1}=\lim_{r\downarrow1}\frac {1}{4\pi i}\oint_{|\xi|=1}\frac{y_2(1+hr\xi)(1+\frac{h}{r\xi})}{(y_2+hr\xi)(y_2+\frac{h}{r\xi})}\cdot
(\frac1{\xi-\frac1r}+\frac1{\xi+\frac1r}-\frac2{\xi+\frac h{y_2r}})\mathrm{d}\xi,
\eqn
\bqn
\eqref{mean2}=-\lim_{r\downarrow1}\frac {\Delta_1}{2\pi i}\oint y_1\frac{y_2(1+hr\xi)(1+\frac{h}{r\xi})}{(y_2+hr\xi)(y_2+\frac{h}{r\xi})}
\cdot\frac{(1-y_2)^2h}{y_2^3}\cdot\frac{\xi}{(\xi+\frac{h}{y_2r})^3}\mathrm{d}\xi,
\eqn
\bqn
\eqref{mean3}=\lim_{r\downarrow1}\frac {\Delta_2}{4\pi i}\oint\frac{y_2(1+hr\xi)(1+\frac{h}{r\xi})}{(y_2+hr\xi)(y_2+\frac{h}{r\xi})}
\cdot\frac{(y_2-1)(\xi^2-\frac{h^2}{y_2r^2})}{
y_2(\xi+\frac{h}{y_2r})^2}\cdot( \frac{2\xi}{\xi^2-\frac{h^2}{y_2r^2}}-
\frac{2}{\xi+\frac{h}{y_2r}})\mathrm{d}\xi.
\eqn
Based on the residue theorem, \eqref{mean1} has three poles,
$$\frac{1}{r},\quad -\frac{1}{r},\quad -\frac{h}{y_2r},$$
and three residues,
$$\frac{y_2(1+h)^2}{(y_2+h)^2},\quad \frac{y_2(1-h)^2}{(y_2-h)^2},\quad \frac{2h^2y_1}{(y_1+y_2)^2}.$$
Then,
\bqn
\eqref{mean1}=0.
\eqn
Similarly,
\bqn
\eqref{mean2}=-\Delta_1\frac{y_1^2y_2^2h^2}{(y_1+y_2)^4}
~~\mbox{and}~~
\eqref{mean3}=-\Delta_2\frac{y_1^2y_2^2h^2}{(y_1+y_2)^4}.
\eqn
Finally, we obtain the result of the mean part
\bqn
\mu_{n}=-\frac{\Delta_1y_{1}^2y_{2}^2h^2}{(y_{1}+y_{2})^4}-
\frac{\Delta_2y_{1}^2y_{2}^2h^2}{(y_{1}+y_{2})^4},\eqn
which completes the proof.
\end{proof}

\begin{proof} {of the variance part in Theorem \ref{th1}.}
To calculate the variance {part \eqref{covariance1} and \eqref{covariance2}},
we make analogous integral conversions $$z_1=(1+hr_1\xi_1)(1+\frac{h}{r_1\xi_1})/(1-y_2)^2$$ and  $$z_2=(1+hr_2\xi_2)(1+\frac{h}{r_2\xi_2})/(1-y_2)^2.$$ Therefore, the relationship between $\xi_l$ and $m_3(z_l)$, where $l=1$, $2$, is as follows
$$m_3(z_1)=-\frac{1+\frac{h}{r_1\xi_1}}{(1-y_2)}~~\mbox{and}~~
m_3(z_2)=-\frac{1+\frac{h}{r_2\xi_2}}{(1-y_2)}.$$
We assume that $r_1<r_2$ {without loss of generality}. When $y_1>1$ and $y_2>1$, according to the residue theorem, we obtain
\bqn
\eqref{covariance1}=2\lim_{r_2\downarrow1}\oint{\frac{1}{2\pi i}
\cdot[\frac{y_2(1+hr_2\xi_2)(1+\frac{h}{r_2\xi_2})}{(y_2+hr_2\xi_2)
(y_2+\frac{h}{r_2\xi_2})}]}
\eqn
\bqn
{{\cdot\{\lim_{r_1\downarrow 1}}\oint\frac{1}{2\pi i}\cdot
[\frac{y_2(1+hr_1\xi_1)(1+\frac{h}{r_1\xi_1})}{(y_2+hr_1\xi_1)
(y_2+\frac{h}{r_1\xi_1})}]\cdot
\frac{r_1}{(r_1\xi_1-r_2\xi_2)^2}\mathrm{d}\xi_1\}\cdot r_2}\mathrm{d}\xi_2.
\eqn
Only one pole $-\frac{h}{y_2}$ exists in the unit disc for $\xi_1$, and the respective residue point is $$\frac{y_1y_2h(y_2-1)}{(y_1+y_2)(h+y_2\xi_2)^2}.$$ Then, we obtain the formula
\bqn
{\eqref{covariance1}=2\lim_{r_2\downarrow1}\oint\frac{1}{2\pi i}\cdot
[\frac{y_2(1+hr_2\xi_2)(1+\frac{h}{r_2\xi_2})}{(y_2+hr_2\xi_2)
(y_2+\frac{h}{r_2\xi_2})}]\cdot
\frac{y_1(y_2-1)h}{r_2^2y_2(y_1+y_2)}\cdot}\frac{1}{(\xi_2+\frac{h}{y_2r
_2})^2}
r_2\mathrm{d}\xi_2,
\eqn which has only one pole, $-\frac{y_2}{r_2\xi_2}$, in the unit disc. Then, we have
\bqn
{\eqref{covariance1}=}\frac{2y_1^2y_2^2h^2}{(y_1+y_2)^4}.
\eqn
In the same way, we have the following calculation
\bqn
{\eqref{covariance2}
=(y_1\Delta_1+y_2\Delta_2)\cdot\{\lim_{r_1\downarrow1}\oint\frac{1}{2\pi i}\cdot
[\frac{y_2(1+hr_1\xi_1)(1+\frac{h}{r_1\xi_1})}{(y_2+hr_1\xi_1)
(y_2+\frac{h}{r_1\xi_1})}]\cdot}
\frac{(1-y_2)r_1h\mathrm{d}\xi_1}{(y_2r_1\xi_1+h)^2}\}
\eqn
\bqn
{\cdot\{\lim_{r_2\downarrow1}\oint\frac{1}{2\pi i}\cdot
[\frac{y_2(1+hr_2\xi_2)(1+\frac{h}{r_2\xi_2})}{(y_2+hr_2\xi_2)
(y_2+\frac{h}{r_2\xi_2})}]\cdot}
\frac{(1-y_2)r_2h\mathrm{d}\xi_2}{(y_2r_2\xi_2+h)^2}\}
\eqn
\bqn
{=(y_1\Delta_1+y_2\Delta_2)\cdot}
\frac{y_1^2y_2^2h^4}{(y_1+y_2)^6}.
\eqn
Therefore, we can conclude that 
$$\nu_{n}^2=\frac{2y_{1}^2y_{2}^2h^2}{(y_{1}+y_{2})^4}+
\frac{(y_{1}\Delta_1+y_{2}\Delta_2)y_{1}^2y_{2}^2h^4}{(y_{1}+y_{2})^6},$$ which completes the proof.
\end{proof}

We now give the proof for the statistic $T_2$.

\begin{proof} {of the limit part $\widetilde{l}_n$ in Theorem \ref{th2}.}
To calculate the limit part $p\widetilde{l}_n $,
 \begin{align}\label{limitL2}
p\int[c_1(\frac{x}{c_1}-1)^2+c_2(\frac{1-x}{c_2}-1){^2]\cdot}\frac{(1+\alpha)}{\sqrt{(x_r-x)(x-x_l)}}{2\pi y_1x(1-x)}\mathrm{d}x,
\end{align}
where
\bqn
x_l,\:x_r=\frac{y_2(h\mp y_1)^2}{(y_1+y_2)^2}.
\eqn
Performing the transformation
\bqn
x=\frac{y_2|y_1+h\xi|^2}{(y_1+y_2)^2},\
\eqn
we have
\bqn
1-x=\frac{y_1|y_2-h\xi|^2}{(y_1+y_2)^2}.
\eqn

Thus,
\begin{equation}\label{limitL1c2}
p\int c_1(\frac{x}{c_1}-1)^2 \cdot\frac{(1+\alpha)\sqrt{(x_r-x)(x-x_l)}}{2\pi y_1x(1-x)}\mathrm{d}x
\end{equation}
\bqn
=\frac{ph^2i}{4\pi}\oint y_2[\frac{|y_1+h\xi|^2}{(y_1+y_2)}-1]^2
\cdot\frac{(\xi^2-1)^2}{\xi^3|y_1+h\xi|^2|y_2-h\xi|^2}\mathrm{d}\xi
\eqn
\begin{equation}\label{limitL1c21}
=\frac{ph^2i}{4\pi}\oint y_2\frac{|y_1+h\xi|^2}{(y_1+y_2)^2}
\cdot\frac{(\xi^2-1)^2}{\xi^3|y_2-h\xi|^2}\mathrm{d}\xi
\end{equation}
\begin{equation}\label{limitL1c22}
+\frac{ph^2i}{4\pi}\oint y_2\frac{-2}{(y_1+y_2)}
\cdot\frac{(\xi^2-1)^2}{\xi^3|y_2-h\xi|^2}\mathrm{d}\xi
\end{equation}
\begin{equation}\label{limitL1c23}
+\frac{ph^2i}{4\pi}\oint y_2
\frac{(\xi^2-1)^2}{\xi^3|y_1+h\xi|^2|y_2-h\xi|^2}\mathrm{d}\xi.
\end{equation}
Under the assumptions $y_1>1$ and $y_2>1$, two poles exist inside the unit circle in \eqref{limitL1c21} $$0, \quad \frac{h}{y_2},$$ and the respective residues are
$$\frac{y_2}{(y_1+y_2)^2}\cdot(\frac{2y_1}{y_2}-\frac{(h^2+y_2^2)(y_1+y_2)^2}{h^2y_2^3}), \quad
\frac{y_2}{(y_1+y_2)^2}\cdot\frac{(y_1+y_2)^3(y_2-1)}{y_2^3h^2}.$$
 Applying the residue theorem to \eqref{limitL1c22}, there are two poles $$0,\quad \frac{h}{y_2}$$ and two residues
 $$\frac{-2y_2}{(y_1+y_2)}\cdot\frac{-(h^2+y_2^2)}{y_2^2h^2}, \quad
 \frac{-2y_2}{(y_1+y_2)}\cdot\frac{y_2^2-h^2}{y_2^2h^2}.$$ Following the same method, we obtain three poles $$0,\quad -\frac{h}{y_2},\quad \frac{h}{y_2}$$ and three residues
 $$-y_2\cdot\frac{1}{y_1y_2h^2},\quad \frac{1}{y_1}\cdot\frac{(y_1-1)y_2}{(y_1+y_2)h^2},\quad \frac{1}{y_1}\cdot\frac{(y_2-1)y_1}{(y_1+y_2)h^2}$$ in \eqref{limitL1c23}.
 Then, we have
\bqn
\eqref{limitL1c2}=c_1ph^2[\frac{y_1}{y_2^3}+\frac{1}{y_1y_2}-\frac{1}{y_2^2}-\frac{y_1}{y_2(y_1+y_2)}].
\eqn
In the same way,
\begin{equation}\label{limitL1c3}
p\int c_2(\frac{1-x}{c_2}-1)^2\cdot\frac{(1+\alpha)\sqrt{(x_r-x)(x-x_l)}}{2\pi y_1x(1-x)}\mathrm{d}x
\end{equation}
\bqn
=\frac{ph^2i}{4\pi}\oint y_1[\frac{|y_2-h\xi|^2}{(y_1+y_2)}-1]^2
\cdot\frac{(\xi^2-1)^2}{\xi^3|y_1+h\xi|^2|y_2-h\xi|^2}\mathrm{d}\xi
\eqn
\bqn
=c_2ph^2[\frac{y_2}{y_1^3}+\frac{1}{y_1y_2}-\frac{1}{y_1^2}-\frac{y_2}{y_1(y_1+y_2)}].
\eqn
Thus,
\bqn
\eqref{limitL2}=ph^2[\frac{1}{y_1^2}+\frac{1}{y_2^2}-\frac{1}{y_1+y_2}-\frac{1}{y_1y_2}].
\eqn

When the assumptions are $y_1>1$ and $y_2<1$, we use the same transformation
\bqn
x=\frac{y_2|y_1+h\xi|}{(y_1+y_2)^2}.
\eqn

In this case, the poles of \eqref{limitL1c2} and \eqref{limitL1c3} are $$0,\quad -\frac{h}{y_1},\quad \frac{y_2}{h}.$$
Therefore,
\bqn
\eqref{limitL1c2}=\frac{ph^2i}{4\pi}\oint y_2[\frac{|y_1+h\xi|^2}{(y_1+y_2)}-1]^2
\cdot\frac{(\xi^2-1)^2}{\xi^3|y_1+h\xi|^2|y_2-h\xi|^2}\mathrm{d}\xi
\eqn
\bqn
=\frac{ph^2i}{4\pi}\oint y_2\frac{|y_1+h\xi|^2}{(y_1+y_2)^2}
\cdot\frac{(\xi^2-1)^2}{\xi^3|y_2-h\xi|^2}\mathrm{d}\xi
\eqn
\bqn
+\frac{ph^2i}{4\pi}\oint y_2\frac{-2}{(y_1+y_2)}
\cdot\frac{(\xi^2-1)^2}{\xi^3|y_2-h\xi|^2}\mathrm{d}\xi
\eqn
\bqn
+\frac{ph^2i}{4\pi}\oint y_2
\cdot\frac{(\xi^2-1)^2}{\xi^3|y_1+h\xi|^2|y_2-h\xi|^2}\mathrm{d}\xi
\eqn
\bqn
=\frac{ph^2i}{4\pi}\cdot 2\pi i \cdot \frac{1}{(y_1+y_2)^2}\cdot[\frac{(y_1+y_2)^3(1-y_2)}{y_2^2h^2}+2y_1-\frac{(h^2+y_2^2)(y_1+y_2)^2}{h^2y_2^2}]
\eqn
\bqn
+\frac{ph^2i}{4\pi}\cdot 2\pi i \cdot \frac{2y_2}{h(y_1+y_2)}\cdot[\frac{(y_1+y_2)(y_2-1)}{y_2^2h}+\frac{(h^2+y_2^2)}{y_2^2h}]
\eqn
\bqn
+\frac{ph^2i}{4\pi}\cdot 2\pi i \cdot \frac{-y_2}{y_1h}\cdot[\frac{(y_2-1)y_1}{(y_1+y_2)hy_2}+\frac{(1-y_1)}{(y_1+y_2)h}+\frac{1}{hy_2}],
\eqn
and
\bqn
\eqref{limitL1c3}=\frac{ph^2i}{4\pi}\oint y_1[\frac{|y_2-h\xi|^2}{(y_1+y_2)}-1]^2
\cdot\frac{(\xi^2-1)^2}{\xi^3|y_1+h\xi|^2|y_2-h\xi|^2}\mathrm{d}\xi
\eqn
\bqn
=\frac{ph^2i}{4\pi}\cdot 2\pi i \cdot \frac{1}{(y_1+y_2)^2}\cdot[\frac{(y_1+y_2)^3(y_1-1)}{y_1^2h^2}+2y_2-\frac{(h^2+y_1^2)(y_1+y_2)^2}{h^2y_1^2}]
\eqn
\bqn
+\frac{ph^2i}{4\pi}\cdot 2\pi i \cdot \frac{-2}{(y_1+y_2)}\cdot[\frac{(y_1+y_2)(y_1-1)}{y_1h^2}+\frac{-(h^2+y_1^2)}{y_1h^2}]
\eqn
\bqn
+\frac{ph^2i}{4\pi}\cdot 2\pi i \cdot \frac{1}{y_2}\cdot[\frac{(1-y_2)y_1}{(y_1+y_2)h^2}+\frac{(y_1-1)y_2}{(y_1+y_2)h^2}+\frac{-1}{h^2}].
\eqn
Thus, using the residue theorem, we can conclude that 
\bqn
\eqref{limitL2}=ph^2[\frac{1}{y_1^2}+\frac{1}{y_2h^2}-\frac{1}{y_1+y_2}-\frac{1}{y_1y_2}]
\eqn
under the conditions $y_1>1$ and $y_2<1$.

When $y_1<1$ and $y_2<1$, we have three poles $$0,\quad -\frac{y_1}{h},\quad \frac{y_2}{h}.$$
 Thus,
\bqn
\eqref{limitL1c2}=\frac{ph^2i}{4\pi}\cdot2\pi i\cdot\frac{1}{(y_1+y_2)^2}\cdot[\frac{(y_1+y_2)^3(1-y_2)}{h^2y_2^2}+2y_1-
\frac{(h^2+y_2^2)(y_1+y_2)^2}{h^2y_2^2}]
\eqn
\bqn
+\frac{ph^2i}{4\pi}\cdot 2\pi i \cdot \frac{2y_2}{h(y_1+y_2)}\cdot[\frac{(y_1+y_2)(y_2-1)}{y_2^2h}+\frac{(h^2+y_2^2)}{y_2^2h}]
\eqn
\bqn
+\frac{ph^2i}{4\pi}\cdot 2\pi i \cdot \frac{-y_2}{y_1h}\cdot[\frac{(y_2-1)y_1}{(y_1+y_2)hy_2}+\frac{(y_1-1)}{(y_1+y_2)h}+\frac{1}{hy_2}],
\eqn
\bqn
\eqref{limitL1c3}=\frac{ph^2i}{4\pi}\cdot 2\pi i \cdot \frac{1}{(y_1+y_2)^2}\cdot[\frac{(y_1+y_2)^3(1-y_1)}{y_1^2h^2}+2y_2-\frac{(h^2+y_1^2)(y_1+y_2)^2}{h^2y_1^2}]
\eqn
\bqn
+\frac{ph^2i}{4\pi}\cdot 2\pi i \cdot \frac{-2}{(y_1+y_2)}\cdot[\frac{(y_1+y_2)(1-y_1)}{y_1h^2}+\frac{-(h^2+y_1^2)}{y_1h^2}]
\eqn
\bqn
+\frac{ph^2i}{4\pi}\cdot 2\pi i \cdot \frac{1}{y_2}\cdot[\frac{(1-y_2)y_1}{(y_1+y_2)h^2}+\frac{(1-y_1)y_2}{(y_1+y_2)h^2}+\frac{-1}{h^2}].
\eqn
Under these conditions,
\bqn
\eqref{limitL2}=ph^2[\frac{1}{h^2}-\frac{1}{y_1+y_2}],
\eqn

we can obtain the conclusion for the case where $y_1<1$ and $y_2>1$ in a similar manner.
According to the discussion above, the limit part is
{
 \begin{equation*}
\tilde l_n =\frac{y_{n_1}y_{n_2}}{y_{n_1}+y_{n_2}}+\frac{(1-y_{n_1})y_{n_2}}{y_{n_1}^2}\delta_{(y_{n_1}>1)}+\frac{y_{n_1}(1-y_{n_2})}
{y_{n_2}^2}\delta_{(y_{n_2}>1)}
\end{equation*}}
\end{proof}

\begin{proof} {of the mean part $\tilde{\mu}_{n}$ in Theorem \ref{th2}.}
According to the above discussion, we assume that $y_1>1, y_2>1.$
From the transformation $$m_3=-(1+hr\xi)/(1-y_2),$$ and its relationship with $m_3$, $z$ satisfies $$z=-\frac{m_3(m_3+1-y_1)}{(1-y_2)m_3+1}.$$

We then obtain
\bqn
\frac{z}{\alpha+z}=\frac{y_2|1+hr\xi|^2}{|y_2+hr\xi|^2}~~\mbox{and
}~~
1-\frac{z}{\alpha+z}=\frac{y_1(y_2-1)^2}{|y_2+hr\xi|^2}.
\eqn
In this mean part, we calculate the {integrand $f_{\gamma}(\frac{z}{\al+z})$,} which is equivalent to
\begin{equation}\label{meaninte}
c_1(\frac{1}{c_1}\cdot\frac{z}{\alpha+z}-1)^2+c_2(\frac{1}{c_2}\cdot(1-\frac{z}{\alpha+z})-1)^2.
\end{equation}
According to the above discussion,
\bqn
\eqref{meaninte}=c_1(\frac{1}{c_1}\cdot\frac{y_2|1+hr\xi|^2}{|y_2+hr\xi|^2}-1)^2+
c_2(\frac{1}{c_2}\cdot\frac{y_1(y_2-1)^2}{|y_2+hr\xi|^2}-1)^2
\eqn
\bqn
=c_1(\frac{|1+hr\xi|^2(y_1+y_2)}{|y_2+hr\xi|^2}-1)^2+c_2(\frac{(y_2-1)^2(y_1+y_2)}{|y_2+hr\xi|^2}-1)^2
\eqn
\bqn
=\frac{y_2(|1+hr\xi|^2(y_1+y_2)-|y_2+hr\xi|^2)^2}{(y_1+y_2)|y_2+hr\xi|^4}+
\frac{y_1((y_2-1)^2(y_1+y_2)-|y_2+hr\xi|^2)^2}{(y_1+y_2)|y_2+hr\xi|^4}
\eqn
\bqn
=\frac{y_1y_2(y_1+y_2)(h^2-y_2^2+2y_2+hr\xi+\frac{h}{r\xi})^2}{(y_1+y_2)|y_2+hr\xi|^4}
=\frac{y_1y_2(|1+hr\xi|^2-(y_2-1)^2)^2}{|y_2+hr\xi|^4}.
\eqn
Therefore, the mean part is
\bqn
{\eqref{mean1}+\eqref{mean2}+\eqref{mean3}=}\lim_{r\downarrow1}\frac 1{4\pi i}\oint_{|\xi|=1}\frac{y_1y_2[|1+hr\xi|^2-(y_2-1)^2]^2}{|y_2+hr\xi|^4}\cdot
(\frac1{\xi-\frac1r}+\frac1{\xi+\frac1r}-\frac2{\xi+\frac h{y_2r}})\mathrm{d}\xi
\eqn
\bqn
-\lim_{r\downarrow1}
\frac {\Delta_1}{2\pi i}\oint y_1\frac{y_1y_2[|1+hr\xi|^2-(y_2-1)^2]^2}{|y_2+hr\xi|^4}\cdot\frac{(1-y_2)^2h}{y_2^3}\cdot
\frac{\xi}{(\xi+\frac{h}{y_2r})^3}\mathrm{d}\xi
\eqn
\bqn
=\lim_{r\downarrow1}
\frac {\Delta_2}{4\pi i}\oint \frac{y_1y_2[|1+hr\xi|^2-(y_2-1)^2]^2}{|y_2+hr\xi|^4}\cdot\frac{(y_2-1)(\xi^2-\frac{h^2}{y_2r^2})}{
y_2(\xi+\frac{h}{y_2r})^2}\cdot[ \frac{2\xi}{(\xi^2-\frac{h^2}{y_2r^2})}-
\frac{2}{\xi+\frac{h}{y_2r}}]\mathrm{d}\xi.
\eqn
According to  Cauchy's residue theorem, we have
\bqn
{\tilde{\mu}_{n}=}\frac{y_1y_2h^2}{(y_1+y_2)^2}
+\Delta_1\frac{y_1^2y_2h^2(h^2+2y_2(y_2-y_1))}{(y_1+y_2)^4}
+\Delta_2\frac{y_2^2y_1h^2(h^2+2y_1(y_1-y_2))}{(y_1+y_2)^4},
\eqn which completes the proof of the mean part.
\end{proof}

\begin{proof} {of the variance part in Theorem \ref{th2}.}
 Under the case where $y_1>1$ and $y_2>1$, using the transformation discussed above, we can easily find that 
\bqn
\frac{\mathrm{d}m_3(z_1)\mathrm{d}m_3(z_2)}{(m_3(z_1)-m_3(z_2))^2}
=\frac{(1-y_2)^2}{(\frac{h}{r_1\xi_1}-\frac{h}{r_2\xi_2})^2}\cdot
\frac{h\mathrm{d}\xi_1}
{(1-y_2)r_1\xi_1^2}\cdot\frac{h\mathrm{d}\xi_2}{(1-y_2)r_2\xi_2^2}
=\frac{r_1r_2}{(r_1\xi_1-r_2\xi_2)^2}\mathrm{d}\xi_1\mathrm{d}\xi_2,
\eqn
\bqn
\frac{\mathrm{d}m_3(z_1)}{(m_3(z_1)+1)^2}
=\frac{(1-y_2)^2}{(y_2+\frac{h}{r_1\xi_1})^2}
\cdot\frac{h\mathrm{d}\xi_1}{(1-y_2)r_1\xi_1^2}
=\frac{(1-y_2)r_1h\mathrm{d}\xi_1}{(y_2r_1\xi_1+h)^2},
\eqn
\bqn
\frac{\mathrm{d}m_3(z_2)}{(m_3(z_2)+1)^2}
=\frac{(1-y_2)^2}{(y_2+\frac{h}{r_2\xi_2})^2}\cdot
\frac{h\mathrm{d}\xi_2}{(1-y_2)r_2\xi_2^2}
=\frac{(1-y_2)r_2h\mathrm{d}\xi_2}{(y_2r_2\xi_2+h)^2}.
\eqn

Thus,
\bqn
{\eqref{covariance1}+\eqref{covariance2}=2}\lim_{r_2\downarrow1}\oint\frac{1}{2\pi i}\cdot
\frac{y_1y_2[|1+h_2r_2\xi|^2-(y_2-1)^2]^2}{|y_2+h_2r_2\xi|^4}
\eqn
\bqn
\cdot\{\lim_{r_1\downarrow1}\oint\frac{1}{2\pi i}\cdot
\frac{y_1y_2[|1+h_1r_1\xi|^2-(y_2-1)^2]^2}{|y_2+h_1r_1\xi|^4}\cdot
\frac{r_1}{(r_1\xi_1-r_2\xi_2)^2}\mathrm{d}\xi_1\}\cdot r_2\mathrm{d}\xi_2
\eqn
\bqn
+c\{\lim_{r_1\downarrow1}\oint\frac{1}{2\pi i}\cdot
[\frac{y_1y_2[|1+h_2r_2\xi|^2-(y_2-1)^2]^2}{|y_2+h_2r_2\xi|^4}]\cdot
\frac{(1-y_2)r_1h\mathrm{d}\xi_1}{(y_2r_1\xi_1+h)^2}\}
\eqn
\bqn
\cdot\{\lim_{r_2\downarrow1}\oint\frac{1}{2\pi i}\cdot
[\frac{y_1y_2[|1+h_1r_1\xi|^2-(y_2-1)^2]^2}{|y_2+h_1r_1\xi|^4}]\cdot
\frac{(1-y_2)r_2h\mathrm{d}\xi_2}{(y_2r_2\xi_2+h)^2}\},
\eqn
where $c=y_1\Delta_1+y_2\Delta_2$ for the sake of brevity.
According to  Cauchy's residue theorem,
\begin{align*}
{\tilde{\nu}_n^2=}\frac{4y_1^2y_2^2h^2(h^2+2(y_1-y_2)^2)}{(y_1+y_2)^4}+c\frac{4y_1^2y_2^2h^4(y_1-y_2)^2}{(y_1+y_2)^6}, 
\end{align*}
which completes the proof.
\end{proof}
\section*{Acknowledgments}
The authors would like to thank the anonymous referee and associate editor for their invaluable and constructive comments. The research was partially supported by
NSFC (No. 11571067, 11771073) and Foundation of Jilin Educational Committee (No. JJKH20190288KJ).


\end{document}